\documentclass[12pt]{article}
\usepackage{localeng}
\usepackage[hyphens]{url}
\usepackage{hyperref}
\usepackage[hyphenbreaks]{breakurl}

\usepackage[a4paper,margin=20mm]{geometry}                % See geometry.pdf to learn the layout options. There are lots.

\theoremstyle{remark}

\newcommand{\bibquote}[1]{}

\makeatletter
\newcommand{\nb}[1]{{\color{red}}}
\newcommand*{\gl}{\nobreak\hskip1pt}
\DeclareRobustCommand*{\dash}{\gl\hbox{-}\gl}
\DeclareRobustCommand*{\endash}{\gl\hbox{--}\gl}
\DeclareRobustCommand*{\emdash}{\gl\hbox{---}\gl}
\makeatother

\let\eps=\varepsilon

\begin{document}

\title{Mathematical works of Vladimir A. Uspensky:\\ a commentary\footnote{Vladimir Andreevich Uspensky was my teacher (and undegraduate and Ph.D advisor). Here I concentrate on his mathematical works; I hope to express my deep gratitude to him elsewere.}}
\author{Alexander Shen\footnote{LIRMM, University of Montpellier, CNRS, Montpellier, France. Supported by RaCAF--ANR-15-CE40-0016.}}
\date{}
\maketitle

\begin{abstract}
Vladimir Andreevich Uspensky [1930--2018] was one of the Soviet pioneers of the theory of computation and mathematical logic in general. This paper is the survey of his mathematical works and their influence. (His achievements in linguistics and his organizational role are outside the scope of this survey.) 

\end{abstract}

\subsection*{Harmonic functions}

The first paper of Uspensky~\cite{1949} appeared when he was an undergraduate student. It suggests an elementary approach to harmonic functions that is based on the definition of a harmonic function on $\mathbb{R}^2$ as a function that has the mean value property. The main tool in the following observation: for fixed two poinst $A$, $B$ the oriented angle $ACB$ is a harmonic function of $C$,\footnote{This is easy to see, because both direction angles $CA$ and $CB$ satisfy the mean value property and the same is true for their difference.} and this function is a locally constant function on any circle that goes through $A$ and $B$.

\medskip

\subsection*{Master thesis}

Uspensky's thesis advisor was Andrei Kolmogorov; the thesis~\cite{1952} was written and defended in 1952. It contains the description of the model of computation suggested by Kolmogorov and now known as \emph{Kolmogorov--Uspensky machines}. It is shown that this model is equivalent to partial recursive functions (defined in terms of substitution, recursion and $\mu$-operator). Moreover, this model is used to define relative computability with respect to some function $f$. For that, the graph of $f$ is represented as an infinite graph (a \emph{complex}) that is available to the graph transformation algorithm together with the input [definition (A) on p.~64]. This definition is compared with other definition of relative computability. For that, Uspensky reformulates the Turing--Post definition~\cite{Turing1939,Post1944}, see  definition (T) on p.~63, and shows that definitions (T) and (A) are equivalent. Moreover, Uspensky proves the equivalence with one more definition given in terms of the closure of the function $f$ (and basic functions) under substitution, recursion and $\mu$-operator. 

The thesis advisor (Kolmogorov) writes in his opinion:

\begin{quote}

This paper analyses (in more details than it was done before) the very notion of algorithmic  computability.

(1) The author reproduces the only completely formal definition of algorithmic reducibility (p.~22) that he ascribes to Boris Trakhtenbrot: a function $\gamma$ is reducible to another function $\delta$ if $\gamma$ belongs to the recursive closure of $\delta$. The author shows that in fact such a reduction can be performed \emph{in some simple canonical way, using some fixed primitive recursive functions $\tau(u)$ and $\omega(u)$ and some primitive recursive functions $h(u,v,w)$ and $\varphi(m)$ that depend on $\gamma$ and $\delta$; see theorem on p.~28. This is the main result of the paper from the purely mathematical viewpoint.}

Trakhtenbrot's definition, on the other hand, needs some ``justification'' that shows that it corresponds to the intuitive notion of reduction: there is a ``mechanical'' way of obtaining $\gamma(x)$ \emph{assuming} that the values $\delta(x)$ are somehow made ``accessible'' for each $x$. The general framework for such a justification were given by Post; the translation of the corresponding part of Post's paper~\cite{Post1944} is reproduces in the thesis. \emph{Then the author gives a completely formal definition of reduction that corresponds to this idea} (probably for the first time) and proves its equivalence to Trakhtenbrot's definition. This is also a very significant achievement of the author.

Also the paper contains a good survey of the definitions of algorithmic computability for function $y=\gamma(x)$ with numerical arguments and values. It is centered around the definition suggest by myself. The thesis provides the motivation for this definition and proves that it is equivalent to the previous ones. \emph{In a sense this equivalence can be considered as a ``justification'' for the previous definitions since my definition makes especially clear the idea of an algorithmic computability; the algorithms that are used are models of the real computation devices, only the amount of ``memory'' is assumed to be unbounded}. 

\end{quote}

To understand the value of this paper, one should recall the historical context (now almost forgotten). Let us make few historical comments.

\subsubsection*{Partial recursive functions}

Ask an expert what is a partial recursive function. Most probably the answer would be: this is a function that can be obtained from basic functions (projection function, zero constant and successor function) using substitution, recursion and minimization ($\mu$-operator). This definition can be found in the classical text of Odifreddi~\cite[p.~127]{Odifreddi1989} and other sources (see, e.g., \cite{Malcev1965,Cutland1980} or Wikipedia article~\cite{Wiki2018}).

However, the original definition was different. The traces of this old definition can be found in another classical textbook~\cite[Section 1.5]{Rogers1972} and in Wolfram MathWorld~\cite{Wolfram2018} site. This definition in equivalent (gives the same class of partial function), but it is different, and one should have this difference in mind while reading old papers.

Let us try to clarify the history. Recursive definitions were well known for a long time (recall the Fibonacci sequence). They were systematically used to define arithmetic functions in the paper of Skolem~\cite{Skolem1923}. He realized that in this way one can define not only addition or multiplication (by the recurrent formula like $x+y' = (x+y)'$ or $x\cdot y' = x\cdot y + x$, where $x'$ is the successor of $x$), but also many other functions that appear in the elementary number theory. After these definitions are given, one could prove basic facts of number theory by induction.\footnote{The goal was to show that many mathematical results can be proven by a simple and robust way, just by using recursive definitions and inductive arguments. Now the corresponding theory is known as primitive recursive arithmetic.}

Skolem did not explicitly consider the class of functions that can be defined recursively in this way. However, in his 1925 talk Hilbert~\cite{Hilbert1926} says that ``the elementary means that we have at our disposal for forming functions are \emph{substitution} (that is, replacement of an argument by a new variable or function) and \emph{recursion} (according to the schema of the derivation of the function value for $n+1$ from that for $n$)''. He then considers the sequence of functions
\begin{multline*}
\varphi_1(a,b)=a+b, \ \varphi_2(a,b)=a\cdot b,  \ \varphi_3(a,b)=a^b,\\
\varphi_4(a,b)= \text{[$b$th term in the sequence $a, a^a, a^{(a^a)}, a^{a^{(a^a)}}\ldots$]}
\end{multline*}
that can be defined in general by the recurrent formula
\[
\varphi_1(a,b)=a+b, \ \varphi_{n+1}(a,1)=a,  \ \varphi_{n+1}(a,b+1)=\varphi_n(a,\varphi_{n+1}(a,b)),
\]
and mentions Ackermann's result saying that the function  $\varphi_n(a,b)$ ``cannot be defined by substitutions and ordinary, step-wise recursions'' (this result was later published in~\cite{Ackermann1928}). When stating this negative result, Hilbert implicitly considers the class of function that can be defined by ``substitutions and ordinary, stepwise recursions'' (even though this class is not defined explicitly and there is no name for the functions from this class.)

Such a definition (and name) appeared in the classical work of G\"{o}del~\cite[p.~179]{Godel1931}: a function is called recursive (\emph{rekursiv} in German) if it can be obtained by a sequence of substitution and recursion operation (we construct $\varphi$ assuming that $\psi$ and $\mu$ are already constructed): 
\begin{align*}
\varphi(0,x_2,\ldots,x_n)&=\psi(x_2,\ldots,x_n)\\
\varphi(k+1,x_2,\ldots,x_n)&=\mu(k,\varphi(k,x_2,\ldots,x_n),x_2,\ldots,x_n)
\end{align*}
(scheme (2) p.~179). G\"{o}del proves that these functions could be represented in a formal system, so for him this class of functions is more a tool than an object.

G\"{o}del's definition does not cover more general recursive definition (like the one used by Ackermann). How can they be treated? Herbrand (in a letter to G\"{o}del and in~\cite{Herbrand1932}) suggested that one could consider systems of functional equations that related the functions we define with the already defined ones. He writes in~\cite[p.~5, p.~624 of the English translation]{Herbrand1932}:

\begin{quote}
We may also introduce any number of functions $f_i(x_1,\ldots,x_{n_i})$ together with hypotheses such that
\begin{description}
\item{(a)} The hypotheses contain no apparent variables;
\item{(b)} Considered intuitionistically,\footnote{This expression means: whey they are translated into ordinary language, considered as a property of integers and not as a mere symbol. [Herbrand's footnote]} they make the actual computation of $f_i(x_1,\ldots,x_{n_i})$ of the $f_i(x_1,\ldots,x_{n_i})$ possible for every given set of integers, and it is possible to prove intuitionistically that we obtain a well-determined result.
\end{description}
\end{quote} 

The reference to intuitionism sounds a bit unclear; probably it means that it is not enough to have a functional equation or a system of equations for which we can prove (using arbitrarily powerful tools) that it has a unique solution. We require that the proof is constructive and provides a method to compute the values of the functions starting from the equations. Indeed, later Kalmar~\cite{Kalmar1955} gave an example of a system of functional equations that uniquely defines a non-computable function.
 
G\"{o}del returns to Herbrand's suggestion (Herbrand died in mountains just after sending his paper~\cite{Herbrand1932} to the editors) in his Princeton's lectures. (The lecture notes circulated at that time and later were reprinted, see~\cite{Herbrand1932}.) As before, he considers <<recursive functions>> that can be obtained from basic functions by substitutions and ``ordinary'' recursions; however, in Section~9 he mentions the recursive definitions of more general type. The functions defined in this way are called ``general recursive functions''. He says:

\begin{quote}
One may attempt to define this notion [general recursive function] as follows: if $\phi$ denotes an unknown function, and $\psi_1,\ldots,\psi_k$ are known functions, and if the $\psi$'s and the $\phi$ are substituted in one another in the most general fashions and certain parts of the resulting expressions are equated, then if the resulting set of functional equations has one and only one solution for $\phi$, $\phi$ is a recursive function.'' 
\end{quote}
(and mentions Herbrand's letter as a reference). Then he added some restrictions that clarify Herbrand's idea:
\begin{quote}
We shall make two restrictions on Herbrand's definition. The first is that the left-hand side of each of the given functional equations defining $\phi$ shall be of the form \[\phi(\psi_{i1}(x_1,\ldots,x_n),\psi_{i2}(x_1,\ldots,x_n),\ldots,\psi_{il}(x_1,\ldots,x_n)).\] The second (as stated below) is equivalent to the condition that all possible sets of arguments $(n_1,\ldots,n_l)$ of $\phi$ can be so arranged that the computation of the value of $\phi$ for any given set of arguments $(n_1,\ldots,n_l)$ by means of the given equations requires a knowledge of the values of $\phi$ only for sets of arguments which precede $(n_1,\ldots,n_l)$.
\end{quote}
G\"{o}del does not specify the ordering on the tuples (used as arguments), so the exact meaning of this definition is unclear. But later he specifies the derivation rules that allow to derive an equality from the other ones, and says:
  \begin{quote}
 Now our second restriction on Herbrand's definition of recursive function is that for each set of natural numbers $k_1,\ldots,k_l$ there should be one and only one $m$ such that $\phi(k_1,\ldots,k_l)=m$ is a derived equation.
\end{quote}
In this way G\"{o}del gives a quite formal definition of some class of functions called ``general recursive functions'' (usually translated to Russian as \rus{<<общерекурсивные функции>>}. However, as Kleene explains in~\cite{Kleene1981}, at the time of these lectures (1934) G\"{o}del was not sure that this class of functions is general enough: <<However, G\"{o}del, according to a letter he wrote to Martin Davis on 15 February 1965, ``was, at the time of [his 1934] lectures, not at all convinced that [this] concept of recursion comprises all possible recursions''>>~\cite[p.~48]{Kleene1981}. Davis writes in~\cite[p.~40]{Davis1965}:
\begin{quote}
In the present article [Davis discussed~\cite{Godel1934}] G\"odel shows how an idea of Herbrand's can be modified so as to give a general notion of recursive function $\langle\ldots\rangle$ G\"odel indicates (cf. footnote 3) that he believed that the class of functions obtainable by recursion of the most general kind were the same as those computable by a finite procedure. However, Dr.~G\"{o}del has stated in a letter that he was, at the time of these lectures, not at all convinced that his concept of recursion comprised all possible recursions; and that in fact the equivalence between his definition and Kleene's in Math. Ann.~112~[this is~\cite{Kleene1936} in our list] is not quite trivial. So despite appearances to the contrary, footnote 3 of these lectures is not a statement of Church's thesis.
\end{quote}
Footnote 3~\cite[p.~44]{Davis1965} discusses the claim that every primitive recursive function (obtained by substitutions and ``ordinary recursions, see below) can be computed by a finitary process, and says that ``The converse seems to be true, if, besides recursions according to the scheme (2) [primitive recursion], recursions of other forms (e.g., with respect to two variables simultaneously) are admitted. This cannot be proved, since the notion of finite computation is not defined, but it serves as a heuristic principle''.

R\'{o}sza P\'{e}ter in~\cite{Peter1934} studies the ``ordinary recursions'' and proves, for example, that one may use several values of the function (for smaller arguments) in the recursive definition and still get the same class of functions. She introduces the name ``primitive Rekursion'' for the ``ordinary'' recursions considered by her predecessors.

Then Kleene in~\cite{Kleene1936}  (1936) introduces the name ``primitive recursive functions'' (\rus{<<примитивно рекурсивные функции>>} in Russian) for functions that can be obtained by substitutions and primitive (=``ordinary'') recursion. At the same time, Kleene suggests to consider a bigger class of functions. He calls the functions from this class ``general recursive function'' (the title of his paper is \emph{General recursive functions of natural numbers}). This class is defined following Herbrand and G\"{o}del; Kleene considered different versions of derivation rules for equalities and shows that they lead to the same class of functions.

Kleene also introduces ``$\eps$-operator''. Namely, $\eps x [A(x)]$ is defined as the minimal $x$ such that $A(x)$ or $0$ if such an $x$ does not exists. This operator is used in Theorem~IV that says that every general recursive function can be represented as
\(
\psi(\eps y [R(x,y)]),
\)
for some primitive recursive function $\psi$ and some primitive recursive predicate~$R$ (this means that $R$ can be represented as $r=0$ for some primitive recursive function $r$), such that for every $x$ there exists $y$ such that $R(x,y)$.\footnote{Note that the clause in the definition of $\eps$-operator that lets the value to be $0$ when $x$ does not exists, is not used in Theorem~IV; so one can use the standard $\mu$-operator instead. (For $\mu$-operator the value is undefined if $y$ does not exst.)}  The next Theorem V says that the reverse statement is also true: every function that can be presented in this way is a general recursive function (in the sense of Herbrand and G\"{o}del). There this representation can be considered as an equivalent definition of the class of general recursive functions. Moreover, this definitions can be used to provide some numbering of all general recursive functions if we add an additional argument $e$ to $R$; not all values of $e$ lead to total functions. One could say that it this way we get a numbering of a family of partial functions, but in this paper Kleene does not considers this class (later called ``partial recursive functions'').

Church (also in 1936) publishes his paper~\cite{Church1936} where he defines some other class of functions with natural arguments and values in terms of some calculus (called $\lambda$-calculus) and claims that this class captures the intuitive idea of computability:
\begin{quote}
The purpose of the present paper is to propose a definition of effective calculability${}^3$ which is thought to correspond satisfactorily to the somewhat intuitive notion.
\end{quote}
Here $({}^3)$ is Church's footnote: 
\begin{quote}
As will appear, this definition of effective calculability can be stated in either of two equivalent forms, (1) that a function of positive integers shall be called effectively calculable if it is $\lambda$-definable in the sense of \S2 below, (2) that a function of positive integers shall be called effectively calculable if it is recursive in the sense of \S4 below. The notion of $\lambda$-definability is due jointly to the present author and S.C.~Kleene $\langle\ldots\rangle$ The notion of recursiveness in the sense of \S4 is due jointly to Jacques Herbrand and Kurt G\"odel $\langle\ldots\rangle$ The proposal to identify these notions with the intuitive notion of effective calculability is first made in the present paper\ldots
\end{quote}
Church adds (a footnote in \S7): 
\begin{quote}
The question of the relationship between effective calculability and recursiveness (which it is here proposed to answer by identifying the two notions) was raised by G\"{o}del in conversation with the author. The corresponding question of the relationship between effective calculability and $\lambda$-definability had previously been proposed by the author independently.
\end{quote}

It is clear from this footnote that for Church the suggestion to identify the intuitive notion of effective calculability with the formally defined class of functions (for which two equivalent definitions are given) is an important contribution. This suggested became known as \emph{Church's thesis}.

Almost at the same time Turing publishes his paper~\cite{Church1936} where he defines the model of computation now called \emph{Turing machines}. Turing calls them $a$-machines (`a' for `automatic'). Turing also constructs the universal machine that can simulate any Turing machine when equipped by a suitable problems. Turing uses this type of machines to define the notion of a computable real number (the digits in the positional representation can be computer by a machine), and also gives his proof of the undecidability of the Entscheidungsproblem (there is no algorithm that can tell whether a given first order formula is logically valid, i.e., true in all the interpretations of the language). Earlier similar results (for equivalent definitions of computability) were proven by G\"{o}del and Kleene, as well as Church (see~\cite[p.109]{Davis1965} for details).

In an Appendix (added August 28, 1936) Turing sketches the proof of equivalence between two definitions of computability of a sequence: in terms of $a$-machines and in terms of $\lambda$-calculus. Describing this result in the Introduction, he writes:
\begin{quote}
In a recent paper Alonzo Church has introduced an idea of ``effective calculability'', which is equivalent to my ``computability'', but is very differently defined. Church also reaches similar conclusions about the Entscheidungsproblem. The proof of equivalence between ``computability'' and ``effective calculability'' [i.e., $\lambda$-definability] is outlined in an appendix to the present paper.
\end{quote}

Independently of Turing (and almost simultaneously) Post publishes his paper~\cite{Post1936}, where he introduces the notion of a ``finite combinatory process'' that is very similar to Turing machines. Some technical details are different; one should mention also that Post never speaks about a machine. He describes how a ``problem solver or worker'' follow ``the set of directions'' of a fixed type. Then Post writes: 
\begin{quote}
The writer expects the present formulation to turn out to be logically equivalent to recursiveness in the sense of the G\"{o}del\endash Church development. Its purpose, however, is not only to present a system of a certain logical potency but also, in its restricted field, of psychological fidelity. In the latter sense wider and wider formulations are contemplated. On the other hand, our aim will be to show that all such are logically equivalent to formulation 1 [the definition suggest by Post]. We offer this conclusion at the present moment as a \emph{working hypothesis}. And to our mind such is Church's identification of effective calculability with recursiveness. $\langle\ldots\rangle$ The success of the above program would, for us, change this hypothesis not so much to a definition or to an axiom but to a \emph{natural law}. 
\end{quote}
In a footnote Post adds:
\begin{quote}
Actually the work already done by Church and others carries this identification considerably beyond the working hypothesis stage. But to mask this identification under a definition hides the fact that a fundamental discovery in the limitations of the mathematizing power of Homo Sapiens has been made and blinds us to the need of its continual verification.\footnote{%
Probably now this ``fundamental discovery'' has lost its value and even may be its meaning: when speaking about the equivalence between the intuitive idea of algorithmic computability and a formal definition, we assume that this intuitive idea was developed independently of any model of computation or programming language. But now it would be almost impossible to find anyone who learned the intuitive notion of algorithm before having some programming experience.}
\end{quote}

It is clear that in 1936 the puzzle (as we know it now) was almost completely assembled: there are several definitions of computability that are shown to be equivalent (the classes of computable functions are the same); these definition are considered as reflecting the intuitive notion of an algorithm, and there are some intuitive arguments that support this thesis.

However, there are two points where the picture is different from the modern one. The first is more about terminology: none of the papers that define recursive functions defines this class using substitution, recursion and $\mu$-operator though all the tools to prove the equivalence are ready and this equivalence is mentioned explicitly by Kleene in 1943~\cite[p.~53, Corollary]{Kleene1943}.

Second, more important difference is that all these papers consider only \emph{total} functions (defined for all natural arguments). \emph{Partial} functions appear only later, in Kleene's paper~\cite{Kleene1938} (published in 1938) where the computable notation systems for ordinal are considered (and partial computable functions are essential). Kleene describes the process of derivation in the sense of Herbrand and G\"{o}del and assumes that such a derivation exists only for one function value (for given arguments). Then he writes:
\begin{quote}
If we omit the requirement that the computation process always terminate, we obtain a more general class of functions, each function of which is defined over a subset (possibly null or total) of the $n$-tuples of natural numbers, and possesses the property of effectiveness when defined. These functions we call partial recursive.
\end{quote}
 In this way the notion of a \emph{partial recursive function} is introduced.\footnote{The traditional Russian translation of this name is \rus{частично рекурсивная функция}. It sound even more strange than \rus{общерекурсивная функция} for general recursive functions; one could think that the function is not completely recursive but only partially recursive.} Kleene considers substitutions and recursions (that can be naturally extended to partial functions), and then defines $\mu$-operator for partial functions:
\[
\mu y [R(m,y)=0] = n
\]
for a partial function $R$ if $R(m,n)$ is defined and equals $0$ while all previous values $R(m,0)$,\ldots, $R(m,n-1)$ are defined and are not zeros. It is obvious that $n$ with this property is unique; however, it may not exist, and in this case the $\mu$-operator defines a non-total function (that is undefined on $m$). Kleene notes that the class of partial recursive functions defined in the language of Herbrand and G\"{o}del is closed under all three operations (substitution, recursion and $\mu$-operator). He notes also that for every $n$ there exists a universal function $\Phi_n(z,\mathbf{x})$ of $n+1$ variables such that every partial recursive function of $n$ variables $\mathbf{x}$ can be obtained from $\Phi_n$ by fixing some value of the first argument~$z$. This universal function $\Phi_n$ can be represented as
  \[
\Phi_n (z,\mathbf{x})= S(z,\mu y T_n(z,\mathbf{x},y)),
  \]
where $S$ is some primitive recursive function and $T_n$ is a primitive recursive predicate (saying that some primitive recursive function equals $0$). Informally speaking, $z$ is a natural number that encodes a system of functional equations (in Herbrand -- G\"odel style) that defines some partial recursive  function of $n$ variables, and $y$ is an encoding of a derivation that, starting with these equations, establishes the value of this partial recursive function on~$\mathbf{x}$. The predicate $T_n$ checks the correctness of this derivation, and the function $S$ extracts the function value from it.\footnote{Kleene provides $z$ as the first argument to the function $S$ but this is not necessary.} This result is called ``Kleene's normal form theorem''; it implies that partial recursive function could be equivalently defined as functions that can be obtained by substitution, recursion, and $\mu$-operator. One may also require additionally that the $\mu$-operator is used only once (being applied to a primitive recursive functions). However, this way of defining partial recursive functions is not mentioned by Kleene.

The same framework and terminology is used in a later paper of Kleene~\cite{Kleene1943} (1943, where he consider the arithmetical hierarchy) and in his classical book of 1952~\cite{Kleene1957} that remained a standard reference for logic and computability theory for a long time. Let us mention again a detail that may sound strange in our time: the statement of ``Church's thesis'' (the equivalence between the intuitive notion of computability and formal definitions) mentions only total functions.

\subsubsection*{Relative (oracle) computability}

One can define the notion of computability of a function relative to some other function (or set, if we identify sets with their characteristic functions). This definition was first considered in Turing's Ph.D thesis (1939, see~\cite{Turing1939}); however, it was only a side remark and only reducibility to some specific set was considered. Turing writes:
\begin{quote}
Let us suppose that we supplied with some unspecified means of solving number-theoretic problems; a kind of oracle as it were. We will not go any further into the nature of this oracle than to say that it cannot be a machine.  With the help of the oracle we could form a new kind of machine (call them $o$-machines), having as one of its fundamental processes that of solving a given number-theoretic problem. More definitely these machines are to behave in this way. The moves of the machine are determined as usual by a table except in the case of moves from a certain internal configuration $\mathfrak{o}$. If the machine is in the internal configuration $\mathfrak{o}$ and if the sequence of symbols marked with $l$ is then the well formed formula \textbf{A}, then the machine goes into the internal $\mathfrak{p}$ or $\mathfrak{t}$ according as it is or is not true that \textbf{A} is dual. The decision as to which is the case is referred to the oracle.\par These machines may be described by tables of the same kind as used for the description of $a$-machines, there being no entries, however, for the internal configuration $\mathfrak{o}$.
\end{quote}
The definition of Turing reducibility for the general case was given by Post in his famous article~\cite[Section 11]{Post1944} where he formulated \emph{Post's problem} (asking whether there exists a recursive enumerable non-recursive set $X$ that it is not Turing-complete: not all recursively enumerable sets are reducible to~$X$). Formally speaking, Post considers the case when both sets (the one being reduced and the other to which it is reduced) and recursively enumerable, but the definition is the same for the general case of arbitrary sets of natural numbers. The Post's definition follows the scheme sketched by Turing. Kleene in 1943~\cite{Kleene1943} suggests a different approach: we define general recursive functions using Herbrand -- G\"{o}del derivations but extend the list of ``axioms'' adding the full information about the values of some fixed total functions $\psi_1,\ldots,\psi_k$. The functions that are definable in this way are then called \emph{general recursive functions in $\psi_1,\ldots,\psi_k$}:
\begin{quote}
A function $\phi$ which can be defined from given functions $\psi_1,\ldots,\psi_k$ by a series of applications of general recursive schemata we call \emph{general recursive} in the given functions; and in particular, a function $\phi$ definable ab initio by these means we call \emph{general recursive}.
\end{quote}
However, Kleene does not develop this idea (which remains a side remark), and does not define relative computability for the case of partial functions (only total functions are considered). In 1952 book Kleene extends the definitions to partial functions and proves that the resulting definition (in Herbrand -- G\"{o}del style) is equivalent to the definition of relative computability given by Turing and Post~\cite[\S 69]{Kleene1957}. The oracle is assumed to a be total function (or a tuple of total functions) but no other restrictions are imposed; recall that Post considered only recursively enumerable sets as oracles.

A survey of different definitions of relative computability can be found in~\cite[Section 4.3, ``History of Relative Computability'']{Soare1996}.

Now we can explain what was the Uspensky's contribution in his master thesis~\cite{1952}:\footnote{Unfortunately this paper was not published, though both reports (by Kolmogorov, the thesis advisor, and by Petr Sergeevich Novikov, the reviewer) recommended its publication. So --- alas --- it hardly could play any role in the further developments.}

\begin{itemize}

\item  For the first time, the (now standard) definition of partial recursive functions in terms of substitutions, recursions, and $\mu$-operator was stated explicitly (with a reference to an ``idea of Boris Trakhtenbrot''~\cite[p.~22]{1952}).

\item It was shown (simultaneously with~\cite[\S 69]{Kleene1957} and in much more clear way) that this definition is equivalent to other definitions of (absolute and relative) computability.

\item For the first time, a ``machine-independent'' definition of relative computability was given. Here machine independence means that the definition does not use any model of computation but only the class of computable functions. It was shown that this definition is equivalent to other definitions of relative computability.

\item Finally, it was the first paper that presents the model of computation suggested by Kolmogorov (later it was published in a joint paper by Kolmogorov and Uspensky~\cite{1958}), the definition of relative computability in terms of this model, and the proof of equivalence of this definition to other definitions of relative computability.
\end{itemize}

The third item in this list requires some clarifications. The Turing -- Post definition of relative computability is a modification of the corresponding definition for (absolute) computability: we extend the class of Turing machines by allowing them to get ``answers'' from an oracle. Similarly, the Kleene's definition of the relative computability modifies the definition of the computable (partial recursive) functions. So even if we have already agreed on the definition of (absolute) computability, we still may not left this definition behind when defining relative computability. Instead, in the latter definition we need to return to the model of computation and make some modifications (that allow some kind of ``oracle access'').

On the other hand, Uspensky defines relative computability in terms of a dialog with an oracle, and this dialog should be computable in the sense that some (partial) functions that describe this dialog should be computable. These function should describe the dialog in the following sense: they specify the next question to the oracle (or output if no more questions are needed) given the input and the list of previous questions and oracle answers.

Now the ``machine-independent'' definitions of relative computability are quite standard. For example, one of them can be found in the classic textbook of Rogers~\cite[Section 9.2]{Rogers1972} (without any references to previous work). One can also note that Uspensky's definition has a technical advantage: unlike the definition from~\cite{Rogers1972} it can be naturally generalized to a partial oracles $\psi$, and the class of functions that are obtained in this way is equal to the closure of the partial recursive functions and $\psi$ with respect to substitutions, recursions and $\mu$-operator.
%\nb{[Who proved this? I asked Slaman but he haven't yet answered. Probably this is mentioned in Odifreddi?]} 
However, Uspensky did not consider this generalization and always assumes that oracle is a total function (though the proof could be easily adapted to the case of partial oracles).

\subsection*{G\"{o}del's incompleteness theory and theory of computability}

The G\"{o}del incompleteness theorem and the class of recursive functions appeared not only at the same time but also together like Siamese twins. The classical paper where G\"{o}del proved incompleteness of Principia Mathematica and related systems~\cite{Godel1931} also introduced the notion of a recursive function (a primitive recursive function in modern terminology, see above), and this notion played an important technical role in the proof. Namely, several functions related to the encoding of formulas and proofs by natural numbers (their ``G\"{odel} numbers'') were defined recursively, and this definition was used to embed these notion into the formal system (thus making self-referential statements and formal reasoning about proofs possible).

On the other hand, the first definition of general recursive functions was given in terms of a formal system (calculus of equalities) that goes back to Herbrand and G\"{o}del.

One could that the separation of these Siamese twins was an important achievement both in the theory of computation and in the proof theory. And historically it was not so simple as it may seem now. The first step was done by Turing and Post that suggested models of computation that do not refer to any calculus (formal theory).  And then the general nature of G\"{o}del's incompleteness theorem was realized; this was done in 1940s by Kleene and (later, but independently) Kolmogorov.

In 1943 Kleene noted~\cite{Kleene1943} that G\"odel's theorem essentially claims that the set of true formulas is not recursively enumerable.\footnote{Now people say ``computably enumerable'' instead of ``recursively enumerable''. Since we do not consider other type of enumerations, we call these sets \emph{enumerable} in the sequel.} In 1950 he gave~\cite{Kleene1950} a similar interpretation for the Rosser's version of incompleteness theorem: it corresponds to the existence of two inseparable enumerable sets. So all the crucial observations were made by Kleene before 1950. Still the exposition both in this 1950 paper and in the 1952 textbook~\cite{Kleene1957} is intertwined with the language of primitive recursive function (it is enough to say that the exposition in~\cite{Kleene1950} starts by  ``Let $T_1$ be the primitive recursive predicate so designated in a previous paper by the author''), and the embedding of the inseparable sets into a formal theory is not described explicitly.

Shortly after than (but most probably, independently) Kolmogorov also realized the connection between G\"{o}del's incompleteness theorem and theory of algorithms. As Uspensky writes in ~\cite[p.~323]{2006a},
 
\begin{quote}

At December 2, 1952 Kolmogorov explained me main ideas relating G\"{o}del's incompleteness theorem for general calculi to the existence of [enumerable] sets that are not recursive, and pairs of [enumerable] sets that can not be separated by a recursive set. The explanation was quite concise (maybe, five minutes) but then he gave me a short written note entitled ``G\"{o}del and recursive enumerability'', so I could read and copy it. The note was written just for himself, and it was not easy for me to understand both the note and his oral comments. Then it became more clear, and on May 8, 1953 Kolmogorov submitted my short paper ``G\"{o}del's theorem and the theory of algorithms'' to Soviet Math. Doklady. When Kolmogorov worked with his students, he made them feel that they are the authors (and he became a coauthor of his students much more rarely than he deserved it) $\langle\ldots\rangle$ a paper ``On the definition of an algorithm'' was published in \rus{\emph{Успехи математических наук}}; in this paper my role was essentially technical.
\end{quote}
Here Uspensky speaks about two papers~\cite{1953,1958}. The second paper contains the detailed exposition of a model of computation based on graph transformations that appeared already in Uspensky's master thesis~\cite{1952} and is known as Kolmogorov -- Uspensky machines (see above). The first paper~\cite{1953} explains (without any reference to primitive recursive functions) that G\"{o}del's incompleteness theory (formal arithmetic is incomplete and cannot be completed) is a corollary of two facts: (1)~there exist recursively inseparable enumerable sets; (2)~this pair of inseparable sets can be embedded into the formal arithmetic (in modern language, can be $m$-reduced to the pair (provable formulas, refutable formulas). Moreover, for every enumerable set of additional axioms (that keeps the theory consistent) one can effectively point out a formula that is is neither provable nor refutable in this extended system, and this fact is a corollary of the existence of two \emph{effectively} inseparable sets.

Let me stress again that all these observations were  made already by Kleene in~\cite{Kleene1950}; it seems that Kolmogorov and Uspensky did not see that paper at the time.  Uspensky's paper~\cite{1953} has a reference to Kleene's 1943 paper~\cite{Kleene1943}; however, when speaking about inseparable enumerable sets, Uspensky does not refer to Kleene's 1950 paper~\cite{Kleene1950} where they were constructed and notes only that they were constructed by Novikov (and provides a reference to Trakhtenbrot's paper of 1953).

Generally speaking, there are two complementary views on G\"{o}del's theorem. The original G\"{o}del's argument is a version of the liar's paradox. This self-referential paradox notes that the statement ``This statement is false'' cannot be either true or false. If we consider instead the statement ``This statement is not provable'' (which can be, unlike the previous one, formulated in the language of arithmetic), we get a statement that is true and (therefore) not provable --- or false and provable, but we assume that formal arithmetic is consistent. This reasoning does not rely on the theory of algorithms; however, to show that one can translate finitary arguments into the language of formal arithmetic one can use primitive recursive functions as a technical tool (following G\"{o}del).

On the other hand, G\"{o}del incompleteness theorem is a consequence of the existence of an enumerable undecidable set (or, in a more symmetric version, of the existence of two recursively inseparable enumerable sets). In this way self-referential nature of the argument is hidden. But it is just moved to the proof of the existence of an undecidable enumerable set (or an inseparable pair). Indeed, this proof uses ``diagonal argument'' that goes back to Cantor, and this diagonal argument is of self-referential nature (the ``diagonal'' function appear when we apply a function to its own number, or run a program on its own text).

Much later Uspensky published a popular exposition of the proof of G\"{o}del's theorem based on the algorithms theory (together with the introduction to this theory) in~\cite{1974}. The extended version of this paper was published as a brochure~\cite{1982} (in the series ``Popular lectures in mathematics'' published by Nauka publishing house in Moscow). This work is probably the most accessible (and correct) non-technical exposition of G\"{o}del's incompleteness theorem in Russian literature (at least if we consider its algorithmic side).

In addition to that, these publications~\cite{1974,1982} suggest a way to explain theory of algorithms that was quite unusual at the time (one may compare them to Rogers' textbook~\cite{Rogers1972}). Usually the exposition started with a detailed analysis of some specific model of computation. The choice of this model changed with time. Initially most expositions used partial recursive functions; then Turing machine became the preferred model. In Russia Markov and his school preferred the so-called \emph{normal algorithms}. The analysis of this model required a lot of efforts (and space). Only after that the readers can learn the basic facts like Post's theorem (an enumerable set with enumerable complement is decidable), etc. Of course, the impatient reader could skip the boring first part, but then all the considerations in the rest of the textbook became baseless.

What can be done? Uspensky suggested the following approach used in~\cite{1974} (and before in his 1972/73 lectures, and may be even earlier). We consider the class of computable function assuming that this class satisfies some properties (``axioms''). These properties include the following ones:
\begin{itemize}
\item some specific functions (e.g., the pair numbering functions) are computable; some specific constructions (e.g., the conditional execution or loops) preserve computability;

\item \emph{The tracing axiom}: for every algorithm $A$ there exist a decidable set $R$ whose elements are called ``traces'', and two computable functions $\alpha$ and $\omega$. Informally, elements of $R$ are traces of terminating runs of $A$ on all possible inputs (that include all information about the computation); this set should be decidable since one can check that the trace is indeed a trace of $A$. The function $\alpha$ recovers the input from the trace; the function $\omega$ recovers the output. This is an informal explanation why this axiom is plausible; the formal requirement is only that $A(x)=y$ if and only if there exists $r\in R$ such that $\alpha(r)=x$ and $\omega(r)=y$.

\item \emph{The program axiom}: there exists a decidable set $P$ (whose elements are called ``programs'' and an algorithm $U$ that can be used to apply an arbitrary program $p\in P$ to arbitrary input~$x$ (so the input of $U$ is a pair $\langle p,x\rangle$). The axiom requires that every computable function $f$ has some program $p$ such that $U(p,x)=f(x)$ for every $x$. The last equality sign is understood as follows: either both sides are undefined or both sides are defined and equal.
\end{itemize}

After we agree with these axioms, we can prove results about computability without going into the technical details. On the other hand, it is quite clear what is missing in this picture to get a formally sound mathematical theory:
\begin{itemize}
\item We need to choose some model of computation.
\item We need to be able to program (in this model) some constructions used in the proofs. In fact, they could be not so simple (recall the priority arguments, for example).
\item We need to prove the tracing axiom and the program axiom for this model.
\end{itemize} 

This looks like a good plan for the first introductory course in the theory of algorithms that postpones some things that could be postponed. The model of computation then could be introduced later when proving the undecidability of specific mathematical problems or defining the complexity classes. Still a psychological barrier remains: many people who are quite fluent in mathematics and can easily deal with complicated constructions still have a feeling of uncertainty when they touch the algorithms theory, but at least this barrier becomes more explicit.\footnote{Nowadays the situation is a bit different; one should take into account that most of the people have a lot of programming experience when starting to learn computability. A modern version of Uspensky's approach could be something like that: we start with a programming language that is familiar to the students, and add some library functions: (a)~an interpreter for this language, i.e., a function that gets two inputs, a program string $p$ and some other string $y$, and simulates program $p$ on input $y$; this corresponds to the program axiom; (b) a step-by-step debugger that gets also the number $n$ of steps that should be simulated (a combination of the tracing axiom and program axiom). One can even add a library function without arguments that returns the program text, this would make the fixed point theorem obvious.}

For the proof of incompleteness theorem we need one more axiom (that is not a consequence of the previous ones): the \emph{arithmetization axiom} saying that every computable function can be expressed by an arithmetical formula. (Later this axiom can be proved for some specific computation model.)

If we use this machine-independent approach to the computability theory, we are not allowed to refer to a model of computation when speaking about (say) program transformations or oracle computations. Instead, we should provide all necessary definitions using only the class of computable functions. As we have said, the definition of relative computability that has this form appeared (for the first time) in the master thesis of Uspensky. Then it was done for enumeration reducibility. To deal with program transformations, Uspensky introduced the notion of a ``main numbering'' (see the next section for the enumeration reducibility and main numberings).

One can also note that this axiomatic approach to computability theory provide a formal justification for the following standard observation: most results of the computability theory can be ``relativized'', i.e., remain true if we replace the class of computable functions by the class of $A$-computable function for some oracle $A$. Here $A$ can be a set or a total function. Indeed, one could check that all axioms (except, of course, the arithmetization axiom) for this class. After that we know that all theorems (derived from the axioms) are true for this class.

Uspensky asked whether this observation fully explains the relativization mechanism, i.e., whether a statement that is true for $A$-computable functions for all oracles $A$, is a consequence of his axioms. It turned out that the (positive) answer is easy to get (after the question is stated), see~\cite{Shen1980}.

\subsection*{Computable mappings of sets and enumeration reducibility}

The notion of reducibility introduces by Turing and Post (and considered in the master thesis of Uspensky, see above) can be called ``decision reducibillity''. If $A$ is reducible to $B$, and $B$ is decidable, then $A$ is decidable. One may say that in this definition  we ``reduce the decision problem for $A$ to the decision problem for $B$''. 

In~\cite{1955} Uspensky gives the definition of \emph{enumeration reducibility} where we reduce the task ``enumerate the set $A$'' to the task ``enumerate the set $B$''. This definition uses the notion of a computable operation on sets (introduced in the same paper). Let us describe this notion.

Let us consider the simple case when unary operation is applied to subsets of $\mathbb{N}$ and maps them also to subsets of $\mathbb{N}$. Consider the set $\mathcal{P}(\mathbb{N})$ of all subsets of $\mathbb{N}$ as a topological space. Namely, for each finite set $X\subset \mathbb{N}$ consider the family  $\mathcal{O}(X)$ of all subsets of $\mathbb{N}$ that are supersets of $X$. The families $\mathcal{O}(X)$ and all their unions are considered as open in $\mathcal{P}(\mathbb{N})$. After the topology on $\mathcal{P}(X)$ is defined, we consider all mappings $F\colon \mathcal{P}(\mathbb{N})\to\mathcal{P}(\mathbb{N})$ that are continuous with respect to this topology. It is easy to check that all  continuous $F$ are monotone (if $U\subset V$, then $F(U)\subset F(V)$), and the value $F(U)$ is determined by the values $F(X)$ for finite subsets $X\subset U$ (is the union of $F(X)$ for all finite $X\subset U$. The values of $F$ on finite sets $X$ can be described by the set of pairs $\{\langle n,X\rangle\mid n\in F(X)\}$ (here $n$ is a natural number, and $X$ is a finite set of natural numbers.

Uspensky gives the following definition: a continuous mapping $F\colon \mathcal{P}(\mathbb{N})\to\mathcal{P}(\mathbb{N})$ is a  \emph{computable operation} is the corresponding set of pairs (see above) is an enumerable sets. Note that pairs $\langle n,X\rangle$ are finite objects, so the notion of an enumerable set of pairs makes sense. Now the enumeration reducibility is defined: a set $A\subset \mathbb{N}$ is \emph{enumeration reducible} to a set $B\subset\mathbb{N}$ if there exists a computable operation $F$ that maps $B$ to $A$. Uspensky notes that Turing reducibility can be described in terms of enumeration reducibility: a total  function $\varphi$ is Turing reducible to a total function $\psi$ (i.e., computable with oracle $\psi$) if and only if the graph of $\varphi$ is enumeration reducible to the graph of $\psi$. We can also characterize the Turing reducibility for sets in the same way; for that we consider the graphs of characteristic functions of those sets. He says also that one can characterize partial recursive operators in the sense of Kleene~\cite{Kleene1957}, but here the terminology is confusing (see the discussion below).

Finally, in this paper (\cite{1955}) Uspensky notes that the definition of a computable operations in terms of topology (discussed above) is equivalent to two ``machine-dependent'' definitions. The corresponding notions are called ``Kolmogorov operations'' and ``Post operations'' by Uspensky (though they do not appear explicitly in the works of Kolmogorov and Post).

In another 1955 paper (\cite{1955a}, see also an exposition of its results with some extensions in~\cite{1956}) Uspensky introduces the notion of a numbering (following Kolmogorov's talk given in 1954 at the seminar on recursive arithmetic, Moscow State University mathematics department), introduces the notion of a ``main numbering'' (\rus{<<главная нумерация>>} in Russian) and related the computable operations on enumerable sets (as defined in~\cite{1955}) with algorithmic transformations of their numbers.

Let us explain Uspensky's contribution in more detail. Assume that we want to consider computable transformations of \emph{programs} for computable functions (or enumerable sets). Then it is not enough to know which functions are computable (or which sets are enumerable). We need also to make some assumptions on the ``programming methods'' (or languages, \rus{<<способы программирования>> in Russian} that are used for establish the correspondence between programs and computable functions. Programs are usually strings (words), but one could identify strings with natural numbers via some computable bijective numbering of strings. Then a programming language (method) for computable functions defines a universal function of two arguments: $U(n,x)$ is the output of the $n$th program on input $x$ (we assume that inputs and outputs are also natural numbers). A programming language for enumerable set defines a universal set of pairs $\langle n,x\rangle$ such that $x$ belongs to the $n$th enumerable sets. In a different (but equivalent) language one may say that a programming method for computable functions (resp. enumerable sets) is a \emph{natural numbering} of the set of all computable functions (enumerable sets), i.e., a (total) mapping of $\mathbb{N}$ onto the set of all computable functions (enumerable sets): a number $n$ is mapped to a computable functions (enumerable set) that corresponds to the $n$th program.

Not all programming methods (numbering) are equally good. A reasonable theory that describes the algorithmic transformations of programs needs some additional assumptions. These assumptions essentially appeared in Kleene's work under the name of ``$s$-$m$-$n$-theorem'', but appeared explicitly for the first time in~\cite{1955a} where Uspensky defines the notion of a \emph{main} numbering. This definition consists of two requirements. First, to be main, a numbering should be computable. This means that the corresponding universal function is a computable partial function of two arguments (for the case of sets: the corresponding universal set of pairs is enumerable). Second, any other computable numbering should be \emph{reducible} to the main numbering.\footnote{The definition of reducibility for numbering also was published in~\cite{1955} with a reference to Kolmogorov' seminar talk, also probably for the first time.} This means that for any other computable numbering of the same family there exists a computable translation functions that transforms a number in this other numbering into a number of the same function (set) in the main numbering.

Fix some main numbering for the family of enumerable sets. Then we may define computable mappings of this family into itself. Here computability of a mapping $P$ means that there exist an algorithm that, given a number of some enumerable set $X$, returns (some) number for the set $P(X)$. In other words, we consider computable transformations of programs (or numbers) that preserves the equivalence relation: if two program $p$ and $p'$ are equivalent, i.e., are programs of the same set, then they are transformed into two equivalent programs. Uspensky proved~\cite[Section 6]{1955a} that computable mappings of the family of enumerable sets are exactly computable operations on the family of all sets, restricted to the subfamily of enumerable sets. He also proved a similar statement for a subfamily of function graphs: every computable mapping of the family of computable functions into itself is a restriction of a computable operation on the family on all function graphs. 

Let us describe the connections of this work of Uspensky to the other research of that time.\footnote{Unfortunately (see below the quote from Uspensky's memoirs) all three publications of him~\cite{1955,1955a,1956} are short notes in the \emph{Soviet Math. Doklady}\cite{1955,1955a} and a resume of a talk in the Moscow Mathematical Society~\cite{1956}; they contain only the statements of the theorems and lemmas used in the proofs. The full proofs were published in Uspensky's PhD thesis~\cite{1955b}. Formally speaking, this thesis was publicly available (it can be ordered and accessed in few libraries in the USSR), but it hardly could influence the developments in the field. Probably the short notes~\cite{1955,1955a,1956} were not read outside the USSR, too. Later Uspensky wrote a monograph~\cite{1960} that become his ``habilitation text'' (\rus{<<докторская диссертация>>}); this book was translated into French. Unfortunately, it included only the definition of main numberings, but not the results on computable transformations and mappings.} Rice~\cite{Rice1953} considered \emph{completely recursively enumerable} classes of enumerable set. A family $X$ of enumerable set is called completely recursively enumerable if the set of \emph{all} programs for all elements of $X$ is enumerable. Rice formulated a conjecture~\cite[p.361]{Rice1953}: every completely recursively enumerable family is the family of all supersets of finite sets from some enumerable family of finite sets. This conjecture becomes Theorem 5 in Uspensky paper~\cite[Theorem 5]{1955a} (1955) and is a crucial point in the proofs of his results about computable transformations. This conjecture also was proven in 1956 paper of Rice~\cite{Rice1956} where it is mentioned that the same result was obtained by McNaughton, Myhill and Shapiro (and there is a reference only to a short note of Myhill~\cite{Myhill1955}). Also in the first (1953) paper of Rice it was shown that no non-trivial property of enumerable set can be decided if a program for this set is given (the generalization of this result appeared in~\cite{1955a} as a corollary to Theorem 5). So this statement is usually called ``Rice theorem'', and the result about completely enumerable classes (Rice conjecture proven by Uspensky, McNaughton, Myhill and Shapiro) is usually called ``Rice -- Shapiro theorem'' (see, e.g., Cutland's book~\cite[Chapter 7, \S 2]{Cutland1980}). The connection between computable transformations of programs and computable operations on partial functions was proven (also in 1955) by Myhill and Sheperdson~\cite{MyhillSheperdson1955}, so it is usually called ``Myhill -- Sheperdson theorem''  (see, e.g.,\cite[Chapter 10, \S 2]{Cutland1980}). Since the Rice--Shapiro theorem is its special case, it is also sometimes called ``Myhill--Sheperdson theorem'' (see, e.g.,~\cite[Theorem II.4.2 or Proposition II.5.19]{Odifreddi1989}).

It is hard to tell how the notion of enumeration reducibility was rediscovered. In Rogers' textbook~\cite{Rogers1972} is given without any references (to Uspensky or anybody else). In the 1971 paper ``Enumeration reducibility and partial degrees'' of Case~\cite{Case1971} the references to Rogers' book and Myhill paper~\cite{Myhill1961} are given. However, Myhill's paper (as well as Davis' book~\cite{Davis1958} referenced by Myhill) does not consider enumeration reducibility (it considers only different definitions of relative computability for functions). Modern survey by Soskova~\cite{Soskova2013} does not mention Uspensky's works at all; it contains a reference to a paper of Friedberg and Rogers~\cite{FriedbergRogers1959} that in its turn refers to notes of Rogers' lectures at MIT in 1955--1956 (distributed in 1957) that were a starting point for his book~\cite{Rogers1972}. One may guess that Rogers rediscovered the notion of enumeration reducibility and its name (that is close to the Russian name \rus{<<сводимость по перечислимости>>} used by Uspensky).

The notion of a main numbering (\rus{<<главная нумерация>>} in Uspensky's terminology) was also rediscovered by Rogers (see~\cite{Rogers1958}) under the name of ``G\"{o}del numbering''. Rogers starts with a ``machine-dependent'' definition: ``A G\"{o}del numbering is a numbering equivalent to the standard numbering'' (p.~333); however, later he provides a machine-independent characterization (as the maximal element with respect to reducibility --- as in the Uspensky definition, though without references to Uspensky). Nowadays the names ``admissible numbering'' (see, e.g., Soare's book~\cite{Soare2016}) and ``acceptable numbering'' (see, e.g.,~\cite[Definition II.5.2]{Odifreddi1989}) are used; in both cases a ``machine-dependent'' definition is given.

When comparing Uspensky's work to the similar publications of others, one should have in mind that there are different (and often mixed) notions of reducibility for partial functions. Assume that $f$ and $g$ are two partial functions (with natural arguments and values). Consider the following three definitions of ``$f$ is reducible to $g$'' (=$f$ is computable relative to $g$); each of them is strictly stronger than the previous ones:

\begin{enumerate}

\item The graph of $f$ is enumeration reducible to the graph of $g$.

\item Consider (following Uspensky) the family $\mathfrak{U}$ of all partial functions with natural arguments and values, and consider the following topology in $U$: the basic open sets are sets of all extensions of some finite partial functions. Call a continuous mapping $F\colon\mathfrak{U}\to\mathfrak{U}$ a \emph{computable operation} if its restriction to finite functions has an enumerable graph, i.e., if the set of all pairs $\langle \langle x,y\rangle, u\rangle$, where $x$ and $y$ are natural numbers, $u$ is a finite partial function and $[F(u)](x)=y$,  is enumerable. Then we require that there exists a computable operation that maps $g$ to $f$.

\item We may extend Trakhtenbrot's definition (see the discussion of Uspensky's master thesis above) to partial function and require that $f$ belongs to the closure of the family of all partial recursive functions with $g$ added under substitution, recursion and $\mu$-operation. (This requirement appears, for example, in~\cite{Malcev1965}.)
\end{enumerate}

The third condition in this list can be equivalently reformulated in the oracle computations language. This reformulation repeats the definition from Uspensky's master thesis but allows partial functions (that were not considered by Uspensky). Namely, an algorithm, given $x$, computes $f(x)$; it is allowed to ask questions about $g(y)$ for arbitrary $y$ --- but it should be done sequentially and as soon as it asks for $g(y)$ that is undefined, the computation hangs without providing any result (so $f(x)$ remains undefined for the corresponding $x$). The second requirement also can be reformulated in terms of oracle computations if we allow asking questions about several values $g(y)$ in parallel (the computations continues while waiting for the oracle's answers; it is required that the result of the computation does not depend on delays before the oracle answers are provided).

To see why each requirement is stronger than the previous one, we may consider two examples. The first example separates the first two requirements.

Let $f$ be an arbitrary total function with natural arguments and values. Let $g$ be a partial function whose values are all zeros, and whose domain is the set of all numbers of pairs $\langle n,f(n)\rangle$ for all $n$. (We assume that some computable numbering of pairs is fixed.) Then the first requirement is true for these $f$ and $g$ while the second one is false unless $f$ is computable itself (a computable mapping that maps $g$ to $f$ should map the zero function to $f$, since the zero function extends $g$). This example is mentioned in the Uspensky's footnote to the Russian translation of Rogers' book~\cite[p.362]{Rogers1972} with a reference to D.G.~Skordev; the original argument of Rogers is much more complicated.

The second example~\cite[Proposition II.3.20, with a reference to Sasso's 1971 thesis]{Odifreddi1989} shows that the third property is stronger than the second one. Let $g$ be an arbitrary partial function with natural arguments that has only zero values. Construct another partial function $f$, also with zero values, in the following way: the value $f(n)$ is defined (and equals $0$) if and only if at least of the one values $g(2n)$ and $g(2n+1)$ is defined. Then the second requirement is satisfied for sure: for input $n$ we ask in parallel what are the values $g(2n)$ and $g(2n+1)$; as soon as one of the answers is given, we return $0$. However, if we have to ask the oracle sequentially, this argument does not work: if we first ask for $g(2n)$ and $g(2n)$ is undefined, then $f(n)$ is undefined even if $g(2n+1)$ is defined. (Of course, this is only an explanation why the previous construction is no more valid; to show that indeed the third requirement may be false we need a simple diagonal argument.)

The first requirement corresponds to the notion that is called ``partial recursive operators'' in Rogers' book~\cite[\S 9.8]{Rogers1972}. The second requirement corresponds to what is called ``recursive operators'' in the same book.

Myhill and Sheperdson~\cite[\S 9.8]{Rogers1972} consider ``partial recursive functionals'' and refer to Thesis~I$^{*\dagger}$ from Kleene's book~\cite[p.~332]{Kleene1957}. However, this Thesis (see the top of p.~332) does not use the name ``partial recursive functional'' that does not appear on p.~332 at all. The subject index refers to page~326 for ``partial recursive functional'', but this page does not mention such a notion. It defines the notion of a partial function $\varphi$ that is partial recursive relative to partial functions $\psi_1,\ldots,\psi_k$ that corresponds to our first requirement (enumeration reducibility of graphs) and mentions some ``scheme'' $F$ but does not say whether this scheme $F$ should define a function for all possible $\psi_1,\ldots,\psi_k$ or only for the specific functions. (The numberings of all functions that are computable with an oracle are considered only for the case when the oracle is total.)  Still Myhill and Sheperdson clarify the situation and say that for their result they need partial recursive functionals that are defined (and produce functions) for all arguments that are functions, so essentially they consider the second requirement (as well as Uspensky in his 1955 papers).

Odifreddi in~\cite[Definition II.3.6]{Odifreddi1989} defines partial recursive functionals with reference to Kleene~\cite{Kleene1957}; however, he uses the third version of the definition (a composition of substitutions, recursions and $\mu$-operators applied to partial recursive functions and input functions) --- one that does not appear in~\cite{Kleene1957}. He uses the names ``effectively continuous functional''  or ``recursive operator'' for the second requirement and the name ``partial recursive operator'' for the first one. He uses topological notions in his definitions (as Uspensky did).

% Odifreddi refers to Uspensky's 1955 works and a paper of Nerode of 1957 

Let us summarize the contribution of Uspensky's papers~\cite{1955,1955a,1956}:

\begin{itemize}
\item the historically first definition of enumeration reducibility;

\item the definition of a numbering and reducibility of numbering was published for the first time (with reference to Kolmogorov's talk);

\item the analysis of the properties of numberings of computable functions and enumerable sets needed for the results about program transformation; the definition of main numberings (later rediscovered by Rogers);

\item the proof of Rice's conjecture about completely recursive enumerable classes of enumerable sets (and similar results for functions, including the undecidability of all non-trivial properties of computable functions);
 
\item the definition of a computable operation (in topological terms) and the proof that algorithmic transformations of programs for computable functions or enumerable sets can be described as restrictions of computable operations on functions or sets.

\end{itemize}

As we have said, these achievements were unavailable to the international community and the corresponding results were independently obtained by other researchers (at the same time or a bit later). Let us note, to avoid possible misunderstanding, that Uspensky does \emph{not} consider algorithms that are defined on all programs of \emph{total} functions and give the same results for equivalent programs. The corresponding work of Kreisel, Lacombe and Shoenfild (1959, see~\cite{KreiselLacombeShoenfield1959}) later generalized by Tseitin~\cite{Tseitin1962} to constructive metric spaces, have no intersections with Uspensky's work.

In the following quote from Uspensky's memoirs (\cite[p.~905--907, 912]{2018b}) he recalls his 1955 results and the Third All-Union Mathematical Congress (1956) where these results were presented:

\begin{quote}
In the survey talk (June, 26) ``On algorithmic reductions'' I spoke about four kinds of reductions and relations between them. These four notions are the following: First, \emph{computability reduction} where the task ``compute $f$'' for some function $f$ is reduced to the task ``compute $g$'' for some other function~$g$. Second, the \emph{decidability reduction}: the task ``construct a decision procedure for $A$'', where $A$ is some set, is reduced for the same task for some other set $B$. Third, the enumerability reduction: the task ``enumerate $A$'' for some set $A$ is reduced to the same problem for some other set $B$. Finally, this is \emph{reduction of mass problems} that reduces one mass problem to another one  $\langle\ldots\rangle$ The notion of mass problems was introduced by Yury Medvedev, who was Kolmogorov's student, who defined also the corresponding reductions. $\langle\ldots\rangle$

Another talk of mine (July, 2) was named ``The notion of a program and computable operators'', and a short communication (July, 3) ``Computable operations, computable operators and effectively continuous functions'' was closely related to that talk.

In the last communication I formulated (without proof, of course) the result which now I consider as my main mathematical achievement and still remember the circumstances when it came to my min; it was called ``Theorem~3''\footnote{Theorem~3 was interesting for me from the semiotic viewpoint, even if I did not know the word ``semiotics'' at that time. I remember how I was walking along Moscow streets thinking about this question only. The insight came when I was at my mother-in-law apartment (on Big Spasoglinitschevskii lane in Moscow). My son was not born yet, my wife and her mother went to their jobs in the morning, there was no phone in the apartment (and, of course, no mobiles!). Suddenly I've understood how it works. [Uspensky's footnote]} (see below). This result was the core of my Ph.D. thesis that was defended in October 1955.  I never published the proof of this result, except for the thesis itself; this thesis is available (or at least \emph{was} available) in the math department library. Why? Mostly due to my laziness (shame on me). Another reason, may be less embarrassing, but stupid, was my desire to present this result in the most general form (but one cannot reach the limits of generalization).
$\langle\ldots\rangle$

\textbf{Theorem 3}. \emph{Let $g$ be a function with natural arguments and values. Assume that this function has the following property: if $m$ and $n$ are programs of the same computable $s$-ary function, then $g(m)$ and $g(n)$ are programs of the same unary function. Then there exists a computable operator $V$ such that for every function $\theta$ with program $n$ the value $V(\theta)$ is a function with program $g(n)$.}

\textbf{A philosophical comment}: a semiotic interpretation of Theorem~3 goes as follows: a ``well-behaved'' computable transformation of names is accompanied by a computable transformation of named objects.

\end{quote}

\subsection*{Constructivism and classical mathematics}

The idea of a constructive interpretation of mathematical statements (and, more general, logical connective) goes back to Brouwer and his ``intuitionistic'' school; later it was developed in a different way by Andrei Markov, jr., and his students under the name of ``constructivism''.  In particular, the constructive interpretation of the statement ``for every $x$ there exists $y$ such that\ldots'' is that there exist a way to get this ``existing'' $y$ for every value of $x$.

Usually this constructive approach was combined with the change in the understanding of logical connectives (that makes the excluded middle law invalid). Still there is another possibility that initially was not very popular: consider the ``effective'' versions of classical notions and results as a part of usual (``classical'', ``non-constructive'') mathematics that uses standard mathematical tools. Many people thought that if we are studying algorithms, this should be done in some ``constructive'' or ``finitistic'' way. Uspensky stressed that this is not the only option and one can study constructive notions inside the classical universum of mathematics.

Here are two examples that he considered. The first is the notion of a computable real number. There are different construction of real numbers (Dedekind cuts, fundamental sequences,  common points of intervals of decreasing lengths, decimal expansions, etc.). For each of the constructions one can consider its effective version. For example, we can consider Dedekind cuts such that there exists an algorithm that says for a rational number whether it belongs to the left or right part. For a fundamental sequence $x_n$ of rational numbers one may require this sequence to be computable (given $n$, one can compute $x_n$), and also require the existence of a computable modulus of convergence (an algorithm that, given rational $\eps>0$, computes some $N$ such that $|x_k-x_l|<\eps$ for all $k,l>N$). For an infinite decimal fraction one may require the computability of the function $n\mapsto \text{($n$th digit)}$, and so on.

Each of these definitions leads to some subset of $\mathbb{R}$ that consists of the numbers that have effective representations in the corresponding sense. One can ask (still working in the framework of classical mathematics) whether these definitions lead to the same subset of to different ones. It is not difficult to see that they define the same subset (in different ways), and the elements of this subset can be called computable real numbers (following Turing~\cite{Turing1937}).

This example can be used to illustrate the difference with Markov-style constructivism. For constructivists there are no such things as ``real numbers'' in the usual sense, so they cannot consider the set of computable real numbers as a subset of the set of all real numbers. For the a (computable) real number is a pair of algorithms: one, given $n$, computes $x_n$, and the other computes the modulus of convergence. Note that not all definitions mentioned above are equally good. For example, the definition with decimal fractions has problems: we cannot define addition, i.e., there is no algorithm that transforms two constructive real numbers (i.e., the algorithms for their representations) into their sum (i.e., the corresponding algorithm).

However, as Uspensky notes, the same problem can be analyzed in the framework of classical mathematics. For that, we consider numberings of computable reals that correspond to different definitions. We may ask then whether these numberings are equivalent (whether one can algorithmically transform the number of a computable real in one numbering into a number of the same real in another numbering). And here the same problem with decimal fractions reappears --- and the other positional systems also have this problem. In~\cite{1960} Uspensky provides necessary and sufficient conditions for the reducibility of two numberings of computable reals (with different bases).

Another example studied by Uspensky~\cite{1960a}: the effective versions of the notion of an infinite set of natural numbers. We may say that a set $X$ is infinite if for every natural $n$ the set $X$ contains at least $n$ different elements. Or: $X$ is infinite if it differs from any finite set $X$: for every finite $F$ there exists some number that belongs to the symmetric difference $F\bigtriangleup X$. Both definitions lead to natural effective versions. In the first case we require that there is an algorithm that, given $n$, produces a list of $n$ different elements of $X$. In the second case we require that there is an algorithm that, given a finite set $X$, produces some element of $F\bigtriangleup X$. It is easy to see that these two effective definitions are equivalent (and we may even modify the second definition requiring only that the algorithm gives an element of $X\setminus F$ for finite subsets $F$ of $X$). Using the terminology from Post's paper~\cite{Post1944} all these properties are equivalent to non-immunity of $X$ (i.e., to the existence of an enumerable infinite subset of $X$).

On the other hand, not all definitions of infinity lead to equivalent effective versions. For example, we may say that $x$ is infinite if for every $n$ there exists an initial segment $[0,N]$ that contains at least $n$ elements if $X$. The effective version of this definition would be: there exists an algorithm that for every $n$ computes some $N$ with this property. This is a weaker property of ``effective infiniteness'': as Uspensky noted in~\cite{1957a}(answering the question of Kolmogorov; A.V.~Kuznetsov and Yu.T.~Medvedev independently answered the same question), this requirement means that the set is not hyperimmune in the sense of Post~\cite{Post1944}.

One may also note (though this has no relation to Uspensky's work) that the basic definition in algorithmic randomness, the definition of randomness given by Martin-L\"{o}f in 1966~\cite{MartinLof1966} is also an effective version of the definition of a null set (a set of Lebesgue measure $0$). This classical definition says that a set $X\subset [0,1]$ is a null set if for every $\eps>0$ there exists a covering of $X$ by intervals whose total measure does not exceed $\eps$. For obvious reason we may consider only rational values of $\eps$ and only interval with rational endpoints. Then both $\eps$ and the intervals are constructive objects, and one may consider the effective version of the definition and require that an algorithm gets $\eps>0$ and enumerates the intervals with required properties. This is exactly what Martin-L\"{o}f suggested.

Many topics in algorithmic randomness can be interpreted as effectivization of classical notion and results. For example, the Solovay's criterion of Martin-L\"{o}f randomness is (as Alexander Bufetov noted) the effective version of the Borel -- Cantelli lemma. It turns out that its standard proof (that considers tails of a convergent series) cannot be effectivized and some other argument (also natural and simple) is need, see~\cite{2013} for details. Another instructive example of this type is a proof of an effective version of an ergodic theorem given by Vladimir~Vyugin (a student of Uspensky)~\cite{Vyugin1998}.

\subsection*{Algorithmic information theory}

It is strange that Uspensky, being a student of Kolmogorov and his colleague at the Mathematics Department of the Moscow State University, was not involved in the research initiated by Kolmogorov in 1960s when he introduced the notion of algorithmic complexity of finite objects (now known also as Kolmogorov complexity). I have asked him about that but it still remains a kind a mystery for me. As Uspensky told me, he came into this field only when preparing (with Alexei L. Semenov) the talk for the Urgench conference~\cite{1981,1982a}. In this talk Uspensky and Semenov suggested a general scheme for defining different versions of complexity (or algorithmic entropy, as Uspensky preferred to name them) known at the time: plain, prefix, monotone, decision entropies, as well as conditional versions of entropy. Initially (see~\cite{Shen1984}) this approach used the notions of $f_0$-spaces and their continuous mappings. In a sense this can be considered as an extension of the topological approach to computability suggested by Uspensky long ago.  However, this was definitely an overkill, and Uspensky and Semenov~\cite{1981,1982a} suggested a much more simple version of this scheme that used only the ``compatibility relation'' on objects and descriptions that is enough to cover most of the cases. Later this simplified scheme was explained in~\cite{1992a,1996}; a detailed exposition from the topological viewpoint (but without $f_0$-spaces) can be found in~\cite{2013}

The different notions of randomness are discussed also in a survey~\cite{1990} and in a monograph~\cite{2013}. In 2005 Uspensky gave a talk at the ``Modern mathematics'' school for undergraduates devoted to algorithmic randomness. A brochure based on this talk was published in 2006~\cite{2005} and was reprinted as a part of a monograph~\cite{2013}.

One of the questions asked by Uspensky, Semenov and An.~Muchnik~\cite{1998} remains open. They asked whether the Martin-L\"{of} randomness is equivalent to the absence of a computable strategy in non-monotone games (``non-predictability''). See~\cite{2006,2013} for more details.

\subsection*{Popular science}

There are different ideas about ``popular science'' (in French one says ``vulgarization'', and it sounds embarrassing though partially correct). One may tell stories about life and fate of great scientists. One can try to retell stories found in other popular science books adding more funny jokes. All this may be a good thing, but Uspensky's approach was different. During all his life he tried to explain faithfully the real scientific achievements. These explanation could be easily accessible or technically difficult (depending on the audience); still it was always a serious and honest explanation of a material that can be explained with a clear indication of what remains without proof (or clarification). And he never was afraid of explaining basic and ``well known'' things: as Aristotle wrote in \emph{Poetics}, ``subjects that are known are known only to a few''.

While being a student, Uspensky (with a senior coauthor, Evgeny B. Dynkin) wrote a book~\cite{1952}  that was based on the materials of mathematical circles in Moscow. Uspensky first was a participant of these circles, and later one of the teachers there. The book covers several topics (graphs' coloring, the basics of number theory and probability theory). These topics are presented as a sequence of problems (as it was done in the circles' meetings), and the solutions of these problems are provided. This was not the first problem book based on the materials of mathematical circles, but and important new idea was that these problems, taken together, form a coherent exposition of some mathematical theory. This book for a long time was very hard to find (before it was reprinted in 2004 and before its appearance on the Internet). 

Several popular brochures written by Uspensky were based on his lectures for high school students (in particular, for the participants of the mathematical olympiads) and appeared in the series ``Popular lectures on mathematics''. Some of there were not related to his own mathematical specialty: he wrote a brochure about applications of mechanics to mathematics~\cite{1958a} and about Pascal's triangle~\cite{1966}. The latter includes also a philosophical discussion: what is a combinatorial problem and why do we fix the list of operations that are allowed in the answer for such a problem  (e.g., including factorials but not the notation for binomial coefficients).

Two other brochures in this series written by Uspensky (``The Post machine''~\cite{1979} and ``The G\"{o}del incompleteness theorem'') are covering topics from mathematical logic and algorithms' theory. The first is quite elementary and is based on the lessons given by Uspensky to elementary school students. The other one (as we have mentioned) is based on the article published in \emph{Russian Mathematical Surveys} and assumes significant mathematical culture (but still is accessible to competent high school students). One more popular exposition~\cite{1983} written by Uspensky was devoted to the non-standard analysis where the tools from mathematical logic are used to proved a mathematically correct approach to infinitesimals. The extended version of this brochure was published few years later~\cite{1987}.

Like Josef Knecht (from Hesse's \emph{Das Glasperlenspiel}) Uspensky switched to more and more basic things when becoming older. He started to preach mathematics among humanities students (and researchers). This preaching started in 1960 when he developed and implemented the mathematics curriculum for the Division of Theoretical and Applied Linguistics of the Philology Department of the Moscow State University. However, during the two last decades of his life he addressed to a much wider audience. Several of his lectures during the summer school on mathematics and linguistics (in Dubna, a town near Moscow) were videotaped (thanks to Vitaly Arnold) and are available (see the references in~\url{http://www.mathnet.ru/php/person.phtml?option_lang=rus&personid=20219}). They give some idea about Uspensky's approach to teaching, but one could fully appreciate it only during university courses (first of all, a non-obligatory ones, \rus{<<спецкурсы>>} in Russian). Uspensky always was preaching mathematics, not preaching ``about mathematics''. He explained simple things, but seriously and with proofs. One of his last books~\cite{2009b} is even called ``Very simple examples of mathematical proofs'' (probably not a good name from the advertising viewpoint). The other book~\cite{2000} 	was named ``What is an axiomatic approach?'', and it also contains a lot of examples, including ``school geometry'' --- not the part that is taught in high school but the axiomatic part that is omitted. For example, this book explains how one can derive from the axioms that for every line there is a point that does not belong to this line.

The materials from these two books were included in a collection of Uspensky's paper named ``Mathematics' Apology''\cite{2009}, together with the some other (more general) essays about mathematics. And strangely his preaching was successful --- at least if we interpret success in the same sense as for Saint Anthony of Padua's preaching to the fish: in 2010 Uspensky got the ``Enlightenment'' award established by Dmitry Borisovich Zimin, Russian engineer and philanthropist, the founder and main sponsor of the \emph{Dynasty} foundation.

In addition to his own books, Uspensky organized the translation and publication of many classical textbooks: he translated (following the suggestion of Kolmogorov) R.~Peter's book on recursive functions~\cite{Peter1954}, was the editor for the translations of monographs of Kleene~\cite{Kleene1957}, Rogers~\cite{Rogers1972}, Davis~\cite{DavisNonStandard1980} (the latter translation probably was the first Russian-language book about non-standard analysis), Church's logic textbook~\cite{Church1960}, the first volume of the ``Elements of mathematics'' by Bourbaki~\cite{BourbakiSetTheory1965}, and Ashby's book on cybernetics~\cite{Ashby1959}.


\clearpage
\begin{thebibliography}{99}

\item[]\hspace{-\labelwidth}\hspace{-\labelsep}\textsl{Uspensky's papers in the Internet:}

\bibitem{mathnet}
Uspensky's page at \texttt{mathnet.ru}: \url{http://www.mathnet.ru/rus/person20219}

\bibitem{kafedra}
Uspensky's page at the Logic and Theory of Algorithms Division of the Mathematical Department of Moscow Lomonosov State University:
\url{http://lpcs.math.msu.su/~uspensky/}

\item[]\hspace{-\labelwidth}\hspace{-\labelsep}\textbf{Publications}

\bibitem{1949}
A geometric approach to proving main properties of harmonic functions [\rus{Геометрический вывод основных свойств гармонических функций}, in Russian] \emph{\rus{Успехи математических наук}}, 1949, vol.~IV, issue 2(30), p.~201--205,
\url{http://lpcs.math.msu.su/~uspensky/bib/Uspensky_1949_UMN_Geometr_vyvod.pdf},
\url{http://mi.mathnet.ru/umn8612}

\bibitem{1952}
\emph{A general definition of algorithmic computability and algorithmic reducibility. \rus{Общее определение алгоритмической вычислимости и алгоритмической сводимости}, in Russian}. Master thesis (advisor A.~Kolmogorov) Moscow State Lomonosov University, Math. Department  \rus{механико\dash математический факультет}.  A typescript. 90 pp. \url{http://lpcs.math.msu.su/~uspensky/bib/Uspensky_1952_Diploma.pdf}
\url{https://archive.org/details/uspensky-1952-master-thesis-and-reviews}

\bibquote{[Minutes of the Division of History of Mathematics Meeting, May 10, 1952. Participants: A.N.~Kolmogorov, P.S.~Novikov, S.A.~Yanovskaya, I.G.~Bashmakova.  ``After hearing the talks of the student, V.A.~Uspensky, of his thesis advisor, A.N.~Kolmogorov, of the reviewer, P.S.~Novikov, and of S.A.~Yanovskaya, the Division considers the work of V.A.~Uspensky exceptional. It contains several new results and shows a deep understanding of the difficult topic of algorithms theory. The exposition is excellent; the paper should be published.  The Division's head Pr. S.A.~Yanovskaya. May 12, 1952. [\rus{<<Заслушав выступления студента В.\,А.\,Успенского, руководителя работы ак. А.\,Н.\,Колмогорова, рецензента П.\,С.\,Новикова и С.\,А.\,Яновской, кафедра постановила: признать работу В.\,А.\,Успенского выдающейся. Отметить, что работа содержит ряд значительных новых результатов и свидетельствует о глубоком владении автором всей трудной проблематикой теории алгоритмов. Отметить также прекрасное оформление работы и признать необходимым опубликование её. Зав.~кафедрой проф.~C.\,А.\,Яновская, 12 мая 1952 года.>>]} The advisor and reviewer's opinions: \url{http://lpcs.math.msu.su/~uspensky/bib/Uspensky_1952_Diploma_reviews.pdf}, see also \url{https://archive.org/details/uspensky-1952-master-thesis-and-reviews}]}

\bibquote{\nb{\rus{Из рецензии Новикова: <<Таким образом, автор не только дал методологически правильное\emdash материалистическое\emdash объяснение причин эквивалентности различных определений алгоритма, но и получил возможность включить их в единую теорию>>.}}}

\bibitem{1952a} 
E.B.~Dynkin, V.A.~Uspensky, \emph{Mathematical Conversations: Multicolor Problems, Problems in the Theory of Numbers, and Random Walks}. [\rus{Е.\,Б.\,Дынкин, В.\,А.\,Успенский, \emph{Математические беседы. Задачи о многоцветной раскраске. Задачи из теории чисел. Случайные блуждания}, in Russian]} (Mathematical circles' library, vol.~6. [\rus{Библиотека математического кружка, вып. 6}]). Moscow, Leningrad: State Publisher of Technical and Theoretical Books. [\rus{Москва\endash Ленинград: Государственное издательство технико\dash теоретической литературы}], 1952.  \url{http://ilib.mccme.ru/djvu/bib-mat-kr/besedy.htm}, \url{http://www.math.ru/lib/book/djvu/bib-mat-kr/besedy.djvu}. 2d ed.: Nauka [\rus{Наука},] 2004. English translation published by D.C.~Heath, Boston in 1963 (in three brochures) and later was published as one book: E.B.~Dynkin, V.A.~Uspenskii, \emph{Mathematical Conversations: Multicolor Problems, Problems in the Theory of Numbers, and Random Walks}, Dover books in mathematics, Dover Publications, 2006, ISBN 0-486-45351-0. German translation:  E.B.~Dynkin, W.A.~Uspenski, \emph{\ger{Mathematische Unterhaltungen. I. Mehrfarbenprobleme.  II. Aufgaben aus der Zahlentheorie. III. Aufgaben aus der Wahrscheinlichkeitsrechnung. Irrfahrten (Markoffsche Ketten)}, Deutscher Verlag der Wissenschaften, Berlin, 1955--1956.} Reprinted in 1966, 1976, 1979 (Aulis Verlag Deubner \& Ko KG Cologne), 1983. Turkish translation of parts I and III was published in 1962,  part II in 1963 (Turkish Mathematical Society Publications).

\bibitem{1953} On the notion of algorithmic reducibility. [\rus{О понятии алгоритмической сводимости}, in Russian] A summary of a talk given at the meeting of Moscow Mathematical Society, March 17, 1953.   \emph{\rus{Успехи математических наук}}, vol.~VIII, issue~4(56), 1953, July--August, p.~176,  see~\url{http://mi.mathnet.ru/umn8234}

\bibquote{\nb{A short exposition of the Master Thesis}}

\bibitem{1953a} G\"{o}del' theorem and the theory of algorithms [\rus{Теорема Гёделя и теория алгоритмов}, in Russian]. A summary of a talk given at the meeting of Moscow Mathematical Society, March 24, 1953.  \emph{\rus{Успехи математических наук}}, vol.~VIII, issue~4(56), 1953, July--August, p.~176--178,  see~\url{http://mi.mathnet.ru/umn8234}

\bibitem{1953b}
G\"{o}del' theorem and the theory of algorithms [\rus{Теорема Гёделя и теория алгоритмов}, in Russian], \emph{\rus{Доклады Академии наук СССР}}, vol.~91, issue 4, p.~737--740 (1953), \url{https://istina.msu.ru/publications/article/92662634/}, \url{https://archive.org/details/uspensky-1953-dan-091-4-godel-algorithms}. English translation: G\"{o}del's theorem and the theory of algorithms,  \emph{American Mathematical Society Translations, Series 2, Advances in the Mathematical Sciences}, vol.~23 (1963), 103--107, \url{DOI 10.1090/trans2/023/06}

\bibitem{1955}
On computable operations [\rus{О вычислимых операциях}, in Russian], \emph{\rus{Доклады Академии наук СССР}}, vol.~103, issue~5 (1955), p.~773--776, \url{https://istina.msu.ru/publications/article/92662640/},  \url{https://archive.org/details/uspensky-1955-dan-103-5-computable-operations}. English translation: \url{https://archive.org/details/uspensky-1955-dan-103-5-eng}

\bibitem{1955a}
Systems of enumerable sets and their numberings [\rus{Системы перечислимых множеств и их нумерации}, in Russian], \emph{\rus{Доклады Академии наук СССР}}, vol.~105, issue~6 (1955), p.~1155--1158, \url{https://istina.msu.ru/publications/article/92662649/}, \url{https://archive.org/details/uspensky-1955-dan-105-6-enumerable-sets-numerations}.

\bibitem{1955b} On computable operations [\rus{О вычислимых операциях}, in Russian], Ph.D thesis, Moscow State Lomonosov University, Mathematics Department, October 1955. \url{https://archive.org/details/uspensky_thesis_1955}.

\bibitem{1956} Computable operations and the notion of a program. [\rus{Вычислимые операции и понятие программы}, in Russian]. A summary of a talk given at the meeting of Moscow Mathematical Society, February 28, 1956,  \emph{\rus{Успехи математических наук}}, vol.~XI, issue~4(70), 1956, July--August, p.~172--176, \url{http://mi.mathnet.ru/umn7861}

\bibquote{\nb{\rus{Потенциально вычислимая нумерация вычислимых функций: универсальная функция вычислима; вполне накрывающая --- если сводится всякая потенциально вычислимая, главная второго рода --- если потенциально вычислима и вполне накрывающая. Существуют потенциально вычислимые нумерации, являющиеся вполне накрывающими, а также не являющиеся таковыми. <<Соображения этого пункта дают основания предложить понятие главной нумерации второго рода в качестве уточнения понятия ``способ программирования''>>. Конструктивные операторы: преобразования вычислимых функций в вычислимые, для которых существует вычислимое преобразование номеров в главной нумерации. Вычислимые операторы: определены на всех частичных функциях, соответствуют операторам перечисления на графиках (рекурсивные операторы в смысле Роджерса~\cite{Rogers1972}). Теорема 1: оператор, продолжаемый до вычислимого, является конструктивным. Теорема 2: всякий конструктивный оператор продолжается до вычислимого.  Теорема 3: если потенциально вычислимая нумерация такова, что всякий вычислимый оператор является относительно неё конструктивным, то эта нумерация главная. (Теорема отсутствует в заметке~\cite{1955a}.) Теорема 4: всякое нетривиальное разбиение множества функций на две части задаёт неразрешимое разбиение номеров в главной нумерации второго рода. Интерпретация как связности пространства.}}} 

\bibitem{1956a}
Third All-Union Mathematical Congress. Plenary talk ``On algorithmic reducibility''. [\rus{<<Об алгоритмической сводимости>>}, in Russian], June 26, 1955.  Talk ``The notion of a program and computable operators'' [\rus{<<Понятие программы и вычислимые операторы>>}, in Russian], July 2,1955. Short talk ``Computable operations, computable operators and constructively continuous function'' [\rus{<<Вычислимые операции, вычислимые операторы и конструктивно-непрерывные функции>>}, in Russian], July 3, 1955. The resumes of the talks is published: Proceedings of the Third All-Union Mathematical Congress [\emph{\rus{Труды третьего всесоюзного Математического съезда}}, in Russian]. Moscow: Academy of Science Publications, 1956. Vol.~2, p.~66--69 (plenary talk), vol.~1, p.~ 186 (talk), vol.~1, p.~185 (short talk).

\nb{no scan?}

\bibitem{1957}
On the uniform continuity theorem. [\rus{К теореме о равномерной непрерывности}, in Russian]. \emph{\rus{Успехи математических наук}}, vol.~XII, issue~1(73), 1957, January--February, p.~100--142, \url{http://mi.mathnet.ru/umn7524}

\bibitem{1957a}
Some remarks on [recursively] enumerable sets [\rus{Несколько замечаний о перечислимых множествах}, in Russian]. \emph{Zeitschrift f\"ur mathematische Logik und Grundlagen der Mathematik}, Bd.~3, Heft 12, S.~157--170 (1957), \url{https://onlinelibrary.wiley.com/toc/15213870/1957/3/12}, \url{https://istina.msu.ru/publications/article/92666188/}, \url{https://archive.org/details/uspensky-1957-zml-zamechanie-perechisl-mnozhestv}. English translation: Some remarks on recursively enumerable sets, \emph{American Mathematical Society Translations, Series 2, Advances in the Mathematical Sciences}, vol.~23 (1963), 89--101, \url{DOI 10.1090/trans2/023/05}

\bibquote{\nb{\rus{Система всех бесконечных линейных перечислимых множеств не допускает вычислимой нумерации, множество нижних точек перечислимого множества может не быть перечислимым, классификация перечислимых множеств, гипериммунные множества как множества, у которых прямой пересчёт не мажорируется вычислимой функцией.}}}

\bibitem{1958}
A.N.~Kolmogorov, V.A.Uspensky, On the definition of an algorithm [\rus{К определению алгоритма}, in Russian]. \emph{\rus{Успехи математических наук}}, vol.~XIII, issue~4(82), 1958, July--August, p.~3--28, \url{http://mi.mathnet.ru/umn7453}. English translation (Elliott Mendelson): Kolmogorov A.N., Uspenskij V.A., On the definition of an algorithm,  \emph{American Mathematical Society Translations, Series 2, Advances in the Mathematical Sciences}, vol.~29 (1963), 217--245, \url{DOI 10.1090/trans2/029/07}

\bibitem{1958a}
\emph{Some application of mechanics to mathematics} [\emph{\rus{Некоторые приложения механики к математике}}, in Russian] (Popular mathematics lectures [\rus{Популярные лекции по математике}], issue~27). Moscow, State Physics and Mathematics Publishers [\rus{Государственное издательство физико\dash математической литературы}], 1958. 48~pp., \url{https://math.ru/lib/plm/27} English translation: V.A.~Uspenskii, \emph{Some applications of mechanics to mathematics} (Popular lectures in mathematics. Vol.~3). Oxford, London, New York, Pergamon Press, 1961. An English translation was also published by Mir publishing house (Moscow) in 1979. They also published a French translation in a collection: \fra{N.~Vilenkine, G.~Chilov, V.~Ouspenski, J.~Lioubitch, L.~Chor, Quelques applications des math\'{e}matiques. M\'{e}thode des approximations successives. Gamme simple (structure de l'\'{e}chelle musicale). Quelques applications de la m\'{e}canique aux math\'{e}matiques. M\'{e}thode cin\'{e}matique dans les probl\`{e}mes de g\'{e}om\'{e}trie, Traduit du russe par Djilali Embarek, Initiation aux Math\'{e}matiques, Editions Mir, M., 1975 , 279 pp.} A Spanish translation was published by Mir publisher in 1979. 

\bibitem{1960}
On the relations between different systems of constructive real numbers [\rus{К вопросу о соотношении между различными системами конструктивных действительных чисел}, in Russian], \emph{Communication of High Education Institutions [\rus{Известия высших учебных заведений}}, mathematics,  1960, issue 2(15), p.~199--208, \url{http://mi.mathnet.ru/ivm2028}

\bibitem{1960a}
\emph{Lectures on computable functions} [\emph{\rus{Лекции о вычислимых функциях}}]. Moscow State Physics and Mathematics Publishers [ \rus{Государственное издательство физико\dash математической литературы}], 1960. 492 pp., \url{https://archive.org/details/uspenskij-1960ru}. French translation: Ouspenski V.A., Le\c cons sur les fonctions calculables. Paris, Hermann, 1966. 412 p.

\bibitem{1964}
M.V.~Lomkovskaya, E.V.~Paducheva, V.A.~Uspensky, Linguistic calculuses. [\rus{М.\,В.\,Ломковская, Е.\,В.\,Падучева, В.\,А.\,Успенский. Лингвистические исчисления.}. In Russian] 
\emph{Proceedings of the Fourth All-Union Mathematical Congress}, Leningrad, \rus{\emph{Труды IV Всесоюзного математического съезда}, Ленинград}, 1961, vol.~2, 83--90 (1964). 
\nb{нет текста, по описанию в Zentralblatt в обратном переводе, \url{https://zbmath.org/?q=an:0178.32801}}

\bibitem{1966}
\emph{Pascal's Triangle} [\emph{\rus{Треугольник Паскаля}}, in Russian].  (Popular mathematics lectures [\rus{Популярные лекции по математике}], issue~43). Moscow, Nauka Publishers, Physics and Mathematics Division [\rus{Наука, главная редакция физико\dash математической литературы}], 1966. 35 pp. Second extended edition: 1979. 48~pp. \url{http://www.math.ru/lib/book/plm/v43.djvu} English translation: V. A. Uspenskii, \emph{Pascal's triangle}, Translated and adapted from the Russian by David J. Sookne and Timothy McLarnan, Popular Lectures in Mathematics, The University of Chicago Press, Chicago--London, 1974, vii+35 pp. An English translation was also published by Mir Publishers, Moscow, 1979. 

\bibitem{1969}
Reductions of computable and potentially computable numberings [\rus{О сводимости вычислимых и потенциально вычислимых нумераций}, in Russian] \emph{\rus{Математические заметки}}, vol.~6, issue~1 (1969), p.~3--9, \url{http://mi.mathnet.ru/mz6891}/ English translation: V.A.~Uspenskii, Reduction of computable and potentially computable numerations, \emph{Mathematical Notes of the Academy of Sciences of the USSR}, 1969, vol.~6, no.~1, 461--464, \url{DOI10.1007/BF01450246}, \url{https://istina.msu.ru/publications/article/92645436/} 

\bibitem{1974}
An elementary exposition of G\"odel's incompleteness theorem [
\rus{Теорема Гёделя о неполноте в элементарном изложении}, in Russian]. \emph{\rus{Успехи математических наук}}, vol.~XXIX, issue~1(175), p.~3--47 (1974, January--February), \url{http://mi.mathnet.ru/umn4322}. English translation: (E.~Lichfield): Uspenskii V.A., An elementary exposition of G\"odel's incompleteness theorem, \emph{Russian Mathematical Surveys}, vol.~29 (1974), no.~1, 63--106, \url{DOI 10.1070/RM1974v029n01ABEH001280}, full text available at\url{https://istina.msu.ru/publications/article/92645484/} 

\bibitem{1979}
\emph{Post's machine} [\emph{\rus{Машина Поста}}, in Russian].
(Popular mathematics lectures [\rus{Популярные лекции по математике}], issue~54). Moscow, Nauka Publishers, Physics and Mathematics Division [\rus{Наука, главная редакция физико\dash математической литературы}], 1979. 96 pp., \url{http://www.math.ru/lib/book/plm/v54.djvu} Second edition: 1988. English translation: V.A.~Uspensky, Post's machine, Translated from the Russian by R. Alavina, with a supplement by Emil L.~Post, Little Mathematics Library, Mir (Moscow), 1983 , 88 pp. Portuguese translation was published by Mir Publishers, Moscow, 1985, 

\bibitem{1980}
V.A.~Uspenski, “Nonstandard analysis”, Fiz.-Mat. Spis. Bulgar. Akad. Nauk., 23(56):3 (1980/81), 219--231 (Bulgarian) (Translated from the Russian by Vl. Sotirov), \url{http://www.mathnet.ru/eng/person20219}. 
\nb{Текста нет, переводом чего является --- непонятно.}

\bibitem{1981}
Uspensky V.A., Semenov A.L., What are the gains of the theory of algorithms: Basic developments connected with the concept of algorithm and with its application in mathematics, \emph{Algorithms in Modern Mathematics and Computer Science, Proceedings, Urgench, Uzbek SSR, September 16--22, 1979}. Edited by A.P.~Ershov and D.~Knuth, Lecture Notes in Computer Science, 122, Springer, 1981, p.~100--234, \url{https://archive.org/details/uspensky-semenov-urgench}

\bibitem{1982}
\emph{G\"odel's incompleteness theorem} [\emph{\rus{Теорема Гёделя о неполноте}}, in Russian]. (Popular mathematics lectures [\rus{Популярные лекции по математике}], issue~57). Moscow, Nauka Publishers, Physics and Mathematics Division [\rus{Наука, главная редакция физико\dash математической литературы}], 1982,  \url{http://www.math.ru/lib/book/plm/v57.djvu}. English translation by N.~Kolblitz was published in1987 (Mir Publisher) and later in \emph{Theoretical Computer Science}, \textbf{130} (1994),  239--319, \url{https://www.sciencedirect.com/science/article/pii/0304397594902224}

\bibitem{1982a}
Uspensky V.A., Semenov A.L. Theory of algorithms: main discoveries and applications. [\rus{Успенский В.А., Семёнов А.Л., Теория алгоритмов: основные открытия и приложения}, in Russian]. Published in the proceedings volume:  \emph{Algorithms in modern mathematics and its applications. The proceedings of the international symposium, Urgench, Uzbek SSR, September 16--22, 1979.} [\emph{\rus{Алгоритмы в современной математике и её приложениях. Материалы международного симпозиума, Ургенч, УзССР, 16--22 сентября 1979 г.}}] Edited by A.P.~Ershov, D.~Knuth, Part I, p.~99--342, \url{https://archive.org/details/uspensky-semenov-urgench-rus}

\bibitem{1983} \emph{Nonstandard, or non-Archimedian, analysis}. [\emph{\rus{Нестандартный, или неархимедов, анализ}}, in Russian]. Moscow, \rus{Знание} publishers, 1983. 61~pp. \url{https://archive.org/details/uspensky-nonstandard-znanie}

\bibitem{1983a}
Uspensky V.A., Kanovei V.G., Luzin's problems on constituents and their fate
[\rus{Проблемы Лузина о конституантах и их судьба}, in Russian]. \emph{\rus{Вестник Московского университета. Серия 1: Математика, механика}}, 1983, issue~6, p.~73--87, \url{https://istina.msu.ru/publications/article/93856782/}. English translation: Uspenskii V.A., Kanovei V.G., Luzin's problems on constituents and their fate, \emph{Moscow University Mathematics Bulletin}, vol.~38, no.~6 (1983), 86--102. (Allerton Press, inc.)

\nb{no scan}

\bibitem{1985} 
 Luzin's contribution to the descriptive theory of sets and functions: concepts, problems, predictions [\rus{Вклад Н.\,Н.\,Лузина в дескриптивную теорию множеств и функций: понятия, проблемы, предсказания}, in Russian]. \emph{\rus{Успехи математических наук}}, vol.~40, issue~3(243), p.~85--116, \url{http://mi.mathnet.ru/umn2648}. English translation: Uspenskii V.A., Luzin's contribution to the descriptive theory of sets and functions: concepts, problems, predictions, \emph{Russian Mathematical Surveys}, vol.~40, no.~3 (1985), 97--134, available at \url{https://istina.msu.ru/publications/article/92645592/} 

\bibitem{1986}
A.L.~Semenov, V.A.~Uspensky, Mathematical logic in computer science and computer applications [\rus{Математическая логика в вычислительных науках и вычислительной практике}, in Russian] \emph{\rus{Вестник Академии наук СССР}}, \textbf{56}(7),  93--103 (1986)

\bibitem{1986a} 
V.A.~Uspensky, On nonstandard analysis. Selecta Mathematica Sovietica, \textbf{5}, 357--396 (1986). The translation of the Editor's Preface to~\cite{DavisNonStandard1980}.

\bibitem{1987} \emph{What is nonstandard analysis?} [\emph{\rus{Что такое нестандартный анализ?}}, in Russian]. Moscow,  Nauka Publishers, Physics and Mathematics Division [\rus{Наука, главная редакция физико\dash математической литературы}], 1987. 128~pp., \url{https://archive.org/details/chto-takoe-nestandartny-analiz-djvu}


	\bibitem{1987a} V.A.~Uspensky, A.L.~Semenov, Theory of algorithms: main discoveries and applications. [\rus{В.\,А.\,Успенский, А.\,Л.\,Семёнов, \emph{Теория алгоритмов: основные открытия и приложения}}, in Russian], Moscow,  Nauka Publishers, Physics and Mathematics Division [\rus{Наука, главная редакция физико\dash математической литературы}], 1987.  (Programmer's library, vol.~49 [\rus{Серия <<Библиотечка программиста>>, выпуск 49}].) 288~pp., \url{https://archive.org/details/uspensky-semenov-1987-algoritmy}. English translation: Vladimir Uspensky, Alexei Semenov, \emph{Algorithms: Main Ideas and Applications}, Kluwer Academic Publishers, 1993, \url{https://doi.org/10.1007/978-94-015-8232-2}

\bibitem{1987b}
A.N.~Kolmogorov, V.A.~Uspensky, Algorithms and randomness [\rus{Алгоритмы и случайность}, in Russian], \emph{\rus{Теория вероятностей и её применения}}, vol.~XXXII, issue~3, 1987 (July, August, September), p.~425--455, \url{http://mi.mathnet.ru/tvp1437}. English translation: Kolmogorov A.N., Uspenskii V.A., Algorithms and randomness, \emph{Theory of Probability and Its Applications}, SIAM Publishers, vol.~32, no.~3, 389--412. \url{http://dx.doi.org/10.1137/1132060}, \url{https://istina.msu.ru/publications/article/92647240/} See also: A. N. Kolmogorov, V. A. Uspensky, ``Algorithms and randomness'', Proceedings of the 1st World Congress of the Bernoulli Society (1986, Tashkent), v. 1, VNU Sci. Press, Utrecht, 1987, 3–53. Bulgarian translation: (Pl.~Mateev): Fiz.-Mat. Spis. Bulgar. Akad. Nauk., \textbf{33(66)}:3--4 (1991), 223--247.


\bibitem{1988}
V.A.~Uspensky, V.G.~Kanovei, M.Ya.~Suslin's contribution to set\dash theoretic mathematics [\rus{Вклад М.\,Я.\,Суслина в теоретико\dash множественную математику}, in Russian]. \emph{\rus{Вестник Московского университета. Серия 1: Математика, механика}}, 1988, issue~5, p.~8--12, \url{https://istina.msu.ru/publications/article/93856823/}. English translation: Uspenskii V.A., Kanovei V.G., M.Ya.~Suslin's contribution to set\dash theoretic mathematics, \emph{Moscow University Mathematics Bulletin}, vol.~43, no.~5 (1988), 29--40. (Allerton Press, inc.)

\nb{no scan}

\bibitem{1990} 
V.A.~Uspensky, A.L.~Semenov, A.~Shen, Can an individual sequence of zeros and ones be random? [\rus{Может ли (индивидуальная) последовательность нулей и единиц быть случайной?}, in Russian]. \emph{\rus{Успехи математических наук}}, vol.~45, issue~1(271), p.~105--162 (1990, January--February), \url{http://mi.mathnet.ru/umn4692}. English translation: V.A.~Uspensky, A.L.~Semenov, A.~Shen, Can an individual sequence of zeros and ones be random? \emph{Russian mathematical surveys}, vol.~45, issue 1, 121--189 (1990), \url{https://archive.org/details/uspensky-semenov-shen-1990}

\bibitem{1991}
V.A.~Uspensky, V.E.~Plisko, Diagnostic propositional formulas [\rus{Диагностические пропозициональные формулы}, in Russian]. \emph{\rus{Вестник Московского университета. Серия 1: Математика, механика}}, 1991, issue~3, p.~7--12, available at \url{https://istina.msu.ru/publications/article/92717778/}

\bibitem{1991a}
V.A.~Uspensky, N.K.~Vereshchagin, V.E.~Plisko, \emph{An introductory course of mathematical logic} [\rus{\emph{Вводный курс математической логики}}], Moscow State Lomonosov University Publishers, 1991. 2nd edition: Moscow, \rus{Физматлит}, 2004, \url{https://archive.org/details/uspensky-vereshchagin-plisko}  

\bibitem{1992} 
\emph{Complexity and Entropy: An Introduction to the Theory of Kolmogorov Complexity}. In: \emph{Kolmogorov Complexity and Computational Complexity}, Osamu Watanabe, editor. Springer, 1992, ISBN 3-540-55840-3, p.~85--102, \url{https://archive.org/details/uspensky-1992-watanabe-book}

\bibitem{1992a}
Kolmogorov and mathematical logic, \emph{The Journal of Symbolic Logic}, volume 57, number 2, June 1992, 385--412, \url{https://doi.org/10.2307/2275276}, \url{https://archive.org/details/uspensky-1992-jsl-kolmogorov-mathematical-logic} See also: Sitzungsberichte der Berliner Mathematischen Gesellschaft, Berliner Math. Ges., Berlin, 1992, 41--74.

\bibitem{1995}
Vladimir A. Uspensky and Valery Ye. Plisko, Review: Raymond M. Smullyan, G\"odel's incompleteness theorems, \emph{The Journal of Symbolic Logic},  volume 60, issue 4 (1995), 1320--1324, \url{https://projecteuclid.org/euclid.jsl/1183744885}
	
\nb{scan available}

\bibitem{1996} V.A.~Uspensky, A.~Shen, Relations Between Varieties of Kolmogorov Complexities, \emph{Mathematical Systems Theory}, \textbf{29}, 271--292 (1996), \url{https://link.springer.com/article/10.1007/BF01201280}, 
\url{lpcs.math.msu.su/~uspensky/bib/Uspensky_1996_MST_Shen_Relations_between_varieties_of_Kolmogorov_complexities.pdf}

\bibitem{1996a}
Kolmogorov complexity: recent research in Moscow. In W.Penczek, A.Szalas (eds.), \emph{Proceedings of the 21st International Symposium on Mathematical Foundations of Computer Science 1996 (MFCS96), Cracow, Poland, September 2--6, 1996}  (Lecture Notes in Computer Science, v.~1113), 1996, p.~156--166,
\url{https://link.springer.com/chapter/10.1007/3-540-61550-4_145}, \url{http://lpcs.math.msu.su/~uspensky/bib/Uspensky_1996_LNCS_Kolmogorov_Complexity_Recent_trents_Moscow.pdf}

\nb{trents???}
\nb{file available}

\bibitem{1997}
Mathematical logic in the former Soviet Union: brief history and current trends,
\emph{Logic and Scientific Methods}, M.L.~Dalla Chiare et al., editors, Kluver Academic Publishers,  \url{https://www.springer.com/us/book/9780792343837}, 457--483.

\nb{dile available}

\bibitem{1998}  Andrei A.~Muchnik, Alexei L.~Semenov, Vladimir A.~Uspensky, Mathematical metaphysics of randomness, \emph{Theoretical Computer Science}, \textbf{207}, 263--317 (1998),  \url{https://www.sciencedirect.com/science/article/pii/S0304397598000693} (full text available), \url{http://lpcs.math.msu.su/~uspensky/bib/Uspensky_1998_TCS_Muchnik_Semenov_Math_metaphysics_randomness.pdf}

\bibitem{2000}
\emph{What is an axiomatic method?} [\rus{\emph{Что такое аксиоматический метод?}} Izhevsk, Udmurt university [\rus{Ижевск: издательский дом <<Удмуртский университет>>}], 2000. 100~pp. ISBN 5-7029-0337-4.

\nb{TeX file available}

\bibitem{2001}
Why Kolmogorov complexity? In: E.Goles and C.Martinez (eds.), \emph{Complex systems} (Series: Nonlinear Phenomena and Complex Systems, Vol. 6),  Kluwer Academic Publishers, 2001, p.~201--260. ISBN 0-7923-6830-4, \url{https://link.springer.com/chapter/10.1007/978-94-010-0920-1_5}

\nb{file and TeX source of a preliminary version}

\bibitem{2002}
\emph{Non-Mathematical Works, with an Appendix: Semiotic Letters of A.N.~Kolmogorov to the author and his friends}, [\emph{\rus{Труды по нематематике с приложением семиотических посланий А.\,Н.\,Колмогорова к автору и его друзьям}.}, in Russian] Moscow, \rus{ОГИ}, 2002. 1409 pp.,  ISBN 5-94282-086-4, \url{http://www.math.ru/lib/book/pdf/shen/usp/usp-all.pdf}

\bibitem{2003}
B.~Durand, V.~Kanovei, V.A.~Uspensky, N.K.~Vereshchagin. Do stronger definitions of randomness exist? \emph{Theoretical Computer Science}, v.~290, No.~3, p.~1987--1996 (2001), available as \url{https://www.sciencedirect.com/science/article/pii/S0304397502000403} 

\bibitem{2005}
V.G.~Kanovei, V.A.~Uspenksy, On the equivalence of the two versions of the Continuum Hypothesis, [\rus{Об эквивалентности двух форм континуум\dash гипотезы}, in Russian]. \emph{\rus{Вестник Московского университета. Серия 1: Математика, механика}}, 2005, issue~3, p.~62--64, \url{https://istina.msu.ru/publications/article/100287822/} 

\nb{page photos exist}

\bibitem{2006}
Four algorithmic faces of randomness, [\rus{Четыре алгоритмических лица случайности}, in Russian]. \emph{\rus{Математическое просвещение}}, ser.~3, issue~10, Moscow, MCCME, 2006, p.~71--108, \url{http://mi.mathnet.ru/mp188}. Included as an appendix in the book~\cite{2013}.

\bibitem{2006a}
Kolmogorov as I remember him
[\rus{Колмогоров, каким я его помню}, in Russian]. In: \emph{Kolmogorov in the memories of his students} [\emph{\rus{Колмогоров в воспоминаниях учеников}}, in Russian], edited by A.N.~Shiryaev and N.G.~Khimchenko, Moscow, MCCME, 2006, 272--371.

\bibitem{2006b}
V.G.~Kanovei, V.A.~Uspensky, On the uniqueness of nonstandard extensions [\rus{О единственности нестандартных расширений}, in Russian]. \emph{\rus{Вестник Московского университета. Серия 1: Математика, механика}}, 2006, issue 5, p.~3--10, available at \url{https://istina.msu.ru/publications/article/100285758/}

\bibitem{2008}
V.G.~Kanovei, T.~Linton, V.A.~Uspensky, A game approach to the Lebesgue measure, [\rus{Игровой подход к мере Лебега}, in Russian] \emph{\rus{Математический сборник}}, vol.~199, issue~11 (2008), 21--44, \url{http://mi.mathnet.ru/msb3948} English translation: V.G.Kanovei, Tom Linton and Vladimir A.~Uspensky, Lebesgue measure and gambling, \emph{Sbornik: Mathematics}, volume 199, no.~11, p.~1597, \url{http://dx.doi.org/10.1070/SM2008v199n11ABEH003974}

\bibitem{2009}
\emph{Mathematics' apology} [\emph{\rus{Апология математики}}, in Russian]. Saint-Petersbourg, \rus{Амфора}, 2009. 554~pp.

\nb{file available}

\bibitem{2009a}
On the history of Goldbach's problem [\rus{К истории проблемы Гольдбаха}, in Russian].  In: \rus{\emph{Историко\dash математические исследования. Вторая серия}. РАН, Институт естествознания и техники им. С.\,И.\,Вавилова}. Vol.~13(48). Moscow, \rus{Янус-К}, 2009,  ISBN 978-5-8037-0449-2, p.~273--283.

\nb{file available}

\bibitem{2009b}
\emph{Basic examples of mathematical proofs} [\emph{\rus{Простейшие примеры математических доказательств}}, in Russian]. (<<\rus{Математическое просвещение}>> series, issue~34). Moscow, MCCME, 2009. 56 pp. ISBN 978-5-94057-492-7.  \url{http://www.math.ru/lib/book/pdf/mp-seria/034_uspensky.pdf}

\bibitem{2010}
V.A.~Uspensky, V.V.~Vyugin, The emergence of algorithmic information theory in Russia [\rus{Становление алгоритмической теории информации в России}, in Russian]. \emph{\rus{Информационные процессы}}, vol.~10, issue~2, p.~145--158.

\bibitem{2011}
G\"{o}del's theorem and four ways to it, [\rus{Теорема Гёделя и четыре дороги, ведущие к ней}, in Russian]. \emph{\rus{Математическое просвещение}}, ser.~3, issue 15, Moscow, MCCME, 2011, p.~35--75, \url{http://mi.mathnet.ru/mp309}

\bibitem{2013} 
N.K.~Vereshchagin, V.A.~Uspensky, A.~Shen, \emph{Kolmogorov complexity and algorithmic randomness} [\emph{\rus{Колмогоровская сложность и алгоритмическая случайность}}, in Russian]. Moscow, MCCME, 2013. 575~pp. (English version: A.~Shen, V.~Uspensky, N.~Vereshchagin, \emph{Kolmogorov complexity and algorithmic randomness}, American Mathematical Society, 2017. \url{https://www.lirmm.fr/~ashen/kolmbook-eng-scan.pdf})

\bibitem{2014} A.~Semenov, S.~Soprunov, V.Uspensky, The lattice of definability. Origins, recent developments, and further directions, \emph{CSR 2014: Computer science --- theory and applications}, Lecture Notes in Comput. Sci., v.~8476, Springer, 2014, 23--38. \url{https://doi.org/10.1007/978-3-319-06686-8_3}

\bibitem{2017}
[V.A.~Uspensky, M.S.~Gelfand], Mathematics is a part of humanities [\rus{Математика\emdash это гуманитарная наука (интервью с В.\,А.\,Успенским ведёт М.\,С.\,Гельфанд)}]. In the book: \emph{\rus{Математические прогулки. Сборник интервью}}, Moscow, \rus{Паулсен}, 2017. ISBN 978-5-98797-057-7, p.~198--207. [The English translation of the book: \emph{Mathematical Walks. A Collection of Interviews}, ISBN: 978-5-98797-167-3, Moscow, Paulsen, 2017.]

\bibquote{\nb{files available}}

\bibitem{2018}
Vladimir Uspenskiy and Alexander Shen, Algorithms and Geometric Constructions, \emph{Computability in Europe, 2018}, Lecture Notes in Computer Science, v.~10936, Springer,  p.~410--420 (2018), see also \texttt{arXiv:1805:12579}

\bibitem{2018a}
\emph{Non-Mathematical Works. Second extended and corrected edition. In five parts. Part 5. Memories and observations}. [\emph{\rus{Труды по \textbf{не}математике. Второе издание, исправленное и дополненное. В пяти книгах. Книга пятая. Воспоминания и наблюдения.}}, in Russian]. Moscow, \rus{Объединённое гуманитарное издательство. Фонд <<Математические этюды>>}. 2018. 1118~pp. 

\bibitem{2018b} Third Mathematical Congress. In:~\cite[p.~897--905]{2018a}. Commentaries. Ibidem, p.~905--912.

\item[]\hspace{-\labelwidth}\hspace{-\labelsep}\textbf{Publications translated or edited by Uspensky}

\bibitem{Peter1954}
R.~Peter, \emph{Recursive functions} [\rus{Р.\,Петер, \emph{Рекурсивные функции}}, in Russian], translated from German by V.A.~Uspensky, edited by A.N.~Kolmogorov, with a preface written by A.N.~Kolmogorov. Moscow,  \rus{Издательство иностранной литературы}, 1954. (Original book: \emph{Rekursive Funktionen}, von R\'osza P\'eter, Budapest, 1951.)

\bibitem{Kleene1957} 
S.~Kleene, Introduction to metamathematics [\rus{Стефен К.\,Клини, \emph{Введение в метаматематику}}, in Russian], translated from English by A.S.~Esenin-Volpin, edited by V.A.~Uspensky. Moscow,  \rus{Издательство иностранной литературы}, 1957. (Original book: Stephen Cole Kleene, \emph{Introduction to metamathematics}, D. van Nostrand company, New York, Toronto, 1952.)

\bibitem{Ashby1959}
W.~Ross Ahsby, \emph{An introduction to cybernetics} [\rus{У. Росс Эшби, \emph{Введение в кибернетику}}, in Russian]. Translated from English by D.G.~Lakhuti. Edited by V.A.~Uspensky. With a preface written by A.N.~Kolmogorov. Moscow, \rus{Издательство иностранной литературы}, 1959. 428~pp. (Original book: \emph{An introduction to cybernetics}, by W. Ross Ashby, London, Chapman\&Hall Ltd., 1956.)

\bibitem{Church1960}
A.~Church, Introduction to mathematical logic, volume I. [\rus{А.\,Чёрч, \emph{Введение в математическую логику, I}}, in Russian]. Translated from English by V.S.~Chernyavsky. Edited by V.A.~Uspensky. Moscow, \rus{Издательство иностранной литературы}, 1960. (Original book: \emph{Introduction to mathematical logic} by Alonzo Church. Volume I. Princeton University Press, 1956.)

\bibitem{BourbakiSetTheory1965}
N.~Bourbaki, \emph{Elements of Mathematics. Part 1. Fundamental structures of analysis. Book 1. Set theory} [\rus{Н.\,Бурбаки, \emph{Начала математики. Первая часть. Основные структуры анализа. Книга первая.  Теория множеств}}, in Russian]. Translated from French by G.N.~Povarov and Yu.A.~Shikhanovich. Edited by V.A.~Uspensky. With a preface written by V.A.~Uspensky. Moscow, \rus{Мир}, 1965. (Original book: \'El\'ements de math\'ematique par N.~Bourbaki, XVII, XX, XXII, I. Premiere partie. Les structures fondamentales de l'analyse. Livre I. Th\'eorie des ensembles. Hermann, 1956--1960.)

\bibitem{MathematicsModernWorld1967}
\emph{Mathematics in the modern world}  [\emph{\rus{Математика в современном мире}}, in Russian],  a collection of translations of papers from Scientific American special issue \emph{Mathematics in the modern world}, Scientific American, 1964. Translated from English by N.G.~Rychkova. Edited by V.A.~Uspensky. With a preface written by V.A.~Uspensky. Moscow, \rus{Мир}, 1967.

\bibitem{Rogers1972}
H.~Rogers, Theory of recursive functions and effective computability [\rus{Х.\,Роджерс, \emph{Теория рекурсивных функций и эффективная вычислимость}}, in Russian], translated from English by V.A.~Dushsky, M.I.~Kanovich, E.Yu.~Nogina. Edited by V.A.~Uspensky. Moscow, \rus{Мир}, 1972. (Original book: Hartley Rogers, Jr., \emph{Theory of recursive functions and effective computability}, McGraw-Hill Book Company, 1967.A preliminary version with the same name. Volume I. Mimeographed. Technology Store, Cambridge, Mass., 1957.)

\bibitem{DavisNonStandard1980}
M.~Davis, \emph{Applied nonstandard analysis} [\rus{М.\,Дэвис, \emph{Прикладной нестандартный анализ}}, in Russian]. Translated from English by S.F.~Soprunov. Edited by V.A.~Uspensky. With a preface written by V.A.~Uspensky. Moscow, \rus{Мир}, 1980. (Original book: Martin Davis, \emph{Applied nonstandard analysis}, Wiley\&Sons, 1977.)

\clearpage
\item[]\hspace{-\labelwidth}\hspace{-\labelsep}\textbf{Other references}

\bibitem{Skolem1923}
Th. Skolem, \emph{Begr\"undung der elementaren Arithmetik durch dir rekurrierende Denkweise ohne Anwendung scheinbarer Ver\"anderlichen mit unendlichem Ausdehnungsbereich} (Videnskapsselskapets Scrifter, I. Mat.-naturv. Klasse, 1923, No.~6),  Kristiania, 1923. (English translation: The foundations of elementary arithmetic established by means of the recursive mode of thought without the use of apparent variables ranging over infinite domains. In~\cite[p.~302--333]{vanHeijenoort1967}.)

\bibitem{Hilbert1926}
David Hilbert, \"Uber das Unendliche, \emph{Mathematische Annalen}, Bd.~95, S.~161--190 (1926). (English translation:  On the Infinite,  \cite[p.~367--392]{vanHeijenoort1967}.)

\bibitem{Ackermann1928}
Wilhelm Ackermann in G\"ottingen, Zum Hilbertschen Aufbau der reellen Zahlen, 
\emph{Mathematische Annalen}, Bd.~99, 118--133 (1928). (English translation: On Hilbert's construction of the real numbers,  \cite[p.~493--507]{vanHeijenoort1967}.)

\bibitem{Godel1931}
Kurt G\"odel in Wien, \"Uber formal unentscheidbare S\"atze der Principia Mathematica und verwandter Systeme I, \emph{Monatshefte f\"ur Mathematik und Physik}, \textbf{38}, 173--198 (1931). (English translation: On formally undecidable propositions of \emph{Principia Mathematica} and related systems I, \cite[p.~596--616]{vanHeijenoort1967} or \cite[p.~4--38]{Davis1965}.)

\bibitem{Herbrand1932}
J.~Herbrand \`a Paris, Sur la non-contradiction de l'Arithm\'etique, \emph{Journal f\"ur die reine und angewandte Mathematik}, Bd.~166, S.~1--8 (1932), \url{http://www.digizeitschriften.de/dms/img/?PID=PPN243919689_0166%7Clog4} (English translation: On the consistency of arithmetic~\cite[p.~618--628]{vanHeijenoort1967}.)

\bibitem{Godel1934} 
Kurt G\"odel, \emph{On undecidable propositions of formal mathematical systems}, Lecture notes, Institute for Advanced Study (Princeton), Spring 1934. Reprinted in~\cite[p.~39-74]{Davis1965}

\bibquote{\rus{\nb{Раздел 9. General recursive functions, после примера определения, выходящего за рамки примитивной рекурсии: ``One may attempt to define this notion [of a general recursive function] as follows: if $\phi$ denotes an unknown function, and $\psi_1,\ldots,\psi_k$ are known functions, and if the $\psi$'s and the $\phi$ are substituted in one another in the most general fashions and certain parts of the resulting expressions are equated, then if the resulting set of functional equations has one and only one solution for $\phi$, $\phi$ is a recursive function.'' Примечание к этому определению: This was suggested by Herbrand in a private communication. К нему добавлено (при издании сборника): ``A slightly different definition was given by him in J. r. ang. Math. 166 (1932), p.~5 [это статья~\cite{Herbrand1932}], where he postulated `computability'. However, also in this definition he did not require computability by any definite formal rules (note the phrase `consider\`ee intuitionistiquement' and footnote 5. In intuitionistic mathematics the two Herbrand definitions are trivially equivalent. In classical mathematics the non-equivalence of general recursiveness with the first mentioned concept of Herbrand was proved by L.~Kalm\'ar in Zs. f. math. Log. u. Grundl. d. Math. 1 (1955) p.93. Whether Herbrand's second concept is equivalent with general recursiveness is a largely epistemological question which has not yet been answered.''
 %
Приводится пример определения Аккермана. ``We shall make two restrictions on Herbrand's definition. The first is that the left-hand side of each of the given functional equations defining $\phi$ shall be of the form \[\phi(\psi_{i1}(x_1,\ldots,x_n),\psi_{i2}(x_1,\ldots,x_n),\ldots,\psi_{il}(x_1,\ldots,x_n)).\] The second (as stated below) is equivalent to the condition that all possible sets of arguments $(n_1,\ldots,n_l)$ of $\phi$ can be so arranged that the computation of the value of $\phi$ for any given set of arguments $(n_1,\ldots,n_l)$ by means of the given equations requires a knowledge of the values of $\phi$ only for sets of arguments which precede $(n_1,\ldots,n_l)$ (не сказано, в каком смысле). Дальше определяются правила вывода и говорится: Now our second restriction on Herbrand's definition of recursive function is that for each set of natural numbers $k_1,\ldots,k_l$ there should be one and only one $m$ such that $\phi(k_1,\ldots,k_l)=m$ is a derived equation.}}}

\bibitem{Peter1934}
R\'osza P\'eter,  \"Uber den Zusammenhang der verschiedenen Begriffe der rekursiven Funktionen, \emph{Mathematische Annalen}, Bd.~110,  n.~1, S.~612--632 (1935),  \url{https://doi.org/10.1007/BF01448046},  \url{https://link.springer.com/article/10.1007%2FBF01448046}

\bibquote{\rus{\nb{Вводится термин primitive rekursion для операций, использованных Гёделем в~\cite{Godel1931}: ``Die einfachste Form einer solchen Rekursion ist jene, die G\"odel in seiner zitierten Arbeit verwendet; diese werde ich im folgenden als ,,\emph{primitive Rekursion}" bezeichnen.''  Но сами функции называются просто ``rekursiv''. Доказывается, что разные схемы сводятся к примитивной рекурсии.}}}

\bibitem{Church1936}
Alonzo Church, An unsolvable problem of elementary number theory, \emph{American Journal of Mathematics}, vol.~58, no.~2 (April 1936), 345--363. Reprinted in~\cite[p.~88--107]{Davis1965}

\bibquote{\rus{\nb{приводится (раздел 4) определение ``recursive function'' через исчисление равенств. Доказывается, что минимизация, если даёт всюду определённую функцию, не выводит из этого класса (теорема IV на с.353). Доказывается, что класс совпадает с $\lambda$-определимыми, и доказывается неразрешимость (какого-то вида) Entscheidungsproblem}}}

\bibitem{Kleene1936} S.C.~Kleene, General recursive functions of natural numbers, \emph{Mathematische Annalen}, Bd.~112, S.~727--742 (1936), \url{https://eudml.org/doc/159849}. Reprinted in~\cite[p.~236--253]{Davis1965}.

\bibquote{\rus{\nb{Появляется термин primitive recursive в современном смысле, и говорится о general recursive functions по Эрбрану и Гёделю, в терминах выводов в исчислении равенств (и даже есть два варианта исчисления, которые дают один и тот же класс функций, 2a 2b). Доказывается теорема о нормальной форме, но вместо $\mu$-оператора написан $\eps$-оператор, который на с.728 объясняется как наименьшее число, удовлетворяющее условию, или нуль, если такого числа нет\emdash впрочем, применяется он только к ситуациям, когда решение есть, и вообще рассматриваются только всюду определённые функции. Теорема о нормальной форме позволяет сказать, что определение 2c на с.738, где говорится о нормальной форме, эквивалентно предыдущим.  Отмечается, что класс систем равенств, которые задают функции, не является recursively enumerable (не является областью значений всюду определённой вычислимой функции). }}}

\bibitem{Post1936}
Emil L. Post. Finite combinatory processes. Formulation I. \emph{The Journal of Symbolic Logic}, vol.~1 (1936), p.~103--105.

\bibitem{Turing1937}
A.M.~Turing, On computable numbers, with an application to the Entscheidungsproblem, \emph{Proceedings of the London Mathematical Society}, ser.~2, vol.~42 (1936--7), p.~230--265; Correction in the next volume of the same journal: vol.~43 (1937), p.~544--546. Reprinted in~\cite[p.~116--154]{Davis1965}.

\bibitem{Kleene1938} 
S.C.~Kleene, On notations for ordinal numbers, \emph{Journal for Symbolic Logic}, \textbf{3}, 150--155 (1938), \url{https://www.jstor.org/stable/2267778}

\bibquote{\nb{\rus{Определяются рекурсивные (Herbrand-G\"odel recursive) и частично рекурсивные (partial recursive) функции с помощью исчисления равенств. Отмечается эквивалентность с определениями Тьюринга и с Church--Kleene $\lambda$-definability. Определяется $\mu$-оператор на частичных функциях и утверждается, что он не выводит из класса частично рекурсивных функций, с намёком на доказательство на с. 152--153. Далее определяются системы обозначений для ординалов.}}}

\bibitem{Turing1939}
A.M.~Turing, Systems of logic based on ordinals, \emph{Proc. London Math. Soc.} (2), vol. 45 (1939), pp.~161--228, \url{https://doi.org/10.1112/plms/s2-45.1.161}.  (Reprinted~\cite[p.~154--222]{Davis1965}. Turing's Princeton Ph.D. thesis (with the same name) availables as~\url{http://www.dcc.fc.up.pt/~acm/turing-phd.pdf})

\bibitem{Kleene1943}
S.C.~Kleene, Recursive predicates and quantifiers, \emph{Transactions of the American Mathematical Society}, \textbf{53}, number 1, 41--73 (1943), \url{https://doi.org/10.1090/S0002-9947-1943-0007371-8}. Reprinted in~\cite[p.254--287]{Davis1965}.

\bibquote{\rus{\nb{section 2: General recursive functions. We shall proceed to the Herbrand--G\"odel generalization of the notion of recursive function. Раздел 2, определение на с.~44--45 для всюду определённых функций (в терминах выводимости из равенств). ``A function $\phi$ which can be defined from given functions $\psi_1,\ldots,\psi_k$ by a series of applications of general recursive schemata we call \emph{general recursive} in the given functions; and in particular, a function $\phi$ definable ab initio by these means we call \emph{general recursive}.'' Однако относительная вычислимость явно не рассматривается (вводится мимоходом, как шаг в последовательности операций, и для частичных функций не определяется вовсе).
Раздел 3, вводится $\mu$-оператор (для случая, когда он даёт всюду определённую функцию), доказано, что он не выводит за пределы общерекурсивных (в смысле Эрбрана--Гёделя). Теорема II говорит про арифметическую иерархию.  Раздел 6 начинается с определения partial recursive functions с помощью исчисления Эрбрана--Гёделя, требуется, чтобы было выводимо не более одного утверждение о значении функции. Говоится, что получится замкнутый относительно минимизации класс (теорема III), но, кажется, не объясняется отчётливо, как применяется минимизация, если функция частична. Теорема IV говорит, что частично рекурсивные (в этом смысле) фукнции представимы в нормальной форме (где один оператор минимизации, и он применяется к всюду определённой функции). Её Corollary на с.53 говорит, что можно определить general recursive functions и partial recursive functions с помощью подстановки, рекурсии и минимизации. На с. 60 effective calculability упоминается в интуитивном смысле, и формулируется Thesis I. Every effectively calculable function (effectively decidable predicate) is general recursive. Теорема VIII говорит, что для некоторого предиката (дополнения самоприменимости) нет полной теории: ``This is the famous theorem of G\"odel on formally undecidable propositions, in a generalized form''.}}}

\bibitem{Post1944}
Emil L. Post, Recursively enumerable sets of positive integers and their decision problems, \emph{Bulletin of the American Mathematical Society}, \textbf{5}, 284--316 (1944). \url{https://projecteuclid.org/download/pdf_1/euclid.bams/1183505800}. Reprinted in~\cite[p.304--337]{Davis1965}.

\bibquote{\rus{\nb{определяются и рассматриваются 1-сводимость, m-сводимость, tt-сводимость (в том числе ограниченная), простые, гиперпростые и креативные множества. Раздел 11: General (Turing) reducibility, со ссылкой на Тьюринга~\cite{Turing1939}, в терминах машин с оракулом. Утверждеется, что это столь же окончательное определение относительной вычислимости, как и машины без оракула для (просто) вычислимости. <<A corresponding formulation of ``Turing reducibility'' should then be the same degree of generality for effective reducibility as say general recursive function is for effective calculability.>>  Предлагается план с гипергиперпростыми множествами для построения неполного перечислимого множества (и отмечается, что неизвестно, выйдет ли из этого что-то). Проблема Поста: <<As a result we are left completely on the fence as to whether there exists a recursively enumerable set of positive integers of absolutely lower degree of unsolvability than the complete set $K$, or whether, indeed, all recursively enumerable sets of positive integers with recursively unsolvable decision problems are absolutely of the same degree of unsolvability. On the other hand, if this question can be answered, that answer would seem to be not far off, if not in time, then in the number of special results to be gotten on the way.>>}}}

\bibitem{Kleene1950}
	S.C.~Kleene, A symmetric form of G\"odel's theorem (Presented to the American Mathematical Society, October 29, 1949. Communicated by Prof. L.E.J.~Brouwer at the meeting of April 29, 1950). Koninklijke Nederlandse Akademie van Wetenschappen, Volume 53, deel 6 (1950), 800--802, \url{http://www.dwc.knaw.nl/DL/publications/PU00018825.pdf}

\bibitem{Rice1953} 
H.G.~Rice, Classes of recursively enumerable sets and their decision problems,
\emph{Transactions of the American Mathematical Society}, vol.~74 (1953), p.~358--366.

\bibquote{\rus{\nb{Complete r.e. class of r.e. sets: all indices form a r.e. set. R.e. class: all sets with numbers in some r.e. sets. Описан способ задания c.r.e. класса с помощью перечислимого семейства конечных множеств (все надмножества), высказана гипотеза, что так получаются все, но доказано лишь, что если входит конечное множество, то входят всего его надмножества:  ``We now give a method for constructing c.r.e. classes which seems to be very general'',  ``we venture the conjecture that every c.r.e. class has a key array''.  Есть утверждение о неразрешимости нетривиальных свойств перечислимых множеств: ``If $P$ is any property possessed by some, but not all, r.e. sets, then there exists no effective general method for deciding, given a set $\alpha$ by means of a partial recursive function enumerating it, whether or not $\alpha$ has the property $P$.'' Пишет, что большая часть результатов входит в его диссертацию под руководством Paul Rosenbloom. Представлено 16 ноября 1951 года.}}}

%\bibitem{Markov1954} %!
%А.\,А.\,Марков (младший), Теория алгорифмов, \emph{Тр. МИАН СССР}, \textbf{42}, 3--375 (1954).

\bibitem{Kalmar1955}
L\'aszlo Kalm\'ar in Szeged, Ungarn, \"Uber ein Problem, betreffend die definition des Begriffes der allgemein-rekursiven Funktion, \emph{Zeitschrift f\"ur mathematische Logik und Grundlagen der Mathematik}, Bd.~1, S.~93--95 (1955), \url{https://onlinelibrary.wiley.com/doi/abs/10.1002/malq.19550010204}

\bibitem{Myhill1955}
John Myhill, A fixed point theorem in recursion theory, abstract,  Eighteenth Meeting of the Association of Symbolic Logic, \emph{The Journal of Symbolic Logic}, volume 20, no.2  (June 1955), p.~205.
\bibquote{\nb{``Rice conjectured that co
nversely every c.r.e. class can be written in the form $\Sigma T(\alpha_i)$. We can use the fixed-point theorem to prove this conjecture (which was proved also in another way by MacNaughton and Shapiro).'' (\rus{Доказательство не приводится})}}


\bibitem{MyhillSheperdson1955}
J.~Myhill in Berkeley,  California (USA), J.C.~Sheperdson in Bristol, England,
Effective operations in partial recursive functions, \emph{Zeitschrift f\"ur mathematische Logik und Grundlagen der Mathematik}, Bd.~1, S. 310--317 (1955),
\url{https://onlinelibrary.wiley.com/doi/abs/10.1002/malq.19550010407}

\bibitem{Rice1956}
H.G.~Rice, On completely recursive enumerable classes and their key arrays,
\emph{The Journal of Symbolic Logic}, volume 21, number 3, Sept.~1956, p.~304--308. (Received September 14, 1955) \url{https://www.jstor.org/stable/2269105}

\bibquote{\rus{\nb{Доказывается, что любой c.r.e. class задаётся key array (среди прочего) --- этот результат, пишет автор, получили McNaughton (не опубликовано), Myhill (ссылка на \cite{Myhill1955}) и Norman Shapiro (без ссылки)}}}

\nb{\rus{проверить, в чём разница между ссылками на МакНотона и Шапиро!}}

\bibitem{Davis1958}
Martin Davis, \emph{Computability and Unsolvability}, McGraw-Hill Book Company, 1958.

\bibitem{Rogers1958}
Hartley Rogers, Jr. G\"odel numberings of partial recursive functions, \emph{Journal of Symbolic Logic}, Volume 23, Number 3, Sept. 1958, p.~331--341 (Received July 7, 1958),  \url{https://www.jstor.org/stable/2964292}

\bibquote{\nb{%
``Intuitively, a G\"odel numbering is an association of numbers with partial recursive functions such that the following three condition hold:

i) we are able effectively tell whether or not a number is associated with a partial recursive function, i.e., the set of numbers associated is recursive;

ii) there is an effective procedure such that given any number associated with a function, we can find instructions for effectively computing that function;

iii) there is an effective procedure such that given instructions for effectively computing a partial recursive function, we can find an integer associated with that function''

(p.~331)

``Definition 1. a \emph{numbering} $\pi$ is a mapping of a recursive set of integers $D_\pi$, called the \emph{domain} of $\pi$, onto the set of partial recursive functions. 

Definition 2. A numbering $\pi$ is \emph{semi-effective} if there exists a partial recursive function of two variables $\Phi$ such that for every $i\in D_{\pi}$, $\Phi(i,x)$ is identical, as a partial function of $x$ with $\pi i$. Any such $\Phi$ determines a numbering. We shall say that $\Phi$ \emph{describes} $\pi$.

Definition 3. A numbering $\pi$ is \emph{fully effective} if there is a partial recursive function of two variables $\Phi$ and a recursive function $f$ such that: $\Phi$ describes $\pi$, $f$ takes all values in $D_\pi$; and, for all $i$, $\Phi(f(i), x)$ is identical, as a partial function of $x$, with $\phi_i$.''  

(p.~332; \rus{здесь} $\phi_i$\emdash \rus{частично рекурсивная функция с номером $i$ в стандартной нумерации}).

``Definition 4. Two numberings, $\rho$ and $\pi$, are \emph{equivalent} if there exists a recursive function $g$ mapping $D_\rho$ into $D_\pi$ and a recursive function $h$ mapping $D_{\pi}$ to $D_{\rho}$ such that $\rho=\pi g$ on $D_{\rho}$ and $\pi = \rho h$ on $D_{\pi}$.

[into/to???? \rus{проверить}]

Definition 5. A \emph{G\"odel numbering} is a numbering equivalent to the standard numbering.

While this definition as an equivalence class is \emph{invariant}, it is not \emph{intrinsic}. That is to say, it still depends on the initial choice of \emph{some} member of the equivalence class. It is of interest to find an intrinsic definition, if possible. This is accomplished as follows.

Definition 6. A numbering $\rho$ is \emph{derivable} from numbering $\pi$ if there exists a recursive function $g$ mapping $D_{\rho}$ into $D_{\pi}$, such that $\rho=\pi g$ on $D_{\rho}$.  $\langle \ldots\rangle$

Theorem. The partial order of semi-effective numberings possesses a unique maximal element, and this element is the class of G\"odel numberings.''

(p. 333--334)
}}

\bibitem{FriedbergRogers1959} 
Richard M. Friedberg,  Hartley Rogers jr., Reducibility and Completeness for Sets of Integers, \emph{Zeitschrift f\"ur mathematische Logik und Grundlagen der Mathematik}, Bd.~5, S.~117--125 (1959), \url{https://doi.org/10.1002/malq.19590050703}

\bibitem{KreiselLacombeShoenfield1959}
Kreisel, G., Lacombe, D., Shoenfield, J.R., 
\emph{Partial recursive functionals and effective operations}, in \emph{Constructivity in Mathematics}, A.~Heyting, editior, North Holland, 1959, p.~195--207.

\bibitem{Myhill1961} John Myhill, Note on degrees of partial functions, \emph{Proceedings of the American Mathematical Society}, \textbf{12} (1961), p.~519--521, \url{https://doi.org/10.1090/S0002-9939-1961-0125794-X }

\bibitem{Tseitin1962}
G.S.~Tseitin, Algorithmic operators in constructive metric spaces. [\rus{Алгорифмические операторы в конструктивных метрических пространствах}, in Russian]. Proceedings of Moscow Steklov Institute, LVII, Problems in constructive mathematics, 2. (Constructive analysis) [\rus{\emph{Труды математического института имени В.\,А.\,Стеклова, LVII, Проблемы конструктивного направления в математике, 2 (Конструктивный математический анализ)}}], a collection of papers. Edited by N.A.~Shanin. Moscow, Leningrad, Academy of Science Publishers, \rus{Издательство Академии наук СССР, Москва, Ленинград}, 1962, p.~295--361.

\bibitem{Davis1965}
\emph{Basic Papers On Undecidable Propositions, Unsolvable Problems and Computable Functions}, a collection compiled by Martin Davis, Raven Press, Hewlett, New York, 1965.

\bibitem{Malcev1965}
A.I.~Maltsev, \emph{Algorithms and recursive functions} [\rus{\emph{Алгоритмы и рекурсивные функции}}, in Russian], Moscow, \rus{Наука}, 1965. (2nd edition, 1986)

\bibitem{MartinLof1966}
Per Martin-L\"of, The definition of random sequences, \emph{Information and Control}, volume 9, issue 6, December 1966, p.~602--619, \url{https://doi.org/10.1016/S0019-9958(66)80018-9}

\bibitem{vanHeijenoort1967} 
\emph{From Frege to G\"odel. A Source Book in Mathematical Logic, 1879--1931}, a collection of papers compiled by Jean van Heijenoort,  Harvard University Press, Cambridge, Massachusetts, 1967. 

%\bibitem{Zaliznyak1967} %!
%А.\,А.\,Зализняк, \emph{Русское именное словоизменение}, М.:Наука, 1967.

\bibitem{Case1971}
John Case, Enumeration reducibility and partial degrees, \emph{Annals of mathematical logic}, vol.~2, no.~4 (1971), 419--439. (Received 9 September 1969), \url{https://www.sciencedirect.com/science/article/pii/0003484371900039}

\bibitem{Shen1980}
A.~Shen, Axiomatic approach to the theory of algorithms and relativized computability,
[\rus{Аксиоматический подход к теории алгоритмов и относительная вычислимость}, in Russian], \emph{\rus{Вестник Московского университета. Серия 1: Математика, механика.}}, 1980, issue~2, p.~27--29. (English translation made by the author is available as~\url{https://hal-lirmm.ccsd.cnrs.fr/lirmm-01923123}.)

\bibitem{Kleene1981}
Stephen C. Kleene, The theory of recursive functions, approaching its centennial. (Elementarrekursiontheorie vom h\"oheren Standpunkte aus.) \emph{Bulletin of the American Mathematical Society}, volume 4, number 1, July 1981, 43--61.

\bibquote{\nb{\rus{название primitive recursion, утверждает Клини, введено Петер в 1934 году, до этого говорили просто о рекурсивных функциях, начиная с Гёделя в 1931 году. Термин <<рекурсия>> (но не <<рекурсивные функции>>) есть у Сколема в 1923 и Гильберта в 1926.}}}

\bibitem{Cutland1980}
Nigel Cutland, \emph{An introduction to recursive function theory}, Cambridge University Press, 1980. Russian translation by Albert A.~Muchnik, edited by S.Yu.~Maslov: \rus{Н.\,Катленд, \emph{Вычислимость. Введение в теорию рекурсивных функций.} М.:Мир, 1983.}

\bibquote{\nb{\rus{Частично рекурсивные функции определяеются (глава 3, раздел 2) с помощью рекурсии, подстановки и минимизации, но даётся ссылка на Гёделя и Клини 1936, впрочем, без указания конкретной работы. Теорема Майхилла -- Шепердсона, глава 10, параграф 2, теорема Райса -- Шапиро, глава 7, параграф 2}}}

\bibitem{Shen1984} 
A.~Shen, Algorithmic versions of the notion of entropy [\rus{Алгоритмические варианты понятия энтропии}, in Russian], \emph{\rus{Доклады Академии наук}}, 1984, vol.~276, issue~3, p.~563--566. (English translation: Soviet Math. Doklady, \textbf{29}(3), 1984, 569--573.)

%\bibitem{Muchnik1985} %!
%Ан.\,А.\,Мучник, Об основных структурах дескриптивной теории алгоритмов. \emph{Доклады АН СССР}, \textbf{285}(2), 280--283 (1985)

%\bibitem{Ershov1985} %!
%А.\,П.\,Ершов и др., \emph{Основы информатики и вычислительной техники}, М.: Просвещение, 1985 (ч.~1), 1986 (ч.~2).

\bibitem{Odifreddi1989}
Piergiorgio Odifreddi, \emph{Classical Recursion Theory. The Theory of Functions and Sets of Natural Numbers} (Studies in logic and the foundations of mathematics, volume 125), Elsevier, 1989. xix+668 pages.

\bibquote{\nb{\rus{partial recursive functions p.127, closed under composition, primitive recursion and unrestricted $\mu$-operator ссылка на Клини 1938
recursive (general пропускается), p.22 - примитивно рекурсивные плюс минимизация, если результат (и аргумент) всюду определен. Ссылка на Kleene 1936}}}

\bibitem{Soare1996} 
Robert I. Soare, Computability and Recursion, \emph{Bulletin of Symbolic Logic}, \textbf{2}(3), 284--321 (1996). See also \url{http://www.people.cs.uchicago.edu/~soare/History/compute.pdf}

\bibitem{Vyugin1998}
V.V.~Vyugin, Ergodic theorems for individual random sequences, \emph{Theoretical Computer Science}, volume 207, issue 2,  November 6, 1998, p.~343--361, \url{https://doi.org/10.1016/S0304-3975(98)00072-3}.

\bibitem{Soskova2013}
Mariya I.~Soskova, The Turing Universe in the Context of Enumeration Reducibility. In: Bonizzoni P., Brattka V., Löwe B. (eds), \emph{The Nature of Computation. Logic, Algorithms, Applications. CiE 2013}. Lecture Notes in Computer Science, vol 7921. Springer, Berlin, Heidelberg, p.~371--382, \url{https://doi.org/10.1007/978-3-642-39053-1_44}

\bibitem{Soare2016}
Robert I.~Soare, \emph{Turing Computability. Theory and Applications}. Springer, 2016, ISBN 978-3-642-31932-7, \url{https://doi.org/10.1007/978-3-642-31933-4}

\bibquote{\nb{%
``Definition 1.7.5. (Acceptable Numbering Conditions). Let $\mathcal{P}$ be the class of partial computable funcitons of one variable. 

(i) A \emph{numbering} of a p.c. fuinctions is a map from $\omega$ onto $\mathcal{P}$.

(ii) The numbering $\{\varphi_e\}_{e\in\omega}$ of definition 1.5.1 is called the \emph{standard} numbering or \emph{canonical} numbering of the partial computable functions.

(iii) Let $\hat{\pi}$ be another numbering and let $\psi_e$ denote $\hat{\pi}(e)$. Then $\hat{\pi}$ is an \emph{acceptable numbering} if there are computable functions $f$ and $g$ such that (1)~$\varphi_{f(x)}=\psi_x$, and (2) $\psi_{g(x)}=\varphi_x$. $\langle\ldots\rangle$

Theorem 1.7.6 (Acceptable Numbering Theorem, Rogers). For any acceptable numbering $\{\psi_e\}_{e\in \omega}$ of the partial computable functions, there is a computable permutation $h$ of $\omega$ such that $\varphi_e = \psi_{h(e)}$ for all~$e$.''
}}

\bibitem{Wiki2018} Wikipedia pages: \emph{Primitive recursive function}, \url{https://en.wikipedia.org/wiki/Primitive_recursive_function} and \emph{$\mu$-recursive function}, \url{https://en.wikipedia.org/wiki/%CE%9C-recursive_function}. (Version of November 5, 2018)

\bibitem{Wolfram2018} 
Szudzik, Matthew. \emph{Recursive Function}. From MathWorld --- A Wolfram Web Resource, created by Eric W. Weisstein. \url{http://mathworld.wolfram.com/RecursiveFunction.html}. (Accessed November 5, 2018)

\end{thebibliography}
\end{document}

\title{Математические работы Владимира Андреевича Успенского: комментарии\footnote{Владимир Андреевич Успенский был моим учителем (и научным руководителем на старших курсах и в аспирантуре, как это формально называется на мехмате МГУ). В этом обзоре говорится прежде всего о его математических работах; надеюсь выразить свою благодарность за всё остальное в отдельном тексте.}}
\author{Александр Шень\footnote{LIRMM, University of Montpellier, CNRS, Montpellier, France. Грант RaCAF--ANR-15-CE40-0016}}

\date{}
\maketitle

\begin{abstract}
Мой учитель Владимир Андреевич Успенский (1930--2018) был одним из пионеров теории вычислений и в целом математической логики в СССР. В этом обзоре предпринимается попытка описать его математические работы и их роль в развитии теории алгоритмов и математической логики в СССР. (Его организационная деятельность и достижения в лингвистике выходят за рамки этого обзора.) 
\end{abstract}

\subsection*{Гармонические функции}

Первая (студенческая) работа Успенского~\cite{1949}  предлагает элементарное изложение основных свойств гармонических функций (с минимальным использованием сведений из анализа\emdash определением считается утверждение теоремы о среднем). Изложение основано на таком наблюдении: угол, под которым виден фиксированный отрезок из переменной точки, является гармонической функцией этой точки в том смысле, что для него верна теорема о среднем,\footnote{Это следует из того, что среднее направление из фиксированной точки на переменную точку окружности совпадает с направлением на её центр.} и функция эта является ступенчатой на окружности, проходящей через концы отрезка.
\medskip

\subsection*{Дипломная работа}
В дипломной работе Успенского~\cite{1952} излагается модель вычислений, предложенная А.\,Н.\,Кол\-мо\-го\-ро\-вым и известная теперь как <<машины Колмогорова\endash Успенского>>. Доказывается, что эта модель эквивалентна частично рекурсивным функциям, определённым с помощью подстановки, рекурсии и минимизации (то есть даёт тот же класс вычислимых функций). Кроме того, в рамках этой модели вводится понятие вычислимости относительно некоторой функции $f$: график функции~$f$ представляется в виде бесконечного графа (комплекса), доступного алгоритму вместе со входом [определение (A) на с.~64].  Это определение относительной вычислимости сравнивается с другими.  Для этого Успенский переформулирует определение Тьюринга\endash Поста (машины с оракулом~\cite{Turing1939,Post1944})  в терминах вычислимого протокола взаимодействия с оракулом [определение (T) на с.~63], и доказывает, что полученное определение эквивалентно определению с бесконечным графом, кодирующим оракул (A). Кроме того, в работе доказано, что эти определения относительной вычислимости эквивалентны определению в терминах замыкания относительно операций подстановки, рекурсии и минимизации [определение (R) на с.~64]. Колмогоров (научный руководитель) пишет в своём отзыве:

\begin{quote}
В ней [дипломной работе] подвергается более полному, чем до сих пор делалось, анализу само понятие алгоритмической вычислимости.

   (1) Автор приводит только одно, предлагавшееся до него, формально безукоризненное определение алгорифмической [в этом предложении Колмогоров пишет <<алгорифмической>> через <<ф>>, как это делал А.\,А.\,Марков] сводимости, которое он на с.~22 приписывает \emph{Б.\,А.\,Трахтенброту: функция $\gamma$ сводится к функции $\delta$, если $\gamma$ принадлежит рекурсивному замыканию $\delta$.} Автор показывает, что в действительности такая сводимость может быть всегда осуществлена \emph{очень простым каноническим образом при помощи раз навсегда заданных примитивно\dash рекурсивных функций $\tau(u)$ и $\omega(u)$ и зависящих от пары $\gamma,\delta$ примитивно рекурсивных функций $h(u,v,w)$ и $\varphi(m)$. См. об этом теорему на стр.~28. Это основной новый с чисто математической точки зрения результат работы.}

Определение сводимости по Трахтенброту нуждается в известном <<оправдании>> его соответствия интуитивной идее сводимости в смысле существования <<механического>> способа получения при любом $x$ значения $\gamma(x)$ \emph{в предположении}, что получение значений $\delta(x)$ сделано каким-то способом <<доступным>> для любого $x$. Общие контуры возможной формализации этой идеи были намечены Поустом [так Колмогоров пишет фамилию Поста (Emil Post)]. В дипломной работе полностью воспроизведён перевод соответствующего места статьи Поуста [имеется в виду статья~\cite{Post1944}]. \emph{Автор дипломной работы, по-видимому впервые, даёт соответствующее этой идее определение сводимости с полной отчётливостью} и показывает его эквивалентность определению Трахтенброта. Это тоже весьма существенное достижение автора дипломной работы.

(2) Кроме того в работе содержится хороший обзор различных предлагавшихся ранее определений \emph{алгоритмической вычислимости числовой функции} $y=\gamma(x)$. В центре изложения помещено определение, предложенное мною, интерес которого на мой взгляд убедительно аргументирован автором дипломной работы. Доказана равносильность этого определения прежде предлагавшимся. \emph{В известном смысле слова этот результат можно рассматривать как <<обоснование>> прежних определений, так как в моём определении становится особенно ясной основная идея алгоритмической вычислимости, которая отличается от вычислимости обыкновенным реальным счётным механизмом только неограниченным объёмом <<запоминающего устройства>> механизма.}
\end{quote}

Чтобы оценить содержание работы, важно представлять себе исторический контекст. Сейчас этот контекст почти забыт, и о нём  надо сказать несколько слов.

\subsubsection*{Частично рекурсивные функции}
Если спросить, что такое частично рекурсивная функция (partial recursive function), большинство современных специалистов ответят, что это функция, которая может быть получена из базовых функций (проекция, нулевая функция и функция прибавления единицы) с помощью операций подстановки, рекурсии и минимизации ($\mu$-оператора). Это определение можно найти в классической книге Одифредди~\cite[с.~127]{Odifreddi1989}, в других учебниках~\cite{Malcev1965,Cutland1980} и в википедии~\cite{Wiki2018}.

Но раньше определение было другим, и следы этого старого определения сохранились в другом классическом учебнике~\cite[раздел 1.5]{Rogers1972} и в cправочном ресурсе Wolfram MathWorld~\cite{Wolfram2018}. Хотя это определение и эквивалентно приведённому выше (задаёт тот же класс функций), но разницу между ними важно иметь в виду при чтении старых работ.

История вопроса здесь такова. Рекуррентные определения были известны давно (достаточно вспомнить о Фибоначчи), но их систематическое использование для построения арифметики появилось в работе Сколема 1923 года~\cite{Skolem1923}. Он понял, что таким образом можно определить не только простые функции, скажем, сложение и умножение (чтобы прибавить следующее за $y$ число к $x$, надо прибавить $y$ к $x$ и взять следующее число; чтобы умножить  $x$  на следующее за $y$ число, надо умножить $x$ на $y$ и прибавить $x$), но и много других функций, встречающихся в элементарной теории чисел. После этого базовые результаты этой теории можно доказать по индукции, исходя из рекурсивных определений.\footnote{Мотивацией этой работы было\emdash показать, что многие математические результаты можно обосновать простым и надёжным способом, индуктивно доказывая равенства между рекурсивно определёнными функциями. Теперь соответствующую теорию называют \emph{примитивно рекурсивной арифметикой}.}

Сколем не рассматривал явно класса всех функций, которые можно получить с помощью такого рода рекурсивных определений. Но уже в докладе 1925 года  Гильберт~\cite{Hilbert1926} говорит об определениях функций <<с помощью подстановок и рекурсий>> и различает <<обычные, пошаговые рекурсии>>, где значение функции на каком-то числе определяется через её значение на предыдущем числе, и более сложные схемы. В качестве примера такой более сложной схемы он приводит последовательность функций
\begin{multline*}
\varphi_1(a,b)=a+b, \ \varphi_2(a,b)=a\cdot b,  \ \varphi_3(a,b)=a^b,\\
\varphi_4(a,b)= \text{[$b$-й член в последовательности $a, a^a, a^{(a^a)}, a^{a^{(a^a)}}\ldots$]}
\end{multline*}
и так далее, которую можно задать рекурсивно равенствами
\[
\varphi_1(a,b)=a+b, \ \varphi_{n+1}(a,1)=a,  \ \varphi_{n+1}(a,b+1)=\varphi_n(a,\varphi_{n+1}(a,b)),
\]
и упоминает результат Аккермана о том, что функцию $\varphi_n(a,b)$ как функцию от трёх переменных $n,a,b$ нельзя задать <<обычными>> рекурсиями (этот результат опубликован позже в \cite{Ackermann1928}). Упоминание этого отрицательного результата Гильбертом означает, что у Гильберта было уже понятие о классе функций, которые можно получить <<обычными>> рекурсиями (и подстановками), хотя не было специального названия для функций этого класса и не было явно дано его определение.

Такое название и такое определение появились в знаменитой работе Гёделя~\cite[с.~179]{Godel1931}: функция называется рекурсивной  (rekursiv по-немецки), если она получается последовательным применением нескольких операций подстановки и операции рекурсии такого вида (построение функции $\varphi$, если $\psi$ и $\mu$ уже построены ранее):
\begin{align*}
\varphi(0,x_2,\ldots,x_n)&=\psi(x_2,\ldots,x_n)\\
\varphi(k+1,x_2,\ldots,x_n)&=\mu(k,\varphi(k,x_2,\ldots,x_n),x_2,\ldots,x_n)
\end{align*}
(схема (2) на с.~179). Гёдель использует представление этих функций в формальной системе для <<арифметизации>> утверждений о выводимости, так что для него они являются не предметом исследования, а средством.

Но как быть с более общими видами рекурсии (например, из упомянутого результата Аккермана)? Эрбран предложил (в письме Гёделю, а также в работе~\cite{Herbrand1932}) рассматривать системы функциональных уравнений (связывающих определяемые функции с построенными ранее), которые однозначно определяют новые функции. В его статье это формулируется так~\cite[с.~5, с.~624 английского перевода]{Herbrand1932}:
\begin{quote}
Мы можем также ввести произвольное количество функций $f_i(x_1,\ldots,x_{n_i})$ вместе с утверждениями о них, если:
\begin{description}
\item{(a)} эти утверждения не содержат связанных переменных;
\item{(б)} рассматриваемые с интуиционистской точки зрения, то есть как утверждения о натуральных числах, а не просто как символы, эти утверждения позволяют вычислить значение $f_i(x_1,\ldots,x_{n_i})$  для любого набора числовых аргументов, и можно интуиционистски доказать, что результат однозначно определён.
\end{description}
\end{quote}
Смысл этой оговорки про <<интуиционистскую точку зрения>>, видимо, в том, что нас не устроит само по себе функциональное уравнение, про которое из каких-то общих соображений можно доказать, что его решение существует и единственно; нужно, чтобы это доказательство было в каком-то смысле конструктивно и давало способ вычисления значений интересующих нас функций, исходя из задающих их равенств. (И действительно, впоследствии Кальмар~\cite{Kalmar1955} привёл пример системы функциональных уравнений, однозначно задающей невычислимую функцию.)

Гёдель возвращается к предложению Эрбрана (погибшего в горах сразу после отправки в редакцию статьи~\cite{Herbrand1932}) в своих лекциях в Принстоне (записки которых были размножены ещё тогда, а позднее перепечатаны, см.~\cite{Godel1934}). Он по-прежнему называет <<рекурсивными>> функции, которые получаются из базовых с помощью подстановок и <<обычных>> рекурсий, но в разделе 9  говорит о рекурсивных определениях более общего вида и задаваемых ими функциях, называя их ``general recursive functions''. Гёдель воспроизводит предложение Эрбрана так:
\begin{quote}
One may attempt to define this notion [general recursive function] as follows: if $\phi$ denotes an unknown function, and $\psi_1,\ldots,\psi_k$ are known functions, and if the $\psi$'s and the $\phi$ are substituted in one another in the most general fashions and certain parts of the resulting expressions are equated, then if the resulting set of functional equations has one and only one solution for $\phi$, $\phi$ is a recursive function.'' 
\end{quote}
(В примечании к этому определению Гёдель ссылается на письмо Эрбрана.) Далее он добавляет ограничения, уточняющие замысел Эрбрана:
\begin{quote}
We shall make two restrictions on Herbrand's definition. The first is that the left-hand side of each of the given functional equations defining $\phi$ shall be of the form \[\phi(\psi_{i1}(x_1,\ldots,x_n),\psi_{i2}(x_1,\ldots,x_n),\ldots,\psi_{il}(x_1,\ldots,x_n)).\] The second (as stated below) is equivalent to the condition that all possible sets of arguments $(n_1,\ldots,n_l)$ of $\phi$ can be so arranged that the computation of the value of $\phi$ for any given set of arguments $(n_1,\ldots,n_l)$ by means of the given equations requires a knowledge of the values of $\phi$ only for sets of arguments which precede $(n_1,\ldots,n_l)$.
\end{quote}
 Гёдель не уточняет порядок на наборах аргументов, так что точный смысл этого не очень ясен. Но дальше описываются конкретные правила вывода одних равенств из других и говорится:
  \begin{quote}
 Now our second restriction on Herbrand's definition of recursive function is that for each set of natural numbers $k_1,\ldots,k_l$ there should be one and only one $m$ such that $\phi(k_1,\ldots,k_l)=m$ is a derived equation.
\end{quote}
Тем самым даётся вполне точное определение некоторого класса функций, названных ``general recursive functions''; по-русски обычно переводят этот термин (немного загадочно) как <<общерекурсивные функции>>. Однако, как пишет Клини в~\cite{Kleene1981}, в момент чтения лекций (1934) сам Гёдель ещё не был уверен в том, что этот класс функций достаточно широк: <<However, Gödel, according to a letter he wrote to Martin Davis on 15 February 1965, ``was, at the time of [his 1934] lectures, not at all convinced that [this] concept of recursion
comprises all possible recursions''>>~\cite[p.~48]{Kleene1981}. В~\cite[p.~40]{Davis1965} про это говорится так:
\begin{quote}
In the present article [речь идёт о~\cite{Godel1934}] G\"odel shows how an idea of Herbrand's can be modified so as to give a general notion of recursive function $\langle\ldots\rangle$ G\"odel indicates (cf. footnote 3) that he believed that the class of functions obtainable by recursion of the most general kind were the same as those computable by a finite procedure. However, Dr.~G\"odel has stated in a letter that he was, at the time of these lectures, not at all convinced that his concept of recursion comprised all possible recursions; and that in fact the equivalence between his definition and Kleene's in Math. Ann.~112~[речь идёт о~\cite{Kleene1936}] is not quite trivial. So despite appearances to the contrary, footnote 3 of these lectures is not a statement of Church's thesis.
\end{quote}
Footnote 3~\cite[p.~44]{Davis1965} относится к утверждению о том, что всякая примитивно рекурсивная функция может быть вычислена с помощью конечной процедуры, и говорит ``The converse seems to be true, if, besides recursions according to the scheme (2) [примитивная рекурсия], recursions of other forms (e.g., with respect to two variables simultaneously) are admitted. This cannot be proved, since the notion of finite computation is not defined, but it serves as a heuristic principle''.

Роза Петер в~\cite{Peter1934} изучает возможности <<обычных рекурсий>> (например, доказывает, что разрешение использовать несколько значений функции в меньших точках сводится к схеме с одним предшественником) и вводит термин ``primitive Rekursion'' для этих самых <<обычных рекурсий>>.

Следуя ей, Клини в 1936 году~\cite{Kleene1936}  вводит термин ``primitive recursive functions'' (примитивно рекурсивные функции) для тех функций, которые получаются с помощью подстановок и примитивных рекурсий и которые Гёдель в~\cite{Godel1931} называл просто <<рекурсивными>>. Одновременно Клини предлагает рассмотреть более общий класс функций, элементы которого он называет ``general recursive functions'' (его статья так и называется, \emph{General recursive functions of natural numbers}). Этот класс определяется в духе Эрбрана и Гёделя, при этом рассматриваются разные правила вывода одних равенств из других, которые, однако (как доказывает Клини), задают один и тот же класс функций.

Клини также вводит $\eps$-оператор $\eps x [A(x)]$ как наименьшее число, удовлетворяющее условию $A(x)$, если таковое существует; в противном случае берётся нуль. Этот оператор фигурирует в теореме IV, которая утверждает, что всякая общерекурсивная функция может быть представлена в виде 
\(
\psi(\eps y [R(x,y)]),
\)
где $\psi$\emdash некоторая примитивно рекурсивная функция, а $R$\emdash примитивно рекурсивный предикат (это означает, что предикат~$R$ записывается как равенство нулю некоторой примитивно рекурсивной функции), причём для всякого $x$ существует $y$, при котором $R(x,y)$.\footnote{Та частьß определения $\eps$-оператора, где результат полагается равным нулю, когда искомого $y$ не существует, при этом роли не играет. Таким образом, здесь можно заменить $\eps$-оператор на стандартный $\mu$-оператор, в котором значение считается неопределённым в случае отсутствия искомого~$y$.} Теорема V утверждает, что верно и обратное: всякая функция, представимая в указанном виде, рекурсивна в смысле определений в духе Эрбрана и Гёделя. Тем самым такое представление может рассматривать как эквивалентное определение понятия рекурсивной функции. Кроме того, из этого можно извлечь некоторый способ нумерации всех рекурсивных функций, введя дополнительный параметр $e$ в примитивно рекурсивный предикат $R$ (хотя не при всех $e$ получается всюду определённая функция;  можно было бы сказать, что нумеруются частичные функции, но пока Клини их не рассматривает).

В том же (1936) году Чёрч публикует статью~\cite{Church1936}, в которой приводит другое определение некоторого класса числовых функций (в терминах так называемого $\lambda$-исчисления) как формализацию интуитивной идеи вычислимости:
\begin{quote}
The purpose of the present paper is to propose a definition of effective calculability${}^3$ which is thought to correspond satisfactorily to the somewhat intuitive notion.
\end{quote}
В подстрочном примечании $({}^3)$ Чёрч пишет: 
\begin{quote}
As will appear, this definition of effective calculability can be stated in either of two equivalent forms, (1) that a function of positive integers shall be called effectively calculable if it is $\lambda$-definable in the sense of \S2 below, (2) that a function of positive integers shall be called effectively calculable if it is recursive in the sense of \S4 below. The notion of $\lambda$-definability is due jointly to the present author and S.C.~Kleene $\langle\ldots\rangle$ The notion of recursiveness in the sense of \S4 is due jointly to Jacques Herbrand and Kurt G\"odel $\langle\ldots\rangle$ The proposal to identify these notions with the intuitive notion of effective calculability is first made in the present paper\ldots
\end{quote}
и добавляет (примечание в \S7): 
\begin{quote}
The question of the relationship between effective calculability and recursiveness (which it is here proposed to answer by identifying the two notions) was raised by G\"odel in conversation with the author. The corresponding question of the relationship between effective calculability and $\lambda$-definability had previously been proposed by the author independently.
\end{quote}

Видно, что Чёрч считает важным делом отождествление интуитивного понятия вычислимости с принадлежностью к точно определённому классу функций (для которого есть два эквивалентных определения). Это отождествление и назвали потом \emph{тезисом Чёрча}.

Почти в то же время Тьюринг публикует свою работу~\cite{Turing1937}, в которой он определяет машины, называемые теперь \emph{машинами Тьюринга} (сам Тьюринг употребляет термин $a$-machine, от слова `automatic'), и строит универсальную машину (которая может моделировать любую машину, будучи снабжена подходящей программой). В терминах этих машин Тьюринг определяет класс вычислимых действительных чисел (те, знаки которых могут вычисляться машиной) и даёт своё доказательство неразрешимости  Entscheidungsproblem (нет алгоритма, который распознаёт общезначимость формул языка первого порядка; ранее это\emdash для эквивалентных определений вычислимости\emdash доказали Гёдель и Клини, а также Чёрч, см. подробнее в~\cite[p.109]{Davis1965}). 

 В приложении (Appendix), добавленном 28 августа 1936 года, Тьюринг намечает доказательство эквивалентности вычислимости последовательности с помощью предложенных им машин и $\lambda$-определимости. Описывая этот результат во введении к работе, он пишет:
\begin{quote}
In a recent paper Alonzo Church has introduced an idea of ``effective calculability'', which is equivalent to my ``computability'', but is very differently defined. Church also reaches similar conclusions about the Entscheidungsproblem. The proof of equivalence between ``computability'' and ``effective calculability'' [имеется в виду $\lambda$-определимость] is outlined in an appendix to the present paper.
\end{quote}

Независимо от Тьюринга и почти одновременно с ним Пост публикует работу~\cite{Post1936}, где описывает своё определение <<финитного комбинаторного процесса>>, отличающееся от машин Тьюринга лишь техническими деталями, а также тем, что он не говорит о машине, а описывает, как ``problem solver or worker'' следует указаниям (the set of directions) определённого вида. Далее Пост замечает:
\begin{quote}
The writer expects the present formulation to turn out to be logically equivalent to recursiveness in the sensе of the G\"odel\endash Church development. Its purpose, however, is not only to present a system of a certain logical potency but also, in its restricted field, of psychological fidelity. In the latter sense wider and wider formulations are contemplated. On the other hand, our aim will be to show that all such are logically equivalent to formulation 1 [предложенный Постом вариант определения]. We offer this conclusion at the present moment as a \emph{working hypothesis}. And to our mind such is Church's identification of effective calculability with recursiveness. $\langle\ldots\rangle$ The success of the above program would, for us, change this hypothesis not so much to a definition or to an axiom but to a \emph{natural law}. 
\end{quote}
И добавляет в примечании:
\begin{quote}
Actually the work already done by Church and others carries this identification considerably beyond the working hypothesis stage. But to mask this identification under a definition hides the fact that a fundamental discovery in the limitations of the mathematizing power of Homo Sapiens has been made and blinds us to the need of its continual verification.\footnote{Сейчас, пожалуй, это фундаментальное открытие (fundamental discovery) уже практически утратило смысл: чтобы говорить о соответствии интуитивной идеи вычислимости (идеи алгоритма в неформальном смысле этого слова) и формально определённого класса алгоритмов, нужно, чтобы эта интуитивная идея была сформирована независимо от программистского опыта\emdash а где теперь найти людей, которые бы познакомились с идеей алгоритма, не имея уже программистского опыта?}
\end{quote}

Видно, что к 1936 году уже сложилась почти что современная картина: есть несколько эквивалентных (задающих один и тот же класс функций) определений вычислимости, есть понимание, что они отражают интуитивную идею алгоритма и вряд ли что-то упущено (и даже есть некоторые интуитивные объяснения, почему это так).

Но есть два отличия от современной картины: одно скорее терминологическое, а другое более принципиальное. Терминологическое состоит в том, что ни в одной из работ, говорящих о рекурсивных функциях, они не определяются как функции, получаемых с помощью подстановки, рекурсии и $\mu$-оператора, хотя все ингредиенты для доказательства эквивалентности этого определения другим есть и сама эквивалентность явно упоминается Клини в 1943 году~\cite[p.~53, Corollary]{Kleene1943}.

Во вторых, во всех этих работах говорится о \emph{всюду определённых} функциях (определённых на всех натуральных числах). \emph{Частичные} функции появляются чуть позже, в другой работе Клини~\cite{Kleene1938}, где строятся вычислимые системы обозначений для ординалов (тут без частичных функций уже не обойтись). Описав процесс вывода утверждения о значении функции из равенств в духе Эрбрана и Гёделя и предположив, что он даёт не более одного ответа для искомого значения, он пишет:
\begin{quote}
If we omit the requirement that the computation process always terminate, we obtain a more general class of functions, each function of which is defined over a subset (possibly null or total) of the $n$-tuples of natural numbers, and possesses the property of effectiveness when defined. These functions we call partial recursive.
\end{quote}
Таким образом впервые появляется понятие \emph{частично рекурсивной} функции (partial recur\-sive function).\footnote{Русский термин тут, пожалуй, ещё более странный, чем слово <<общерекурсивные>> для всюду определённых вычислимых функций: создаётся впечатление, что функция лишь отчасти рекурсивна. Английский термин лучше, потому что в нём слово `partial' относится к слову `function'.} Клини рассматривает (естественным образом обобщаемые на частичные функции) операции подстановки и рекурсии, а также определяет действие $\mu$-оператора для случая частичной функции: 
\[
\mu y [R(m,y)=0] = n
\]
для частичной функции $R$, если $R(m,n)$ определено и равно нулю, а все предыдущие значения $R(m,0),\ldots,R(m,n-1)$ определены и не равны нулю. (Очевидно, что такое $n$ единственно, если существует; если же нет, то определяемая с помощью $\mu$-оператора функция не определена на $m$.) Как отмечает Клини, все три операции (подстановка, рекурсия и $\mu$-оператор) не выводят из класса частично рекурсивных функций (определённых по Эрбрану и Гёделю). Он отмечает также, что для класса частично рекурсивных функций от любого числа ($n$) переменных существует универсальная функция $\Phi_n(z,\mathbf{x})$ от $n+1$ переменной. (Универсальность означает, что фиксацией первого аргумента $z$ из $\Phi_n$ можно получить любую частично рекурсивную функцию от $n$ переменных.) Эта универсальная функция может быть представлена в виде
\[
\Phi_n (z,\mathbf{x})= S(z,\mu y T_n(z,\mathbf{x},y)),
\]
где $S$\emdash некоторая примитивно рекурсивная функция, а $T_n$\emdash некоторый примитивно рекурсивный предикат (задаваемый условием обращения в нуль некоторой примитивно рекурсивной функции). В этом представлении, неформально говоря, $z$ кодирует (в виде натурального числа) систему равенств, задающих частично рекурсивную (в смысле Эрбрана и Гёделя) функцию от $n$ переменных, а $y$ является протоколом вывода из этой системы равенств утверждения о значении функции на входе $\mathbf{x}$. Предикат  $T_n$ проверяет корректность этого вывода, а функция $S$ извлекает из этого вывода установленное значение функции.\footnote{В формулировке Клини есть ещё первый аргумент $z$ у функции $S$, но его можно было бы и опустить.} Из этого результата (который называют <<теоремой Клини о нормальной форме>>) уже вытекает, что можно эквивалентно определить частично рекурсивные функции с помощью подстановки, рекурсии и $\mu$-оператора (и даже дополнительно потребовать, чтобы $\mu$-оператор применялся только один раз к примитивно рекурсивной функции). Но такой вариант определения Клини по-прежнему не упоминает.

Примерно эта же система понятий, терминология и способ изложения используются в более поздней статье Клини~\cite{Kleene1943}, посвящённой в основном арифметической иерархии, и в его же учебнике 1952 года~\cite{Kleene1957}, ставшем на многие десятилетия классическим. Отметим, помимо сказанного выше, ещё одну непривычную для нас вещь: формулировка <<тезиса Чёрча>> (отождествление интуитивной идеи вычислимой функции с точно определённым классом функций) относится только к всюду определённым функциям. 

\subsubsection*{Относительная вычислимость (вычислимость с оракулом)}

Можно определить относительную вычислимость одной функции относительно другой (или относительно некоторого множества, которое можно отождествить с его характеристической функцией). Впервые обсуждаемое определение сводимости предложил Тьюринг в своей диссертации (1939, см.~\cite{Turing1939}), но это было там побочной темой и рассматривалось лишь в некотором частном случае (сводимость к некоторому конкретному множеству). Он пишет:
\begin{quote}
Let us suppose that we supplied with some unspecified means of solving number-theoretic problems; a kind of oracle as it were. We will not go any further into the nature of this oracle than to say that it cannot be a machine.  With the help of the oracle we could form a new kind of machine (call them $o$-machines), having as one of its fundamental processes that of solving a given number-theoretic problem. More definitely these machines are to behave in this way. The moves of the machine are determined as usual by a table except in the case of moves from a certain internal configuration $\mathfrak{o}$. If the machine is in the internal configuration $\mathfrak{o}$ and if the sequence of symbols marked with $l$ is then the well formed formula \textbf{A}, then the machine goes into the internal $\mathfrak{p}$ or $\mathfrak{t}$ according as it is or is not true that \textbf{A} is dual. The decision as to which is the case is referred to the oracle.\par These machines may be described by tables of the same kind as used for the description of $a$-machines, there being no entries, however, for the internal configuration $\mathfrak{o}$.
\end{quote}
Общее определение сводимости по Тьюрингу предложил Пост в своей знаменитой статье~\cite[раздел 11]{Post1944}\emdash той самой, в которой он сформулировал \emph{проблему Поста} о существовании перечислимого неразрешимого множества, не являющегося полным (к которому сводятся не все перечислимые множества). Формально говоря, Пост в своём определении ограничивается сводимостью одного перечислимого множества к другому, но реально требование перечислимости в его тексте не используется. Определение Поста следует схеме Тьюринга и использует машины с оракулом. Клини в статье 1943 года~\cite{Kleene1943} предлагает другой вариант определения: к определению общерекурсивных функций с помощью выводов в исчислении равенств можно добавить равенства, выражающие таблицу значений для произвольных (всюду определённых) функций $\psi_1,\ldots,\psi_k$, и назвать те функции, которые можно определить таким способом, \emph{общерекурсивными относительно $\psi_1,\ldots,\psi_k$}:
\begin{quote}
A function $\phi$ which can be defined from given functions $\psi_1,\ldots,\psi_k$ by a series of applications of general recursive schemata we call \emph{general recursive} in the given functions; and in particular, a function $\phi$ definable ab initio by these means we call \emph{general recursive}.
\end{quote}
Но это описание дальше не развивается и не используется, оставаясь лишь пояснением к даваемому дальше определению, и для частично рекурсивных функций (в отличие от общерекурсивных) относительная вычислимость не рассматривается.  В книге 1952 года Клини говорит и о частичных функциях и их вычислимости и доказывает, что такое определение относительной вычислимости (с  выводами в исчислении равенств по Эрбрану и Гёделю) эквивалентно определению Тьюринга\endash Поста~\cite[\S 69]{Kleene1957}. При этом оракул  должен быть всюду определённой функцией (или набором таких функций). Но, в отличие от статьи Поста, где рассматривались лишь перечислимые множества в качестве оракулов, эта функция может быть произвольной всюду определённой функцией.

Обзор различных определений относительной вычислимости можно найти в~\cite[Section 4.3, ``History of Relative Computability'']{Soare1996}.

Теперь мы можем указать, в чём была новизна работы Успенского:\footnote{К сожалению,  дипломная работа не была опубликована, несмотря на рекомендации в отзывах, так что вряд ли повлияла на дальнейшее развитие событий, но тем не менее.}

\begin{itemize}

\item  впервые было явно указано (со ссылкой на Трахтенброта\emdash видимо, на неопубликованное сообщение) простое определение частичной рекурсивности, абсолютной и относительной, с помощью операций подстановки, рекурсии и минимизации ($\mu$-оператора);

\item была доказана (одновременно с книгой  Клини~\cite{Kleene1957} и гораздо более отчётливо) эквивалентность этого определения с другими определениями вычислимости (абсолютной и относительной);

\item было дано (впервые) определение относительной вычислимости, не связанное ни с какой конкретной вычислительной моделью, а использующее лишь класс вычислимых функций, и доказана эквивалентность этого определения другим определениям относительной вычислимости;

\item наконец, в дипломной работе Успенского была впервые изложена конструкция машин Колмогорова (позже опубликованная в совместной статье Колмогорова и Успенского~\cite{1958}), дано определение относительной вычислимости в терминах этой модели и доказана эквивалентность этого определения другим.
\end{itemize}

Третий пункт этого перечня требует некоторых пояснений. Определение Тьюринга\endash Поста для относительной вычислимости является модификацией соответствующего определения для абсолютной вычислимости: мы расширяем класс машин Тьюринга, дополнительно разрешая получать ответы от оракула. Аналогичным образом определение Клини для относительной вычислимости модифицирует определение частично рекурсивной функции. Таким образом, даже если мы договорились, какие функции мы считаем (абсолютно) вычислимыми, нам нельзя ещё забыть про конкретную модель вычислений и не возвращаться к ней: определяя относительную вычислимость, надо снова вспомнить о модели и её модифицировать. В отличие от этой ситуации, в определении из работы Успенского относительная вычислимость определяется с помощью алгоритмов диалога с оракулом, то есть требуется вычислимость некоторых (частичных) функций, задающих этот диалог (функций, указывающих следующий вопрос к оракулу при известном входе, предыдущих вопросах и ответах на них). Таким образом, для понимания смысла этого определения достаточно знать, какие функции считаются вычислимыми. 

Впоследствии такого рода (<<машинно\dash независимые>>) определения относительной вычислимости были даны другими авторами. Например, одно из них можно найти в учебнике Роджерса~\cite[раздел 9.2]{Rogers1972} (без ссылки на каких бы то ни было предшественников). Можно ещё отметить, что определение Успенского имеет то преимущество, что оно (в отличие от определения из книги Роджерса) естественно распространяется на произвольные частичные оракулы и при этом сохраняется эквивалентность с принадлежностью частично рекурсивному замыканию функции, использованной как оракул. 
%\nb{[Где это доказано? Я спросил у Сламана, но он пока не ответил. Видимо, это есть в Одифредди, но надо разбираться, кто это первый доказал]} 
Но сам Успенский этого не делает, рассматривая только функции с перечислимой областью определения (хотя доказательство для общего случая не требует существенных изменений).

\subsection*{Теорема Гёделя и теория алгоритмов}

Теорема Гёделя и теория вычислимых функций появились не только одновременно, но и вместе, как сиамские близнецы. В классической статье Гёделя, где доказан его результат о неполноте формальных теорий, одновременно было введено понятие рекурсивной функции (то, что теперь называется <<примитивно рекурсивными>> функциями, см. выше), и это понятие было важным техническим средством в доказательстве. А именно, различные функции, связанные с кодированием формул и выводов натуральными числами (их <<гёделевыми номерами>>, как раньше говорили), были определены рекурсивно, и это определение использовалось для погружения рассуждений о выводах в формальную систему.

С другой стороны, определение общерекурсивных функций было дано в терминах формальной системы (исчисления равенств), восходящей к Эрбрану и Гёделю.

Можно, пожалуй, сказать, что одним из первых достижений и в области теории вычислений, и в области теории доказательств, было разделение этих <<сиамских близнецов>>, и было это не таким простым делом, как сейчас кажется.   Сначала Тьюринг и Пост предложили модель вычислений (машины Тьюринга\endash Поста), позволившую определить вычислимость безо всякого упоминания формальных теорий и выводов в них. Общая природа теоремы Гёделя и её связи с теорией алгоритмов были осознаны в 1940-е годы, видимо, в первую очередь Клини и Колмогоровым.

 В статье Клини 1943 года~\cite{Kleene1943} была указано, что теорема Гёделя по существу означает  неперечислимость множества истинных формул, а в статье 1950 года~\cite{Kleene1950} аналогичная интерпретация была дана для теоремы Гёделя в форме Россера и указано, что она соответствует построению пары эффективно неотделимых перечислимых множества. Но хотя по существу все необходимые наблюдения были уже сделаны, по форме изложение Клини и в этой статье 1950 года, и в классическом учебнике 1952 года~\cite{Kleene1957} остаётся ещё тесно связанным с языком теории примитивно рекурсивных функций (достаточно сказать, что изложение  в~\cite{Kleene1950} начинается словами ``Let $T_1$ be the primitive recursive predicate so designated in a previous paper by the author''), и процедура погружения неотделимых множеств в формальную теорию скорее подразумевается, чем явно описана.

Примерно в то же время, и, вероятно, независимо к пониманию соотношения между теоремой Гёделя и теорией алгоритмов пришёл Колмогоров. Как рассказывает Успенский в~\cite[с.~323]{2006a},
 
\begin{quote}
2 декабря 1952 г. Колмогоров изложил мне весьма кратко, в течение пяти минут,\emdash но зато дал списать с заготовленной им бумажки, озаглавленной <<Гёдель и рекурсивная перечислимость>>,\emdash основополагающие идеи о связи теоремы Гёделя о неполноте аксиоматических систем (для самых общих исчислений) с существованием множеств, не являющихся рекурсивными, и пар множеств, не являющихся рекурсивно отделимыми. Бумажка была написана им <<для себя>>, и разобраться в ней, а тем более в его сопутствующих устных комментариях, мне было тогда непросто. Потом всё как-то выстроилось, и 8 мая 1953 г. Колмогоров представил в <<Доклады АН СССР>> мою заметку <<Теорема Гёделя и теория алгоритмов>>, написанную на основе его идей. Высокое искусство Колмогорова как учителя состояло в умении создать у ученика впечатление, что именно он, ученик, полноценный автор статьи. Колмогоров во много раз реже, чем имел на это все права, выступал в роли соавтора своих учеников $\langle\ldots\rangle$ В 1958 г. в <<Успехах математических наук>> под двумя нашими фамилиями вышла статья <<К определению алгоритма>>, в которой мне принадлежит, по существу, лишь черновая работа.
\end{quote}
(Статьи, о которых идёт речь:~\cite{1953, 1958}. Вторая из них содержит изложение вычислительной модели с преобразованием графов, которая фигурировала уже в дипломной работе Успенского.)

В статье Успенского 1953 года~\cite{1953} было отчётливо указано, уже без всякого упоминания о примитивно рекурсивных функциях, что теорема Гёделя о том, что достаточно богатая формальная система (скажем, формальная арифметика) неполна и не может быть пополнена, следует из того, что существуют перечислимые неотделимые множества и что эта пара множеств погружается в формальную систему\emdash как сейчас сказали бы, сводится к паре (доказуемые формулы, опровержимые формулы). А эффективная непополнимость (тот факт, что по расширению формальной системы дополнительными аксиомами можно алгоритмически указать формулу, которая остаётся недоказуемой и неопровержимой) следует из существования эффективно неотделимых перечислимых множеств. Но, повторим ещё раз, всё это по существу уже было в работе Клини~\cite{Kleene1950}, о которой Колмогоров и Успенский, судя по всему, тогда не знали. Успенский ссылается на работу Клини 1943 года~\cite{Kleene1943}, говоря о рекурсивных функциях, но говоря о существовании перечислимых неотделимых множеств, не ссылается на~\cite{Kleene1950}, где они построены, а  говорит лишь, что они были построены Новиковым (не указывая никакой публикации Новикова, а лишь давая ссылку на работу Трахтенброта 1953 года).

Говоря о ситуации в целом, можно сказать, что есть два взаимно дополнительных взгляда на теорему Гёделя. С одной стороны, она является реализацией парадокса лжеца (в одном из вариантов он говорит, что утверждение <<Это утверждение ложно>> не может быть ни ложным, ни истинным): если вместо этого сделать утверждение <<Это утверждение недоказуемо>>, то оно будет истинным (и потому неопровержимым), но недоказуемым. Это объяснение не ссылается на теорию алгоритмов, хотя для обоснования возможности записать неформальные рассуждения в формальной арифметике можно, следуя Гёделю, использовать примитивно рекурсивные функции как техническое средство. С другой стороны, теорема Гёделя является следствием существования неразрешимых перечислимых множеств (или, в симметричном варианте, существования неотделимых перечислимых множеств), и в таком изложении никакой <<самоприменимости>> не заметно. Но, конечно, она никуда не делась, переместившись в конструкцию неразрешимого перечислимого множества (или неотделимых множеств). Эта конструкция следует идее <<диагонального аргумента>> Кантора, которая, в свою очередь, является формой проявления самоприменимости (<<диагональ>> состоит из результатов применения вычислимой функции к своему номеру, или алгоритма к его собственному тексту). 

Впоследствии Успенский опубликовал подробное изложение доказательства теоремы Гёделя с помощью средств теории алгоритмов (а также изложение начал этой теории) сначала в виде статьи~\cite{1974}, а затем (в расширенном виде) брошюры~\cite{1982} в серии <<Популярные лекции по математике>>.  Это изложение до сих пор остаётся, пожалуй, наиболее доступным и корректным изложением теоремы Гёделя для неспециалистов (по крайней мере если говорить о её алгоритмическом аспекте).

Помимо этого, в~\cite{1974,1982} была намечена (совсем не очевидная в то время, достаточно сравнить с тем же учебником Роджерса~\cite{Rogers1972}) схема изложения теории алгоритмов. Традиционно (в том числе и в книге самого Успенского~\cite{1960})  изложение теории алгоритмов начиналось с подробного разбора какой-то вычислительной модели (сначала в этом качестве были популярны рекурсивные функции, потом машины Тьюринга; Марков использовал для этого нормальные алгорифмы), и это занимало достаточно много места и времени. Лишь после этого оставшиеся слушатели (читатели) знакомились с простейшими фактами вроде теоремы Поста (перечислимое множество с перечислимым дополнением разрешимо), и т.п. Конечно, можно было пропустить первую часть, с построением конкретной вычислительной модели, и начинать прямо со второй, рассуждая как в анекдоте о беспроволочном телеграфе (<<представьте себе длинную кошку, которую в одном городе дёргают за хвост, а в другом она мяукает\emdash это проволочный телеграф,\emdash а теперь то же самое, но без кошки>>), тогда рассуждения становились простыми и наглядными, но беспочвенными.

Выход, предложенный Успенским в~\cite{1974} (а до этого использованный в его лекциях 1972/3 года, но, вероятно, и в предыдущие годы), состоял в следующем: мы рассуждаем о классе вычислимых функций, приняв на веру (в качестве <<аксиом>>) некоторые свойства этого класса, отчётливо сформулированные, но оставленные без доказательства. Помимо вычислимости конкретных функций (а также сохранения вычислимости при конкретных построениях одних функций из других), Успенский выделяет два таких свойства, называя их <<аксиомой протокола>> и <<аксиомой программы>>. 

\emph{Аксиома протокола} утверждает, что для всякого алгоритма $A$ существует разрешимое множество $R$, элементы которого называются <<протоколами>>, и две вычислимые функции $\alpha$ и $\omega$. Неформально говоря, элементы $R$ являются протоколами  (программисты сказали бы: <<логами>>) вычисления алгоритма $A$, то есть записями всех последовательных шагов его работы на некотором входе, в тех случаях, когда эта работа завершается и даёт какой-то результат. Функция $\alpha$ выделяет из протокола исходное данное (вход), а $\omega$\emdash результат работы (выход). Формально же требуется выполнение такого свойства: выход алгоритма $A$ на входе $x$ равен $y$ тогда и только тогда, когда существует $r\in R$, для которого $\alpha(r)=x$ и $\omega(r)=y$.

\emph{Аксиома программы} утверждает, что есть некоторое разрешимое множество $P$, элементы которого называются <<программами>>, и алгоритм $U$ применения произвольной программы $p\in P$ к произвольному входу $x$ (таким образом, входом алгоритма $U$ является пара $\langle p,x\rangle$). При этом любая вычислимая функция $f$ задаётся некоторой программой $p$ в том смысле, что $U(p,x)=f(x)$ для всех $x$. Последнее равенство понимается так: обе его части одновременно определены или одновременно не определены, и равны в том случае, когда определены.

Приняв эти аксиомы, можно развивать теорию алгоритмов, не вдаваясь в технические детали модели вычислений. Вместе с тем остаётся совершенно понятно, чего недостаёт в этих рассуждениях: мы должны указать конкретную модель вычисления, научиться в ней программировать те конструкции, которые использованы в доказательствах (и которые не так уж и просты, достаточно вспомнить, скажем, метод приоритета), а также проверить выполнение аксиом протокола и программы для этой модели вычислений. Конечно, и после этого некоторый психологический барьер остаётся (многим людям, которые легко ориентируются в достаточно сложных математических конструкциях, теория вычислимости всё же кажется чем-то странным), но по крайней мере он становится более явным и отчётливо видным.\footnote{Сейчас ситуация с методической точки зрения изменилась, прежде всего потому, что большинство приступающих к изучению теории вычислимых функций уже имеют программистский опыт. Возможно, современная реализация педагогических идей Успенского состояла бы в том, что мы предлагаем слушателям представить себе знакомый им язык программирования, к которому добавлены библиотечные функции интерпретатора этого же языка (аргументами которого являются строка, понимаемая как программа, и вход; это соответствует аксиоме программы), а также пошагового отладчика (который комбинирует аксиому программы с аксиомой протокола; на вход ему подаются текст программы, вход и число шагов работы отлаживаемой программы). Если к этому добавить ещё и функцию (без аргументов), выдающую текст текущей исполняемой программы, это облегчит и доказательство теоремы о неподвижной точке, сделав её самоочевидной.}

Для доказательства теоремы Гёделя нужна и третья аксиома, не следующая из этих двух\emdash что всякую вычислимую функцию можно выразить арифметической формулой (и тут снова приходится обращаться к конкретной модели вычислений).

Заметим, что при таком <<машинно-независимом>> изложении теории вычислимости мы не имеем права вновь возвращаться к вычислительной модели, рассуждая, скажем, о преобразованиях программ или о вычислениях с оракулом, а должны давать все необходимые определения, ссылаясь только на понятие вычислимой функции. Как мы уже говорили, такое определение для относительной вычислимости было (впервые) дано Успенским в дипломной работе, затем это было сделано для сводимости по перечислимости, а также для <<способов программирования>> (формализацией которых стало введённое Успенским понятие главной нумерации). Об этих двух последних достижениях Успенского мы говорим в следующем разделе.

Ещё можно отметить, что эта аксиоматизация теории алгоритмов позволяет строго обосновать известное наблюдение о том, что большинство результатов теории алгоритмов <<релятивизуются>>, то есть сохраняют силу, если вычислимые функции заменить функциями, вычислимыми с некоторым фиксированным оракулом (в качестве которого можно взять множество или всюду определённую функцию). В самом деле, достаточно проверить, что для этого релятивизованного класса (функций, вычислимых с данным оракулом) выполнены все аксиомы теории алгоритмов (кроме свойства арифметичности, естественно), и потому выполнены и все теоремы, выводимые из этих аксиом. Успенский поставил вопрос, полностью ли объясняет это наблюдение возможность релятивизации, то есть верно ли, что всякое утверждение, выполненное для всех классов вычислимых с некоторым оракулом функций, является следствием указанных им аксиом. Оказалось, что когда этот вопрос поставлен, получить на него (положительный) ответ уже несложно~\cite{Shen1980}.

\subsection*{Вычислимые отображения множеств и\\ сводимость по перечислимости}

Понятие сводимости, введённое Тьюрингом и Постом и рассмотренное в дипломной работе Успенского (см. выше), можно назвать <<сводимостью по разрешимости>>: сводимость множества $A$ к множеству $B$ гарантирует, что если $B$ разрешимо, то и $A$ разрешимо. Можно сказать, что в этом определении мы <<сводим задачу разрешения множества $A$ к задаче разрешения множества $B$>>. В работе~\cite{1955} Успенский предлагает определение \emph{сводимости по перечислимости}, в котором речь идёт о сведении задачи перечисления множества $A$ к задаче перечисления другого множества $B$. Это определение использует введённой в этой же работе понятие \emph{вычислимой операции} над множествами. 

В простейшем случае (одноместная операция, аргументами и значениями которой являются подмножества натурального ряда) вычислимые операции определяются так. Введём на множестве $\mathcal{P}(\mathbb{N})$ всех подмножеств натурального ряда топологию. Для каждого конечного множества $X\subset \mathbb{N}$ рассмотрим семейство $\mathcal{O}(X)$  всех его надмножеств, и будем считать открытыми в $\mathcal{P}(\mathbb{N})$ все такие семейства и все их объединения. Теперь рассмотрим все непрерывные в смысле этой топологии отображения $F\colon \mathcal{P}(\mathbb{N})\to\mathcal{P}(\mathbb{N})$. Легко проверить, что все такие отображения монотонны (если $U\subset V$, то $F(U)\subset F(V)$), и значение $F$ на любом множестве $U$ определяется значениями $F$ на конечных подмножествах $U$ (надо объединить $F(X)$ для всех конечных $X\subset U$). Значения $F$ на конечных множествах можно описать множеством пар $\{\langle n,X\rangle\mid n\in F(X)\}$ (здесь $n$\emdash натуральное число, а $X$\emdash конечное множество натуральных чисел). Непрерывное отображение $F\colon \mathcal{P}(\mathbb{N})\to\mathcal{P}(\mathbb{N})$ Успенский называет \emph{вычислимой операцией}, если соответствующее ему множество пар перечислимо (заметим, что эти пары являются конечными объектами, так что можно говорить о перечислимости их множества). После этого определяется сводимость по перечислимости: множество $A$ \emph{сводится по перечислимости} к множеству $B$, если существует вычислимая операция $F$, переводящая $B$ к $A$. Отмечается, что сводимость по Тьюрингу можно определить в терминах сводимости по перечислимости: всюду определённая функция $\varphi$ вычислима по Тьюрингу с оракулом для всюду определённой функции $\psi$ тогда и только тогда, когда график $\varphi$ сводится по перечислимости к графику $\psi$. Отсюда же получается и критерий сводимости множеств (переходом к их характеристическим функциям). Указывается, что в этих же терминах можно получить определить и понятие частично рекурсивного оператора в смысле Клини (\cite{Kleene1957}, см. обсуждение ниже).

Наконец, в этой же работе указывается эквивалентность предложенного определения вычислимой операции с двумя <<машинно\dash зависимыми>> определениями. Соответствующие понятия Успенский называет <<операциями Колмогорова>> и <<операциями Поста>> (хотя в явном виде они в публикациях Колмогорова и Поста не встречаются).

В другой работе 1955 года (\cite{1955a}, см. также изложение результатов этой работы с некоторыми добавлениями в~\cite{1956}) вводится (со ссылкой на доклад Колмогорова на семинаре по рекурсивной арифметике в 1954 году) понятие нумерации, определяется понятие главной нумерации системы перечислимых множеств и устанавливается связь вычислимых операций на перечислимых множествах в смысле~\cite{1955} с алгоритмическими преобразованиями номеров. 

Более подробно. Пусть мы хотим говорить об алгоритмических преобразованиях \emph{программ} вычислимых функций (или перечислимых множеств). Тогда нам мало знать, какие функции вычислимы (какие множества перечислимы), но нужно ещё и оговорить класс <<способов программирования>>, используемых для записи их программ. Программы являются обычно словами в некотором алфавите, но их можно отождествить с натуральными числами (при какой-то естественной нумерации слов). Тогда способ программирования вычислимых функций превращается в универсальную функцию двух аргументов: $U(n,x)$ есть результат применения программы с номером $n$ в входу $x$, который мы тоже считаем натуральным числом. Способ программирования перечислимых множеств тогда становится универсальным множеством пар $\langle n,x\rangle$, для которых $x$ принадлежит перечислимому множеству, программа которого имеет номер~$n$. На другом (но эквивалентном) языке можно сказать, что способ программирования вычислимых функций (перечислимых множеств) представляет собой \emph{натуральную нумерацию} множества вычислимы функций (перечислимых множеств), то есть отображение всего натурального ряда на множество вычислимых функций (соответственно перечислимых множеств): числу $n$ ставится в соответствие вычислимая функция (перечислимое множество), программа которого имеет номер $n$.

Если не накладывать на способы программирования (нумерации) никаких ограничений, то оказывается, что не все они одинаково хороши. Разумная теория, описывающая алгоритмические преобразования программ, требует некоторых дополнительных ограничений. По существу они встречались в классических работах Клини ($s$-$m$-$n$-теорема), но явно впервые были сформулированы в~\cite{1955a}. А именно, требуется, чтобы нумерация была \emph{главной}.  Определение главной нумерации включает в себя два требования. Во-первых, нумерация должна быть вычислимой: это значит, что соответствующая универсальная функция является вычислимой (частичной) функцией двух аргументов (вариант для множеств: соответствующее универсальное множество пар должно быть перечислимым множеством пар). Во\dash вторых, к этой нумерации должна \emph{сводиться} любая другая вычислимая нумерация.\footnote{Определение сводимости нумераций также опубликовано в~\cite{1955a} со ссылкой на доклад Колмогорова (видимо, впервые).} Это требование означает, что для любой другой вычислимой нумерации того же множества существует (всюду определённая) вычислимая функция, преобразующая номера в этой второй нумерации в номера в первой (главной).

Используя понятие главной нумерации для перечислимых множеств, можно рассмотреть вычислимые отображения семейства перечислимых множеств в себя. \emph{Вычислимость} означает, что есть алгоритм, который по (любой) программе перечислимого множества даёт (какую-то) программу его образа. Другими словами, мы рассматриваем всюду определённые преобразования программ, которых сохраняют эквивалентность (эквивалентные программы, то есть программы, задающие одно и то же множество, преобразуются в эквивалентные программы).  Успенский доказывает~\cite[раздел 6]{1955a}, что вычислимые отображения семейства перечислимых множеств в себя в точности представляют собой ограничения вычислимых операций (на семействе всех множеств) на класс перечислимых множеств. Аналогичное утверждение делается и для подсемейства униформных множеств пар (графиков функций): вычислимые отображения семейства вычислимых функций в себя представляют собой ограничения вычислимых операций, отображающих семейство униформных множеств в себя.

Опишем связи этих работ Успенского с работами других авторов того же периода.\footnote{К сожалению (см. ниже отрывок из воспоминаний Успенского), все три публикации Успенского~\cite{1955,1955a,1956} представляют собой краткие заметки в Докладах Академии наук СССР (первые две) и резюме доклада в Московском математическом обществе (третья), где приводятся только формулировки теорем и лемм, используемых в их доказательствах. Сами доказательства были опубликованы в кандидатской диссертации Успенского~\cite{1955b}, которая хотя формально и была доступна (её можно было заказать и получить в нескольких библиотеках в СССР), но вряд ли повлияла на дальнейшее развитие области. Да и статьи с кратким изложением~\cite{1955,1955a,1956}, видимо, остались неизвестными вне СССР. Позже книга Успенского~\cite{1960} (учебник по теории вычислимых функций, который стал докторской диссертацией Успенского) была переведена на французский язык\emdash к сожалению, в неё вошло лишь определение главной нумерации, но не результаты о вычислимых операциях и отображениях.} Райс~\cite{Rice1953} рассмотрел \emph{вполне перечислимые классы} перечислимых множеств, то есть такие классы перечислимых множеств, для которых множество \emph{всех} программ всех множеств этого класса перечислимо, и сформулировал гипотезу о том, что всякий такой класс состоит из всех надмножеств множеств из некоторого перечислимого семейства конечных множеств.  Эта гипотеза приводится в~\cite{1955a} как теорема 5 и является ключевым шагом в доказательстве упомянутых результатов о вычислимых отображениях. Она также доказана в статье самого Райса 1956 года~\cite{Rice1956}, где говорится, что независимо тот же результат получили МакНотон, Майхилл и Шапиро (ссылки на их работы не приводятся, кроме ссылки на краткую заметку Майхилла~\cite{Myhill1955}). Кроме того, уже в первой статье Райса~\cite{Rice1953} доказано, что никакое нетривиальное свойство перечислимых множеств нельзя алгоритмически распознать по его номерам (обобщение этого результата сформулировано как следствие теоремы 5 в~\cite{1955a}). Это утверждение поэтому в англоязычной литературе называют обычно <<теоремой Райса>>, а утверждение о строении вполне перечислимых классов перечислимых множеств (и аналогичный результат о вполне перечислимых классах вычислимых функций) называют <<теоремой Райса\endash Шапиро>> (см., например, \cite[глава 7, \S 2]{Cutland1980}).  Результат о связи вычислимых отображений класса вычислимых функций в себя и вычислимых операций на классе функций (под названием ``partial recursive functionals'') был доказан также Майхиллом и Шепердсоном~\cite{MyhillSheperdson1955} в том же 1955 году, когда были опубликованы статьи Успенского~\cite{1955,1955a} и потому в англоязычной литературе называется <<теоремой Майхилла\endash Шепердсона>> (см., например, \cite[глава 10, \S 2]{Cutland1980}. (Его частным случаем является теорема Райса\endash Шапиро, поэтому она иногда тоже называется теоремой Майхилла\endash Шепердсона, см., например, \cite[Theorem II.4.2 или Proposition II.5.19]{Odifreddi1989}.)

Трудно сказать, как было переоткрыто понятие сводимости по перечислимости; в книге Роджерса 1967 года~\cite{Rogers1972} соответствующее определение (``enumeration reducibility'') приводится без каких бы то ни было ссылок на Успенского или кого-либо ещё. В статье~\cite{Case1971} (``Enumeration reducibility and partial degrees'', 1971 года) даны ссылки на книгу Роджерса и статью Майхилла~\cite{Myhill1961}, но в статье Майхилла (как и в книге Дэвиса~\cite{Davis1958}, на которую он ссылается) сводимость по перечислимости не определяется, а рассматриваются различные варианты относительной вычислимости для функций. В современном обзоре~\cite{Soskova2013} работы Успенского вообще не упоминаются, а даётся ссылка на работу Фридберга и Роджерса 1959 года~\cite{FriedbergRogers1959}, которая в свою очередь ссылается на записки лекций Роджерса в MIT 1955--56 годов, опубликованные (размноженные) в 1957 году, из которых потом вышла книга~\cite{Rogers1972}. Видимо, Роджерс чуть позже Успенского независимо пришёл к тому же понятию (и с тем же названием).

Понятие главной нумерации также было переоткрыто Роджерсом (см.~\cite{Rogers1958}). Роджерс называет это понятие ``G\"odel numbering''. Сначала он даёт <<машинно\dash зависимое>> определение:  ``A G\"odel numbering is a numbering equivalent to the standard numbering'' (с.~333), но потом приводит и машинно\dash независимую характеризацию (как максимальый элемент относительно сводимости, как у Успенского\emdash хотя и без ссылки на него). В современной литературе иногда употребляется также термин ``admissible numbering'', см., например, недавнюю книгу Соара~\cite{Soare2016}, а также ``acceptable numbering''~\cite[Definition II.5.2]{Odifreddi1989}; в обоих случаях приводится <<машинно\dash зависимое>> определение.

Сравнивая работы Успенского с другими публикациями на близкие темы, нужно иметь в виду, что с определениями сводимости (относительной вычислимости) для частичных функций имеется путаница, как терминологическая, так и по существу. Есть три разных понятия вычислимости одной частичной функции относительно другой. Пусть $f$ и $g$\emdash две частичные функции (с натуральными аргументами и значениями). Сводимость $f$ к $g$ (вычислимость $f$ относительно $g$) можно понимать в трёх смыслах (каждый следующий сильнее предыдущего).

\begin{enumerate}

\item График $f$ сводится по перечислимости к графику $g$.

\item Рассмотрим, следуя Успенскому, семейство $\mathfrak{U}$ всех частичных функций из $\mathbb{N}$ в $\mathbb{N}$ с топологией, в которой базовыми открытыми множествами являются семейства всех продолжений некоторой конечной функции. Непрерывное отображение $F\colon\mathfrak{U}\to\mathfrak{U}$ мы будем называть вычислимой операцией, если его ограничение на конечные функции имеет перечислимый график, то есть если множество всех пар $\langle \langle x,y\rangle, u\rangle$, где $x$ и $y$\emdash натуральные числа, а $u$\emdash конечная частичная функция из $\mathbb{N}$ в $\mathbb{N}$, и при этом $[F(u)](x)=y$, перечислимо. Теперь вычислимость $f$ относительно $g$ можно понимать как существование вычислимой операции, переводящей $g$ в $f$.

\item Можно распространить определение Трахтенброта (см. обсуждение дипломной работы выше) на частичные функции и говорить, что частичная функция $f$ вычислима относительно частичной функции $g$, если $f$ принадлежит частично рекурсивному замыканию множества частично рекурсивных функций, к которому добавлена функция~$g$. (Это определение, например, приводится в~\cite{Malcev1965}.)
\end{enumerate}

Третье определение имеет эквивалентную переформулировку в терминах машины с оракулом. Эта переформулировка по существу повторяет определение из дипломной работы Успенского, но для частичных функций. А именно, значение $f(x)$ вычисляется алгоритмом, который получает на вход $x$ и может запрашивать значения функции $g$ в произвольных точках, но как только одно из запрошенных значений $g$ окажется неопределённым, вычисление <<зависает>> (прерывается без результата) и $f(x)$ остаётся неопределённым. Второй вариант определения можно тоже переформулировать в терминах машин, дополнительно разрешив параллельные запросы нескольких значений функции $g$; каждый из этих запросов не останавливает вычисление, которое продолжается и получает информацию о запрошенных значениях, если они определены, через какое-то время. При этом требуется, чтобы результат вычисления не зависел от того, через какое время поступят запрошенные значения.

Эти три определения различаются: каждое следующее сильнее предыдущего (более ограничительно). Разницу между этими определениями можно пояснить двумя примерами. Первый пример, разделяющий первый и второй варианты определения, таков. Пусть $f$\emdash произвольная всюду определённая функция, а $g$\emdash функция, принимающая только нулевые значения, и область определения функции $g$ состоит из всех номеров пар $\langle n, f(n)\rangle$ (для какой-то вычислимой нумерации пар). Тогда $f$ сводится к $g$ в смысле первого определения, но не обязательно сводится в смысле второго. (Это рассуждение приведено в примечании Успенского на с.~362 русского перевода книги Роджерса~\cite{Rogers1972} со ссылкой на Д.\,Г.\,Скордева; приведённое Роджерсом рассуждение существенно сложнее.)

Второй пример~\cite[Proposition II.3.20, ссылка на диссертацию Sasso 1971 года]{Odifreddi1989} показывает разницу между вторым и третьим определениями. Пусть $g$\emdash произвольная частичная функция натурального аргумента, принимающая только нулевые значения. Построим другую частичную функцию $f$, которая тоже принимает только нулевые значения, при этом $f(n)$ определено (и равно нулю) тогда и только тогда, когда хотя бы одно из значений $g(2n)$ и $g(2n+1)$ определено. Тогда функция $f$ вычислима относительно $g$ в смысле второго определения, но не обязательно вычислима в смысле третьего. (В терминах машин: если разрешено параллельно запрашивать $g(2n)$ и $g(2n+1)$, ожидая, пока один из этих запросов будет удовлетворён, то вычислить $f$ легко, но если можно лишь запрашивать их последовательно, то ничего не выйдет, потому что неудовлетворённый первый запрос помешает перейти ко второму. Конечно, это лишь неформальное пояснение, для доказательства различия нужно строить соответствующий пример функции $g$ диагональным методом, и это легко сделать.)

Первое определение соответствует тому, что в книге Роджерса~\cite[\S 9.8]{Rogers1972} названо <<частичнорекурсивными операторами>> (partial recursive operators). Второе соответствует тому, что названо там же <<рекурсивными операторами>> (recursive operators).

Майхилл и Шепердсон~\cite{MyhillSheperdson1955} говорят о ``partial recursive functionals'', ссылаясь на Thesis I$^{*\dagger}$ на с.~332 книги Клини~\cite{Kleene1957}, но этот тезис (начало страницы 332) не использует термин ``partial recursive functional'' и вообще этот термин на с.~332 не встречается. Предметный указатель этой книги~\cite{Kleene1957} отсылает к с.~326 по слову ``partial recursive functional'', но и эта страница не содержит соответствующего упоминания. Правда, на этой странице даётся определение частичной рекурсивности частичной функции $\varphi$ относительно частичных функций $\psi_1,\ldots,\psi_k$, соответствующее сводимости графиков по перечислимости (первый вариант из приведённых выше трёх). и говорится о <<схеме>>  (scheme) $F$, но какие требования предъявляются к этой схеме (должна ли она давать функцию только в применении к функциям $\psi_1,\ldots,\psi_k$, или к любым $k$ функциям), из текста не ясно. (А нумерации класса вычислимых с данным оракулом функций рассматриваются только для случая, когда оракул представляет собой множество или всюду определённую функцию.) Но Майхилл и Шепердсон уточняют, что в их результате речь идёт от частично рекурсивных функционалах, определённых (и дающих функции) для всех функций в качестве аргументов, что эквивалентно второму определению, как и должно быть.

Одифредди~\cite[Definition II.3.6]{Odifreddi1989} определяет частично рекурсивные функционалы со ссылкой на Клини~\cite{Kleene1957}, но следует третьему варианту определения (композиция операций подстановки, минимизации и рекурсии, применённых к частично рекурсивным функциям и аргументам), который в~\cite{Kleene1957} не встречается. Понятие, соответствующее второму варианту определения, он называет ``effectively continuous functional''  или ``recursive operator'', а первому\emdash ``partial recursive operator''.

% Одифредди даёт определения в топологических терминах, ссылаясь на работы Успенского 1955 года и Нероуда 1957 года (последнюю я пока не видел)

Возвращаясь к работам Успенского~\cite{1955,1955a,1956}, можно отметить следующие достижения:
\begin{itemize}
\item впервые было дано определение сводимости по перечислимости;

\item впервые было опубликовано (восходящее к Колмогорову) определение нумерации и сводимости нумераций;

\item были проанализированы свойства нумераций вычислимых функций и перечислимых множеств, необходимые для рассуждений о номерах программ, введено понятие главной нумерации (впоследствии переоткрытое Роджерсом);

\item была доказана (до того, как это сделал сам Райс) гипотеза Райса об описании вполне перечислимых классов перечислимых множеств; получены также аналогичные результаты для вычислимых функций (вместо перечислимых множеств); в частности, доказана невозможность распознавания нетривиальных свойств вычислимых функций по их номерам в главной нумерации;

\item было дано определение вычислимой операции в терминах топологического подхода (вычислимость как некоторый специальный случай непрерывности) и доказано (одновременно с Майхиллом и Шепердсоном), что ограничения вычислимых операцией на вычислимые функции  можно эквивалентно описать как алгоритмические преобразования программ, а также доказана аналогичная теорема для операций над перечислимыми множествами.
\end{itemize}

Во избежание недоразумений отметим, что Успенский не рассматривает алгоритмы, определённые на всех программах \emph{всюду определённых} функций и дающие одинаковые значения на эквивалентных программах: эти работы Крайзеля, Лакомба, Шёнфилда~\cite{KreiselLacombeShoenfield1959} (позднее обобщённые Цейтиным~\cite{Tseitin1962} на конструктивные метрические пространства) никак не пересекаются с работами Успенского.

Приведём отрывок из воспоминаний Успенского~\cite[с.~905--907, 912]{2018b}, где он пишет о своих результатах 1955 года и об их представлении на Третьем всесоюзном математическом съезде (1956):

\begin{quote}
В сделанном 26 июня обзорном докладе <<Об алгоритмической сводимости>> я рассказывал о четырёх видах сводимости и связях между ними. Это \emph{сводимость по вычислимости}, состоящая в сведении вычисления одной функции к вычислению другой. Это \emph{сводимость по разрешимости}, состоящая в сведении задачи построения разрешающего алгоритма для одного множества к задаче построения разрешающего алгоритма для другого множества. Это  \emph{сводимость по перечислимости}, состоящая в сведении перечисления одного множества к перечислению другого. Это \emph{сводимость проблем}, состоящая в сведении решения одной проблемы к решению другой.  $\langle\ldots\rangle$  [Говоря о понятии сводимости проблем, Успенский указывает на его источник:] В 1955 году интересную разновидность проблем ввёл ученик Колмогорова Ю.\,Т.\,Медведев; он же определил понятие сводимости для таких проблем.  $\langle\ldots\rangle$ 

Наименование моего доклада 2 июля было <<Понятие программы и вычислимые операторы>>, а сообщения 3 июля\emdash <<Вычислимые операции, вычислимые операторы и конструктивно\dash непрерывные функции>>. Доклад и сообщение были тесно связаны тематически.

В сообщении 3 июля был изложен (разумеется, без доказательства) результат, который я считаю своим главным математическим результатом; я помню обстоятельства его получения\footnote{Теорема 3 интересовала меня с семиотической точки зрения, $\langle\ldots\rangle$ хотя слова `семиотика' я тогда, скорее всего, ещё не знал. Помню, как я бродил по улицам и думал только об этом. Озарение пришло, когда в течение некоторого времени я проводил дневные часы в квартире моей тёщи в Большом Спасоглинищевском переулке. Сын ещё не родился, жена и тёща уходили на работу, телефона в квартире не было, мобильные телефоны ещё не были изобретены. Вот в этой обстановке меня внезапно осенило.>> [Примечание В.\,А.\,Успенского]}\emdash теорема 3 (см. ниже). Он составил основу моей кандидатской диссертации, защищённой в октябре 1955 года. Доказательства этого результата, кроме текста диссертации, хранящейся\emdash или только уже хранившейся\emdash в библиотеке мехмата, я так никогда и не опубликовал. Основная причина, как ни стыдно в этом признаться, банальная лень. Дополнительная причина, не столь постыдная, но глупая, это преследовавшее [в тексте: преследующее] меня, пока я не повзрослел, желание изложить всё в максимально общем виде, но достичь предела в обобщении нереально. $\langle\ldots\rangle$

\textbf{Теорема 3}. \emph{Пусть функция $g$ с натуральными аргументами и значениями обладает следующим свойством. Если $m$ и $n$ служат программами одной и той же вычислимой функции от $s$ аргументов, то $g(m)$ и $g(n)$ также служат программами одной и той же функции от одного аргумента. Тогда существует вычислимый оператор $V$ со следующим свойством. Для всякой функции $\theta$ с программой $n$ значением $V(\theta)$ оператора $V$ на функции $\theta$ является функция с программой $g(n)$}.

\textbf{Философский комментарий}. Семиотический смысл теоремы~3 таков: <<хорошее>> вычислимое преобразование имён сопровождается вычислимым преобразованием соответствующих объектов. 
\end{quote}

\subsection*{Конструктивность в классической математике}

Идея о том, что можно понимать математические утверждения конструктивно, была известна давно (<<интуиционизм>> Брауэра и его последователей и позже <<конструктивизм>> Маркова и его учеников). В частности, утверждения вида <<для всех $x$ существует такой $y$, что\ldots>> при  конструктивном понимании предполагают, что существует некоторый способ получения этого самого <<существующего>> $y$ по любому данному значению~$x$. 

Обычно вместе с этим предлагали изменить и саму логику, понимая конструктивно, в частности, связку <<или>> и не пользуясь законом исключённого третьего. Другое (на первый взгляд, напрашивающееся) предложение, а именно, рассматривать <<эффективные>> аналоги классических определений и результатов внутри обычной (<<неконструктивной>>) математики, удивительным образом сначала было менее популярно, и часто считалось, что если уж мы рассуждаем об алгоритмах, то это почему-то обязывает нас рассуждать <<конструктивно>>,  <<финитно>> и т.п.  В отличие от этой традиции, Успенский систематически пропагандировал <<классический>> подход к алгоритмическим понятиям. Вот два примера из его работ.

Есть разные конструкции действительных чисел (сечения Дедекинда, фундаментальные последовательности, вложенные отрезки, бесконечные десятичные дроби). Для каждой из них можно рассмотреть её эффективный вариант. Скажем, для сечений Дедекинда можно требовать существования алгоритма, который по рациональному числу говорит, в каком из двух множеств сечения оно лежит. Для фундаментальной последовательности рациональных чисел естественно требовать существования алгоритма, вычисляющего её члены, а также <<регулятора сходимости>>\emdash алгоритма, находящего по (рациональному) $\eps>0$ то место, начиная с которого члены последовательности отличаются менее чем на~$\eps$. Для бесконечной десятичной дроби естественно требовать существования алгоритма, который по $n$ указывает $n$\dash ю цифру дроби, и так далее.

В каждом из этих случаев возникает некоторое подмножество множества действительных чисел (соответствующее тем числам, для которых имеются такие эффективные представления). Можно поставить, оставаясь в рамках <<классической математики>>, вопрос о том, дают ли перечисленные варианты определений одно и то же подмножество или разные. Нетрудно убедиться, что одно и то же, и его элементы называют \emph{вычислимыми действительными числами}.

Здесь хорошо видна разница с конструктивистским подходом (скажем, в смысле школы А.\,А.\,Маркова). Для конструктивистов никаких <<обычных действительных чисел>> не существует, и множество конструктивных действительных чисел не является подмножеством никакого большего множества. Вместо этого конструктивным действительным числом называется пара алгоритмов (один вычисляет члены последовательности, другой является регулятором сходимости), и рассуждать о таких парах предлагается в рамках конструктивной логики. При этом не все варианты определения действительных чисел равнозначны с точки зрения их конструктивизации. Скажем, подход с десятичными дробями неудачный\emdash потому, например, что для так определённых конструктивных действительных чисел нельзя конструктивно определить сложение (как преобразование, которое по алгоритмам для двух дробей давало бы алгоритм для их суммы).

Но этот же дефект, замечает Успенский, можно проанализировать и в рамках классической математики. Будем интересоваться не только тем, один и тот же класс действительных чисел возникает в рамках разных определений или разные, но также и более тонким вопросом: эквивалентны ли нумерации множества вычислимых действительных чисел, которые получаются из этих определений (можно ли по номеру вычислимого действительного числа в одной нумерации алгоритмически получить номер того же числа в другой). И тут возникает та же самая проблема с десятичными дробями (и с дробями в любых системах счисления). Подробно этот вопрос разобран в~\cite{1960}, где даны необходимые и достаточные условия, при которых переход от одного основания к другому (в позиционной записи) эффективен.

Другой пример, разобранный Успенским в~\cite{1960a}\emdash эффективизация определения бесконечного множества. Можно определить бесконечность множества $X$ так: для всякого $n$ в множестве $X$ есть не менее $n$ элементов. Или так: для всякого конечного множества $F$ есть число, которое отличает $F$ от $X$, то есть принадлежит симметрической разности $F\bigtriangleup X$. Оба эти определения можно эффективизировать, потребовав существования соответствующих алгоритмов. В первом случае речь идёт об алгоритме, который по $n$ указывает список из $n$ элементов множества $X$; во втором случае речь идёт об алгоритме, который применим к любому конечному множеству $F$ и даёт какой-то элемент разности $F\bigtriangleup X$. Можно проверить, что эти свойства (существование того и другого алгоритма) эквивалентны (и даже можно без нарушения эквивалентности требовать, чтобы по $F$ давался элемент разности $X\setminus F$). В терминологии Поста~\cite{Post1944} оба эти свойства равносильны тому, что множество $X$ не является иммунным (содержит бесконечное перечислимое подмножество).

Не все естественные определения бесконечности приводят к эквивалентным эффективизациям. Скажем, можно сказать, что множество $X$ бесконечно, если для всякого $n$ можно указать начальный отрезок $[0,N]$ натурального ряда, содержащий по крайней мере $n$ элементов множества $X$. Эффективный вариант этого определения требует, чтобы существовал алгоритм, указывающий $N$ по $n$. Это более слабое определение эффективной бесконечности, которое, как доказал Успенский в~\cite{1957a}, равносильно тому, что множество не является гипериммунным в смысле Поста~\cite{Post1944}. (Параллельно и независимо, отвечая на вопрос Колмогорова, это же доказали А.\,В.\,Кузнецов и Ю.\,Т.\,Медведев.)

Интересно отметить, хотя это и не имеет отношения к работам Успенского, что важнейшее достижение алгоритмической теории случайности, а именно, определение случайности по Мартин\dash Лёфу, данное им в 1966 году~\cite{MartinLof1966}, тоже можно рассматривать как естественную эффективизацию классического (во всех смыслах этого слова) определения нулевого множества (множества нулевой меры в смысле Лебега).  В этом классическом определении (скажем, для подмножеств отрезка) говорится, что множество $X\subset[0,1]$ является нулевым, если для всякого $\eps>0$ существует покрытие множества $X$ интервалами с суммой длин не больше $\eps$. По очевидным причинам можно ограничиться рациональными значениями $\eps$ и интервалами с рациональными концами. Тогда они будут конструктивными объектами и можно эффективизировать определение, потребовав, чтобы существовал алгоритм, который, получив на вход $\eps$, перечисляет интервалы покрытия с требуемыми свойствами. Это и предложил Мартин\dash Лёф.

Можно добавить, что многие вопросы и результаты алгоритмической теории случайности можно интерпретировать как вопросы об эффективизации классических понятий и теорем. Скажем, критерий Соловея случайности по Мартин\dash Лёфу, как заметил Александр Буфетов, является эффективным вариантом классической леммы Бореля\endash Кантелли. При этом интересно, что стандартное её доказательство (хвосты сходящегося ряда могут быть сколь угодно малы) не эффективизируется, и приходится использовать другое (тоже естественное и несложное, см. подробности в~\cite{2013}). Другой поучительный пример того же рода\emdash обнаруженное В.\,В.\,Вьюгиным (учеником Успенского) доказательство эффективного варианта эргодической теоремы~\cite{Vyugin1998}.

\subsection*{Алгоритмическая теория информации}

Странным образом Успенский, будучи учеником Колмогорова и работая рядом с ним на мехмате МГУ, был в 1960-е годы в стороне от исследований Колмогорова, связанных с определением понятия сложности конечного объекта. По его словам, он впервые вплотную занялся этой областью, готовя доклад (с А.\,Л.\,Семёновым) на конференции в Ургенче~\cite{1981,1982a}. В этом докладе была предложен подход к классификации различных видов сложности (или, как предпочитал говорить Успенский, <<энтропии>>) для конечных объектов, определённых к тому времени (простая, префиксная, условная, монотонная энтропии, а также энтропия разрешения). Изначально (см.~\cite{Shen1984}) этот подход был предложен в терминах $f_0$-пространств и операций над ними, что можно рассматривать как развитие  подхода Успенского к определению вычислимых отображений как частного случая непрерывных. Однако для целей классификации различных видов колмогоровской сложности можно было обойтись и без $f_0$\dash пространств, и Успенский с Семёновым в~\cite{1981,1982a} предложили более простой вариант определения (в терминах отношения согласованности на описаниях и объектах), достаточный для целей классификации. Впоследствии это упрощённое описание было изложено в~\cite{1992a,1996}; подробное изложение с топологической точки зрения (но без $f_0$\dash пространств) можно найти в~\cite{2013}.

Алгоритмической теории информации (точнее, различным определениям понятия случайности) посвящён также обзор~\cite{1990} и монография~\cite{2013}. Популярная лекция (для студентов\dash младшекурсников), посвящённая различным определениям случайности, была прочитана Успенским на летней школе <<Современная математика>> в Дубне в 2005 году, и материалы этой лекции были изданы~\cite{2006} и вошли в качестве приложения в монографию~\cite{2013}.

До сих пор остаётся открытым вопрос, поставленный Успенским, Семёновым и Мучником в~\cite{1998} о том, совпадают ли понятия случайности по Мартин\dash Лёфу и <<непредсказуемости>> (отсутствия выигрышной вычислимой стратегии в немонотонных играх,  см. подробнее в~\cite{2006,2013}).

\subsection*{Популяризация}

Есть разные представления о том, что такое <<популяризация науки>> (по-французски это называют vulgarisation, что для русского уха звучит обидно, хотя  и не совсем незаслуженно). Можно рассказывать байки о трудной судьбе или особенностях личной жизни выдающихся учёных. Можно пересказывать недавно прочитанное в другой популярной книге, добавляя <<оживляж>>. Видимо, это полезное дело\emdash но Успенский всю жизнь занимался другим, пытаясь честно объяснить <<суть дела>>. При этом сложность темы, естественно, зависела от того, на кого рассчитано объяснение, но всегда это было настоящее объяснение того, что объяснить можно, с отчётливым указанием, что именно оставлено без доказательства или уточнения. При этом он не боялся объяснить очевидное, помня, что известное известно немногим (<<Что мы знаем о лисе? Ничего, и то не все>>\emdash писал сосед Колмогорова по даче, знаменитый детский писатель Борис Заходер.)

Ещё будучи студентом, Успенский (вместе со старшим соавтором, Евгением Борисовичем Дынкиным) написал книгу <<Математические беседы>>~\cite{1952a} по материалам математических кружков, где он сначала был участником, а потом руководителем. В ней несколько тем (раскраска графов, начала теории чисел и теории вероятностей) представлены в виде последовательности задач, как это и делалось на кружке,  и приведены решения этих задач. И до этого издавались книжки с задачами математических кружков, но тут идея была в том, что эти задачи в целом образуют изложение некоторой математической теории. Книга Дынкина и Успенского была библиографической редкостью, пока не была переиздана уже сравнительно недавно,  в 2004~году (и не появилась в интернете).

Успенский читал лекции для школьников (в частности, для участников математических олимпиад) и написал несколько популярных брошюр в серии <<Популярные лекции по математике>>, никак не связанных с его собственными научными интересами  (математической логикой и теорией алгоритмов):  про применение механики в математике~\cite{1958a} и  про треугольник Паскаля~\cite{1966}. Впрочем, последняя брошюра затрагивает и логический вопрос: что означает решить комбинаторную задачу и почему нужно фиксировать список разрешённых операций (скажем, включив в него факториалы, но исключив обозначения для биномиальных коэффициентов).

Две другие брошюры в этой серии (<<Машина Поста>>~\cite{1979} и <<Теорема Гёделя о неполноте>>~\cite{1982}) уже посвящены темам из математической логики и теории алгоритмов. Первая из них совсем элементарна и основана на занятиях с младшеклассниками, вторая, наоборот, написана на основе статьи в <<Успехах математических наук>>~\cite{1974} и предполагает некоторую математическую культуру, но вполне может быть прочитана и понята продвинутыми старшеклассниками. Популярному изложению ещё одной темы, так называемого <<нестандартного анализа>>, в котором методы математической логики используются для математически корректного рассмотрения бесконечно малых и бесконечно больших величин, посвящена брошюра~\cite{1983}; её расширенный вариант был опубликован затем издательством <<Наука>>~\cite{1987}.

Как герой <<Игры в бисер>> Гессе, становясь старше, Успенский объяснял всё более и более базовые вещи, занявшись проповедью математики среди <<гуманитариев>>. Впрочем, началось это уже давно, в 1960-е годы, когда он разрабатывал и осуществлял курс математики для отделения теоретической и прикладной лингвистики филологического факультета МГУ, но в последние два десятилетия своей жизни он обращался к гораздо более широкой аудитории.  Несколько его лекций на летних школах по математике и лингвистике в Дубне, к счастью, сохранились как видеозаписи (прежде всего благодаря Виталию Арнольду), и по ним (см. ссылки в \url{http://www.mathnet.ru/php/person.phtml?option_lang=rus&personid=20219}) можно составить представление об Успенском как лекторе\emdash хотя, конечно, в полной мере оценить его  можно было только на университетских лекциях, особенно спецкурсах. При этом проповедь Успенского была именно проповедью математики, а не <<о математике>>. Он рассказывал простые вещи, но всерьёз, с определениями, примерами и доказательствами. Одна из его последних книжек~\cite{2009b} так и называется:   <<Простейшие примеры математических доказательств>>. Другая книжка~\cite{2000} называется <<Что такое аксиоматический метод?>>\emdash и там тоже подробно разобрано множество примеров (в частности, из школьной геометрии, точнее, из той части школьной геометрии, которая в школах пропускается). Например, объясняется, как вывести из аксиом, что для всякой прямой найдётся точка, на ней не лежащая. Материалы из этих двух книг вошли в сборник <<Апология математики>>~\cite{2009} (вместе с другими статьями, уже более общего характера). И удивительным образом проповедь Успенского имела успех (по крайней мере в том же смысле, что у Антония Падуанского): Успенскому была присуждена премия <<Просветитель>> (учреждённая Дмитрием Борисовичем Зиминым и фондом <<Династия>>) за 2010 год в области естественных и точных наук. 

Помимо собственных книг, большой заслугой Успенского (то, что называют по\dash английски community service) является организация издания многих классических книг по математической логике и теории алгоритмов на русском языке: он переводил (по инициативе Колмогорова) книгу Р.\,Петер о рекурсивных функциях~\cite{Peter1954}, редактировал перевод классических монографий Клини~\cite{Kleene1957} и Роджерса~\cite{Rogers1972}, книги Дэвиса о нестандартном анализе~\cite{DavisNonStandard1980} (видимо, первого изложения идей нестандартного анализа, появившегося по-русски), фундаментального учебника Чёрча по логике~\cite{Church1960}, первого тома <<Начал математики>> Бурбаки~\cite{BourbakiSetTheory1965}, а также книги Эшби о кибернетике~\cite{Ashby1959}.